\documentclass[11pt]{article}

\usepackage[T1]{fontenc}
\usepackage{lmodern}
\usepackage[utf8]{inputenc}
\usepackage{amsmath,amssymb,amsthm,mathtools}
\usepackage{mathrsfs}
\usepackage{enumitem}
\usepackage[dvipsnames]{xcolor}
\usepackage[a4paper,margin=3.2cm]{geometry}
\usepackage{microtype}
\usepackage[hidelinks,pagebackref,colorlinks=true,linkcolor=blue,citecolor=blue]{hyperref}

\usepackage[nottoc,notlot,notlof]{tocbibind}

\mathtoolsset{showonlyrefs}
\theoremstyle{plain}
\newtheorem{theorem}{Theorem}[section]
\newtheorem{proposition}[theorem]{Proposition}
\newtheorem{lemma}[theorem]{Lemma}
\newtheorem{corollary}[theorem]{Corollary}
\newtheorem{mainthm}{Theorem}

\theoremstyle{definition}
\newtheorem{definition}[theorem]{Definition}
\newtheorem{example}[theorem]{Example}

\theoremstyle{remark}
\newtheorem{remark}[theorem]{Remark}

\newcommand{\WM}[0]{\hyperlink{WM}{{\rm WM}_{p,\alpha}}}

\newcommand{\R}{\mathbb R}
\newcommand{\1}{\mathbf 1}
\newcommand{\cE}{\mathcal E}
\newcommand{\cB}{\mathcal B}
\newcommand{\cF}{\mathcal F}

\newcommand{\cN}{\mathcal N}
\newcommand{\cR}{\mathcal R}
\newcommand{\Lip}{\operatorname{Lip}}
\newcommand{\Per}{\operatorname{Per}_G}
\newcommand{\Vol}{\mathrm V}

\title{Weak monotonicity, Critical exponents, \\ and Gagliardo-Nirenberg inequalities on weighted graphs}
\author{Patricia Alonso Ruiz\thanks{Institute of Mathematics, Friedrich Schiller University Jena,
07737 Jena, Germany. \href{mailto:patricia.alonso.ruiz@uni-jena.de}{\texttt{patricia.alonso.ruiz@uni-jena.de}}}
\and Fabrice Baudoin\thanks{Department of Mathematics, Aarhus University,
Ny Munkegade 118, 8000 Aarhus C, Denmark. \href{mailto:fbaudoin@math.au.dk}{\texttt{fbaudoin@math.au.dk}}}
\and Li Chen\thanks{Department of Mathematics, Aarhus University,
Ny Munkegade 118, 8000 Aarhus C, Denmark. \href{mailto:lchen@math.au.dk}{\texttt{lchen@math.au.dk}}}}
\date{}

\begin{document}

\maketitle

\begin{abstract}
We study weak monotonicity of heat kernel and ball averaged functionals on weighted graphs, and  the associated critical exponents. Under sub-Gaussian heat kernel estimates, weak Bakry–Émery estimates for the semigroup are proved, with explicit values of the H\"older regularity exponent in the strongly recurrent case.  Applications of the theory include Gagliardo–Nirenberg, Sobolev, Nash and isoperimetric inequalities on strongly recurrent graphs and their products.
\end{abstract}

\tableofcontents
\section{Introduction}\label{sec:intro}

We study weak monotonicity of heat kernel and ball averaged functionals on weighted graphs, related critical exponents, and the resulting Gagliardo--Nirenberg and isoperimetric inequalities. Our guiding principle   goes back to Varopoulos~\cite{Varopoulos1985}, Coulhon and Saloff-Coste~\cite{CoulhonSaloffCoste1993}, and Coulhon~\cite{Coulhon1995,Coulhon1996}: 
$\{P_t\}_{t\geq 0}$ implies Sobolev type inequalities as soon as a
pseudo-Poincar\'e inequality like 
\begin{equation}\label{E:pPI_intro}
\|f-P_tf\|_p\le C\sqrt t\,\|\nabla f\|_p    
\end{equation}
is available. Then, all Gagliardo-Nirenberg inequalities become ``Sobolev inequalities in disguise'' as described by Bakry, Coulhon, Ledoux and Saloff-Coste~\cite{BCLS1995}, who developed this principle into a general method for deriving Gagliardo--Nirenberg inequalities; see also~\cite{SaloffCoste2002}, and~\cite{Ledoux1996} for its isoperimetric form.

\medskip

On a Dirichlet space $X$ with heat kernel $p_t$ and reference measure $\mu$, where a gradient need not be available, the heat kernel Besov seminorms
\begin{equation}\label{E:def_Besov_seminorm_intro}
    \|f\|_{p,\alpha}:=\sup_{t>0}\Big(t^{-p\alpha}\int_X\int_X|f(x)-f(y)|^p\,p_t(x,y)\,d\mu(y)\,d\mu(x)\Big)^{1/p}
\end{equation}
provide a substitute for the Sobolev norms appearing in~\eqref{E:pPI_intro}, see~\cite{AlonsoEtAl2020I,AlonsoEtAl2020II,AlonsoEtAl2021III}. Jensen's inequality gives the pseudo-Poincar\'e estimate 
\[
\|f-P_tf\|_p\le t^{\alpha}\|f\|_{p,\alpha},
\]
which, together with ultracontractivity, allows to apply the methods~\cite{BCLS1995} as long as the heat kernel Besov functionals 
\begin{equation*}
    \cN_{p,t}(f):=\Big(\int_X\int_X|f(x)-f(y)|^p\,p_t(x,y)\,d\mu(y)\,d\mu(x)\Big)^{1/p}
\end{equation*}
satisfy the following property: For every $t>0$ and every $f$ in the Besov class,
\begin{equation}\label{eq:intro-WM-continuous}
 t^{-\alpha}\cN_{p,t}(f)\le C\liminf_{s\downarrow0}s^{-\alpha}\cN_{p,s}(f)\,;
\end{equation}
see~\cite[Definition~6.7 and Theorem~6.9]{AlonsoEtAl2020I} and~\cite{AlonsoRuizBaudoin2021}. In the literature, this property has sometimes been referred to as weak monotonicity, see e.g.~\cite{MuruganShimizu2023,KajinoShimizu2026,GaoYuZhang2026}. We will however call it $\mathrm{P}(p,\alpha)$, as originally done in~\cite{Baudoin2024KS}, since we will introduce a stronger condition, which we shall consider ``actual'' weak monotonicity and compares two different scales, see Definition~\ref{def:WM}. 

\medskip

By identifying the resulting classes with Sobolev or BV spaces, one recovers the classical Sobolev and isoperimetric inequalities on complete Riemannian manifolds with nonnegative Ricci curvature~\cite{AlonsoEtAl2020II}, and also obtains Gagliardo--Nirenberg inequalities in the generality of~\cite{AlonsoRuizBaudoin2021}. The method already yielded new Sobolev and isoperimetric inequalities on fractals in~\cite{AlonsoEtAl2021III}.

\medskip

For $p=2$, the property~\eqref{eq:intro-WM-continuous} follows from spectral theory. For $p\ne2$, it is a genuine additional property of the heat semigroup. In the sub-Gaussian setting of~\cite{AlonsoEtAl2021III}, with volume and walk dimensions $d_f$ and $d_w$, respectively, the $p=1$ estimate in \eqref{eq:intro-WM-continuous} follows from the weak Bakry-\'Emery estimate
\[
 |P_tf(x)-P_tf(y)|\le C\left(\frac{d(x,y)}{t^{1/d_w}}\right)^{\kappa}\|f\|_\infty.
\]
On strongly recurrent spaces, this bound holds with $\kappa=d_w-d_f$, as a consequence of the resolvent estimates of Barlow~\cite{Barlow1998}; see~\cite[Theorem~3.7]{AlonsoEtAl2021III}. Weak monotonicity also plays a central role in the study of $p$-energies and Korevaar--Schoen spaces on fractals and metric measure spaces~\cite{Baudoin2024KS,MuruganShimizu2023,KajinoShimizu2026,GaoYuZhang2026}. The goal of the paper is to develop the theory in the discrete setting.

\subsection*{The discrete setting}

Let $G=(V,E)$ be a connected, locally finite graph with symmetric conductances $c_{xy}$, vertex measure $m(x)=\sum_y c_{xy}$, and graph distance $d$. Write $P$ for the Markov operator with transition probabilities $\pi(x,y)=c_{xy}/m(x)$, and $p_n(x,y)$ for the density of its $n$-step transition probabilities with respect to $m$; see Section~\ref{sec:graphs}. For $1\le p<\infty$, $n\ge1$ and $f\in\ell^p(V,m)$, set
\begin{equation}\label{eq:intro-N}
 \cN_{p,n}^{\rm H}(f)
 :=\left(\sum_{x,y\in V}|f(x)-f(y)|^p
 p_n(x,y)m(x)m(y)\right)^{1/p},
\end{equation}
and $\|f\|_{p,\alpha}:=\sup_{n\ge1}n^{-\alpha}\cN_{p,n}^{\rm H}(f)$ for $\alpha\ge0$. At time one, $p_1(x,y)m(x)m(y)=c_{xy}$, so
\begin{equation}\label{eq:intro-gradient}
 \cN_{p,1}^{\rm H}(f)=\|\nabla_pf\|_p
 :=\left(\sum_{x,y\in V}c_{xy}|f(x)-f(y)|^p\right)^{1/p}.
\end{equation}

The first we will ask is for which $\alpha\ge0$ there is a constant $C$, independent of $f$, with
\begin{equation}\label{eq:intro-equivalence}
 \|\nabla_pf\|_p\le\|f\|_{p,\alpha}\le C\|\nabla_pf\|_p,
 \qquad f\in\ell^p(V,m).
\end{equation}
The lower bound follows from the term $n=1$. The upper bound is another formulation of what we call property $\mathrm{P}(p,\alpha)$ for the heat kernel functionals in Definition~\ref{def:P}; its optimal exponent is studied in Section~\ref{sec:critical-exponents} and plays the same role in the discrete theory as the critical exponents introduced in~\cite[Section~5.2]{AlonsoEtAl2020I}. We shall also consider what we call the stronger weak monotonicity property $({\rm WM}_{p,\alpha})$, that is
\begin{equation}\label{eq:intro-WM}
 r^{-\alpha}\cN_{p,r}^{\rm H}(f)
 \le C\,s^{-\alpha}\cN_{p,s}^{\rm H}(f)
 \qquad\text{for all integers }r\ge s\ge1
\end{equation}
and every $f\in\ell^p(V,m)$, with $C$ independent of $r,s,f$, c.f. Definition~\ref{def:WM}. This property persists when $\alpha$ increases. Taking $s=1$ gives $\mathrm{P}(p,\alpha)$, whereas the converse requires an additional argument. By Jensen's inequality, the property $\mathrm{P}(p,\alpha)$ alone yields the pseudo-Poincar\'e bound $\|f-P^nf\|_p\le Cn^\alpha\|\nabla_pf\|_p$ needed to apply the method of~\cite{BCLS1995}.

\vspace{0.6cm}

\subsection*{Assumptions, Main results}
Compared to the framework of continuous time and strongly local Dirichlet forms, discrete time introduces a specific issue: the spectrum of $P$ may contain $-1$, preventing analyticity of $(P^n)_{n\ge1}$. For heat kernel comparisons at arbitrary times, we assume the \emph{diagonal lower bound}~\eqref{eq:LB}, $\pi(x,x)\ge\varepsilon_{\rm L}>0$, which removes this obstruction; this is a classical assumption, see~\cite{CoulhonSaloffCoste1990,Dungey2006,Feneuil2018}. Without such a condition, periodicity can destroy weak monotonicity: for the walk exchanging two vertices, $\cN_{p,2}^{\rm H}(f)=0$ but $\cN_{p,3}^{\rm H}(f)=\cN_{p,1}^{\rm H}(f)>0$ for nonconstant $f$. A graph satisfying~\eqref{eq:LB} has \emph{controlled weights} if also $\pi(x,y)\ge\varepsilon_0>0$ for all neighbours $x\ne y$ (Definition~\ref{def:weight-assumptions}). We write $({\rm SG}_{d_f,d_w})$ for two-sided sub-Gaussian heat kernel estimates with volume exponent $d_f$ and walk exponent $d_w$ (Definition~\ref{def:sub-gaussian}). Set $\Vol(x,r):=m(B(x,r))$. Under~\eqref{eq:LB}, these imply the volume growth condition~\eqref{eq:volume-growth} (Proposition~\ref{prop:volume-growth}). Our main results are the following.

\begin{mainthm}\label{thm:main-p-ge-2}
Let $1\le p<\infty$. Under~\eqref{eq:LB}, $(\cN_{p,n}^{\rm H})_{n\ge1}$ satisfies $({\rm WM}_{p,1/\min\{p,2\}})$, with a constant depending only on $p$ and $\varepsilon_{\rm L}$. Without assuming~\eqref{eq:LB}, the heat kernel functionals satisfy $\mathrm{P}(p,\alpha)$ for every $\alpha\ge1/\min\{p,2\}$; equivalently, for every $f\in\ell^p(V,m)$,
\begin{equation}\label{eq:intro-p-ge-2-equivalence}
 \|\nabla_pf\|_p\le\|f\|_{p,\alpha}\le C_p\|\nabla_pf\|_p\,.
\end{equation}
For $1\le p\le2$, one has $\|f\|_{p,\alpha}=\|\nabla_pf\|_p$ for every $\alpha\ge1/p$; in particular, $\|f\|_{2,\alpha}^2=2\cE(f,f)$ for $\alpha\ge1/2$.
\end{mainthm}

For $1\le p<2$, semigroup regularity gives a smaller weak monotonicity exponent. Following~\cite{AlonsoEtAl2021III}, for $0<\kappa<d_w$ the weak Bakry--\'Emery condition $({\rm wBE}(\kappa))$ is defined by
\begin{equation}\label{eq:intro-wBE}
 |P^nf(x)-P^nf(y)|
 \le C\min\left\{1,
 \left(\frac{d(x,y)}{n^{1/d_w}}\right)^\kappa\right\}\|f\|_\infty,
 \qquad n\ge1,\quad x,y\in V,
\end{equation}
for every bounded $f$ (Definition~\ref{def:wBE}). For $1\le p<\infty$, set
\begin{equation}\label{eq:intro-beta}
 \beta_p(\kappa):=\frac1p+\left(1-\frac2p\right)\frac\kappa{d_w}\,.
\end{equation}

\begin{mainthm}\label{thm:main-wBE}
Assume $({\rm wBE}(\kappa))$ for some $d_w\ge2$ and $0<\kappa<d_w$. Then, for every $1\le p\le2$, $(\cN_{p,n}^{\rm H})_{n\ge1}$ satisfies $\mathrm{P}(p,\beta_p(\kappa))$. If, in addition,~\eqref{eq:LB} and $({\rm SG}_{d_f,d_w})$ hold, then the family satisfies $({\rm WM}_{p,\beta_p(\kappa)})$.
\end{mainthm}

The exponent $\beta_p(\kappa)$ interpolates between $1-\kappa/d_w$ at $p=1$ and $1/2$ at $p=2$, and is the one obtained by the continuous-time argument of~\cite[Section~3.4]{AlonsoEtAl2021III}. For property $\mathrm{P}(p,\beta_p(\kappa))$, it suffices to assume the regularity estimate only across adjacent vertices; see Theorem~\ref{thm:HK-P-wBE}. Under controlled weights and $({\rm SG}_{d_f,d_w})$, condition $({\rm wBE}(\kappa))$ holds for some $\kappa\in(0,1]$ by parabolic Harnack regularity (Proposition~\ref{prop:wBE-sub-gaussian}). 

Strong recurrence supplies an explicit value of $\kappa$.

\begin{mainthm}\label{thm:main-strongly-recurrent}
Assume~\eqref{eq:LB} and $({\rm SG}_{d_f,d_w})$ with $1\le d_f<d_w$. Then $({\rm wBE}(d_w-d_f))$ holds. Consequently, for every $1\le p\le2$, $(\cN_{p,n}^{\rm H})_{n\ge1}$ satisfies $({\rm WM}_{p,\beta_p(d_w-d_f)})$ with
\begin{equation}\label{eq:intro-WM-strong}
 \beta_p(d_w-d_f)=\frac1p+\left(1-\frac2p\right)\frac{d_w-d_f}{d_w}\,,
\end{equation}
so that $\beta_1(d_w-d_f)=d_f/d_w$ and $\beta_2(d_w-d_f)=1/2$.
\end{mainthm}

The proof combines a local Poisson estimate obtained from resistance bounds with the temporal estimate $\|LP^n\|_{\ell^\infty\to\ell^\infty}\le Cn^{-1}$. This gives the discrete counterpart of~\cite[Theorem~3.7]{AlonsoEtAl2021III}.

We also consider the ball averaged functionals
\begin{equation}\label{eq:intro-B}
 \cB_{p,r}(f)
 :=\left(\sum_{x\in V}\frac{m(x)}{\Vol(x,r)}
 \sum_{y\in B(x,r)}|f(x)-f(y)|^p m(y)\right)^{1/p},
 \qquad r\ge1,
\end{equation}
the discrete counterparts of the Korevaar--Schoen functionals in~\cite{Baudoin2024KS}.

\begin{mainthm}\label{thm:main-KS}
Assume $\inf_{x\in V}m(x)>0$ and~\eqref{eq:volume-growth} for some $d_f>0$. Then, for every $1\le p<\infty$ and $f\in\ell^p(V,m)$,
\begin{equation}\label{eq:intro-BL-WM}
 r^{-\sigma_p}\cB_{p,r}(f)
 \le C s^{-\sigma_p}\cB_{p,s}(f),
 \qquad r\ge s\ge1,
 \qquad
 \sigma_p:=1+\frac{d_f-1}{p}\,.
\end{equation}
\end{mainthm}

The proof uses only the graph metric and vertex measure, combining Kanai's discretization method~\cite{Kanai1985} with the path counting argument of~\cite{Chen2019}. Conditions~\eqref{eq:p0} and~\eqref{eq:volume-growth} imply the required lower bound on vertex masses. Under~\eqref{eq:LB} and $({\rm SG}_{d_f,d_w})$, Theorem~\ref{thm:ball-heat-comparison} gives
\begin{equation}\label{eq:intro-ball-heat}
 \cN_{p,n}^{\rm H}(f)\simeq\cB_{p,n^{1/d_w}}(f),
 \qquad n\ge1,\quad 1\le p<\infty\,,
\end{equation}
and hence heat kernel weak monotonicity with exponent $\sigma_p/d_w$ for every finite $p\ge1$. On strongly recurrent graphs, the weak Bakry--\'Emery estimate improves this exponent when $1<p\le2$ and $d_w-d_f<1$; see Remark~\ref{R:comparison_ball_vs_path}.

The infima of the admissible exponents for domination by $\|\nabla_pf\|_p$ are denoted by $\alpha_p^{\rm H}$ for the heat kernel functionals and $\alpha_p^{\rm B}$ for the ball averaged functionals. Under~\eqref{eq:LB} and $({\rm SG}_{d_f,d_w})$, the comparison above gives
\begin{equation}\label{eq:intro-critical-comparison}
 \alpha_p^{\rm B}=d_w\alpha_p^{\rm H}.
\end{equation}
We prove general bounds and stability under tensor powers and rough isometries. Weak Bakry--\'Emery regularity also gives $\alpha_p^{\rm H}\ge\beta_p(\kappa)$ for $p\ge2$. Consequently, whenever $d_w=d_f+1$, both critical exponents are attained and explicit for $p\ge2$. On planar Vicsek graphs, the same formulas hold for every finite $p\ge1$:
\begin{equation}\label{eq:intro-critical-Vicsek}
 \alpha_p^{\rm H}=\frac{d_f+p-1}{pd_w},
 \qquad
 \alpha_p^{\rm B}=\frac{d_f+p-1}{p},
 \qquad 1\le p<\infty.
\end{equation}
On the gasket, we identify the values at $p=1,2$ and obtain bounds for the other exponents. On infinite Cayley graphs of polynomial growth, equipped with counting measure and the walk $P=(I+P_0)/2$ where $P_0$ is simple random walk, we obtain $\alpha_p^{\rm H}=1/2$ and $\alpha_p^{\rm B}=1$ for every $1\le p<\infty$, with attainment and critical weak monotonicity; see Section~\ref{sec:examples-cayley}. In general, admissibility at the critical infimum and two-scale weak monotonicity are separate questions.

\subsection*{Functional inequalities and applications}

An abstract theorem of Section~\ref{sec:prelim} transfers the previous estimates to Gagliardo-Nirenberg inequalities by the method of~\cite{BCLS1995}. In the heat kernel setting, assume~\eqref{eq:LB} and $({\rm SG}_{d_f,d_w})$, and put $\delta=d_f/d_w$. If the heat kernel functionals satisfy $\mathrm{P}(p,\alpha)$ for some $\alpha>0$, then
\begin{equation}\label{eq:intro-GN-gradient}
 \|f\|_a\le C\|f\|_b^{1-\theta}\|\nabla_pf\|_p^\theta,
 \qquad
 \frac1a=\theta\left(\frac1p-\frac\alpha\delta\right)+\frac{1-\theta}{b},
\end{equation}
for $f\in\ell^p(V,m)$, $a,b\in(0,\infty]$, and $\theta\in(0,1]$, with $a<\infty$ if $p\alpha=\delta$ and $p>1$. When $\theta=1$, the factor $\|f\|_b^{1-\theta}$ is omitted. Thus $\mathrm{P}(p,\alpha)$ suffices for the graph inequalities. Under $({\rm WM}_{p,\alpha})$, the abstract theorem permits replacing the scale-one functional by $\inf_{s\ge1}\cN_{p,s}(f)$.

The same theorem applies to ball kernels, with decay exponent $\delta=d_f$. Its consequences include Sobolev, Nash, and $L^\infty$ interpolation inequalities. When $p\alpha=\delta$ and $p>1$, the truncation argument of~\cite[Corollary~6.4]{BCLS1995} also gives exponential integrability for functions of finite support and hence, since $\inf_xm(x)>0$, a logarithmic $L^\infty$ estimate (Corollary~\ref{cor:abstract-Trudinger}). At $p=1$, coarea relates the Sobolev estimate to an isoperimetric inequality; the kernel perimeters introduced in Section~\ref{sec:prelim} are comparable to the weighted edge perimeter under the respective properties $\mathrm{P}(1,\alpha)$, with the ball comparison taken under~\eqref{eq:p0} and~\eqref{eq:volume-growth}, or under~\eqref{eq:LB} and sub-Gaussian estimates. Smaller admissible exponents yield stronger global discrete Sobolev inequalities, since $\ell^q(V,m)\subset\ell^{q'}(V,m)$ for $q\le q'$ when $m$ is bounded below.

\paragraph{Cayley graphs.}
Let $\Gamma$ be an infinite group generated by a finite symmetric set $S$ and satisfying $(\mathrm V_D)$, that is, $|B(e,r)|\simeq r^D$ for $r\ge1$; by Gromov's theorem~\cite{Gromov1981} and the Bass-Guivarc'h formula~\cite{Bass1972,Guivarch1973}, $D$ is an integer.  Theorem ~\eqref{thm:main-wBE}  recovers the classical Sobolev, Nash and isoperimetric inequalities of Varopoulos and of Coulhon and Saloff-Coste on groups of polynomial growth~\cite{Varopoulos1985,CoulhonSaloffCoste1993,VaropoulosSaloffCosteCoulhon1992}; see Section~\ref{sec:examples-cayley}.

\paragraph{Strongly recurrent graphs.}
Theorem~\ref{thm:main-strongly-recurrent} gives a critical $L^1$ Sobolev estimate, $L^\infty$ interpolation inequalities, and Nash inequalities; see Corollary~\ref{cor:strongly-recurrent-inequalities}. The supported estimates improve the bounds obtained from~\eqref{eq:volume-growth} alone when $1<p\le2$ and $d_w-d_f<1$.

\paragraph{Products.}
Under controlled weights, $({\rm SG}_{d_f,d_w})$, and $({\rm wBE}(\kappa))$, products preserve the walk and weak Bakry--\'Emery exponents while replacing $d_f$ by $Nd_f$ (Proposition~\ref{prop:tensor-product}). For strongly recurrent factors, Corollary~\ref{cor:strongly-recurrent-products} gives Sobolev and Nash inequalities, including
\begin{equation}\label{eq:intro-product-isoperimetry}
 m_N(A)^{(N-1)/N}\le C_N\operatorname{Per}_{G^{\otimes N}}(A),
 \qquad N\ge2,\quad 0<m_N(A)<\infty\,,
\end{equation}
where $m_N=m^{\otimes N}$; compare~\cite[Proposition~3.8 and Remark~4.19]{AlonsoEtAl2021III}.

\subsection*{Organization of the paper}

In Section~\ref{sec:prelim}, the abstract Gagliardo--Nirenberg principle and its isoperimetric form are proved. Section~\ref{sec:graphs} outlines the assumptions on the graphs and their associated heat kernel, while Section~\ref{sec:HK-Besov} establishes the weak monotonicity for heat kernel functionals. In Section~\ref{sec:wBE}, weak Bakry--\'Emery estimates, including the explicit strongly recurrent exponent, are proved and, Section~\ref{sec:BL} treats ball averaged functionals, comparing them with heat kernel functionals. The critical exponents for comparison with the edge gradient are treated in Section~\ref{sec:critical-exponents}, and Section~\ref{sec:examples} offers several examples including Cayley graphs of polynomial growth, strongly recurrent graphs and fractals. Functional inequalities on these graphs and their products are also derived.

\medskip

Throughout, the letters $C,c$ will denote positive constants depending only on parameters and other constants involved in the hypotheses; their values may change between occurrences. We write $A\simeq B$ when $cB\le A\le CB$.

\section{Gagliardo--Nirenberg inequalities}\label{sec:prelim}

In this section, we prove Gagliardo--Nirenberg inequalities under ultracontractivity and the property $\mathrm{P}(p,\alpha)$, which compares every scale with scale one. Weak monotonicity gives the stronger conclusion with the infimum over all reference scales. The proof adapts the method of~\cite{BCLS1995}, as used in~\cite{AlonsoRuizBaudoin2021,Baudoin2024KS}, to discrete time.

\subsection{Property \texorpdfstring{$\mathrm{P}(p,\alpha)$}{P(p,alpha)} and weak monotonicity}\label{subsec:WM}
Let $V$ be a countable set equipped with a measure $m$ such that $0<m(x)<\infty$ for every $x\in V$. For $A\subset V$, we have  $m(A)=\sum_{x\in A}m(x)$.

\medskip

In addition, functions will be assumed to take real values and, unless otherwise explicitly mentioned, all function spaces $\ell^p(V)$ will be taken with respect to the measure $m$. For $0<p<\infty$, we will use the notation
\[
 \|f\|_p:=\left(\sum_{x\in V}|f(x)|^pm(x)\right)^{1/p}.
\]

\medskip

Let $(k_n)_{n\ge1}$ be a sequence of nonnegative kernels  $V\times V \to \mathbb{R}$ with the property that
\begin{equation}\label{eq:averaging-normalization}
 \sum_{y\in V}k_n(x,y)m(y)=1,
 \qquad x\in V,\quad n\ge1.
\end{equation}
The corresponding averaging operators are defined as
\begin{equation}\label{eq:averaging-operator}
 K_nf(x):=\sum_{y\in V}k_n(x,y)f(y)m(y)
\end{equation}
whenever the series is absolutely convergent. We do not assume that the kernels $k_n$ are symmetric in $x,y$.

\begin{definition}\label{def:averaging-functional}
Let $1\le p<\infty$. For $n\ge1$ and $f\in\ell^p(V)$, define, with values in $[0,\infty]$,
\begin{equation}\label{eq:Npn-general}
 \cN_{p,n}(f)
 :=\left(\sum_{x,y\in V}|f(x)-f(y)|^p
 k_n(x,y)m(x)m(y)\right)^{1/p}
\end{equation}
and
\begin{equation}\label{eq:tail-functional}
 \|f\|_{p,\alpha}
 :=\sup_{n\ge1}n^{-\alpha}\cN_{p,n}(f)
\end{equation}
for $\alpha\ge0$.
\end{definition}

\begin{remark}
Our definition is related to \cite{Kaimanovich1992}. More precisely, for symmetric $k_n$, the functional $2^{-1/p}\cN_{p,n}$ is the $p$-Dirichlet norm associated with $K_n$ in the sense of ~\cite[equation~(1.5)]{Kaimanovich1992}.
\end{remark}
\medskip

Averaging functionals satisfy the following elementary estimates which will repeatedly be used in the sequel.

\begin{lemma}
\label{lem:abstract-pseudo-Poincare}
Let $1\le p<\infty$, $n\ge1$, and $f\in\ell^p(V,m)$ such that $\cN_{p,n}(f)<\infty$. Then,
\begin{enumerate}[wide=0em,label={\rm (\roman*)}]
\item the average $K_nf$ is well defined and
\begin{equation}\label{eq:abstract-pseudo-Poincare}
 \|f-K_nf\|_p\le\cN_{p,n}(f);
\end{equation}
\item for any Lipschitz function $\Phi\colon\R\to\R$ with $\Phi(0)=0$,
\begin{equation}\label{eq:abstract-normal-contraction}
 \cN_{p,n}(\Phi\circ f)
 \le \Lip(\Phi)\cN_{p,n}(f).
\end{equation}
\end{enumerate}
\end{lemma}

\begin{proof}
\begin{enumerate}[wide=0em,label={\rm (\roman*)}]
\item Since $m(x)>0$ and $\cN_{p,n}(f)<\infty$, one has for every $x$
\[
\sum_{y\in V}|f(x)-f(y)|^p k_n(x,y)m(y)<+\infty.
\]
Moreover, since $k_n(x,y)m(y)$ is by assumption~\eqref{eq:averaging-normalization} a probability measure, H\"older's
inequality implies
\[
 \sum_{y\in V}|f(y)|k_n(x,y)m(y)
 \le |f(x)|+\left( \sum_{y\in V}|f(x)-f(y)|^p k_n(x,y)m(y) \right)^{1/p}<\infty.
\]
Thus, $K_nf(x)$ is well defined.  Jensen's inequality and Tonelli's
theorem now give
\[
 \|f-K_nf\|_p^p
 \le\sum_{x,y\in V}|f(x)-f(y)|^p
 k_n(x,y)m(x)m(y)
 =\cN_{p,n}(f)^p.
\]
\item The assertion follows pointwise from
$|\Phi(f(x))-\Phi(f(y))|\le\Lip(\Phi)|f(x)-f(y)|$.
\end{enumerate}
\end{proof}

The following property is a discrete analogue of~\cite[Definition~6.7]{AlonsoEtAl2020I} and~\cite[Definition~4.5]{Baudoin2024KS}, with scale one replacing the infinitesimal limit.

\begin{definition}\label{def:P}
Let $1\le p<\infty$ and $\alpha\ge0$. The family $(\cN_{p,n})_{n\ge1}$ satisfies $\mathrm{P}(p,\alpha)$ if there exists $C_{\rm P}>0$ such that
\begin{equation}\label{eq:P}\tag{$\mathrm{P}(p,\alpha)$}
 r^{-\alpha}\cN_{p,r}(f)\le C_{\rm P}\cN_{p,1}(f)
\end{equation}
for every $f\in\ell^p(V,m)$ and every integer $r\ge1$.
\end{definition}

Equivalently,
\begin{equation}\label{eq:P-seminorm-comparison}
 \cN_{p,1}(f)\le\|f\|_{p,\alpha}\le C_{\rm P}\cN_{p,1}(f)\,.
\end{equation}
Property $\mathrm{P}(p,\alpha)$ is preserved when $\alpha$ increases. By Lemma~\ref{lem:abstract-pseudo-Poincare}, it implies
\begin{equation}\label{eq:P-pseudo-Poincare}
 \|f-K_nf\|_p\le C_{\rm P}n^\alpha\cN_{p,1}(f),
 \qquad n\ge1\,,
\end{equation}
whenever $\cN_{p,1}(f)<\infty$. To compare two arbitrary scales we introduce the following concept of weak monotonicity.

\begin{definition}\label{def:WM}
Let $1\le p<\infty$ and $\alpha\ge0$. The family of functionals $(\cN_{p,n})_{n\ge1}$
satisfies the \hypertarget{WM}{\emph{weak monotonicity property}}
$(\WM)$ if there is $C_{\rm WM}>0$ such that
\begin{equation}\label{eq:WM}\tag{${\rm WM}_{p,\alpha}$}
 r^{-\alpha}\cN_{p,r}(f)
 \le C_{\rm WM}s^{-\alpha}\cN_{p,s}(f)
\end{equation}
for every $f\in\ell^p(V,m)$ and all integers $r\ge s\ge1$.
\end{definition}

Note that taking $s=1$ gives
\begin{equation*}
 (\WM)\Longrightarrow\mathrm{P}(p,\alpha)\,.
\end{equation*}
The advantage of weak monotonicity is that the functional inequalities below hold uniformly in the reference scale $s$, allowing the infimum over $s\ge1$.

\subsection{Gagliardo--Nirenberg inequalities}
Throughout this section, we make the following two quantitative additional assumptions:
\begin{equation}\label{eq:positive-atoms}
     m_*:=\inf_{x\in V}m(x)>0,
\end{equation}
and
\begin{equation}\label{eq:abstract-UC}
 k_n(x,y)\le C_{\rm U}n^{-\delta}
 \qquad\qquad \forall\;x,y\in V,\quad n\ge1,
\end{equation}
for some $\delta>0$. Note that the estimate above implies
\begin{equation}\label{eq:abstract-smoothing}
 \|K_nf\|_\infty\le C_{\rm U}n^{-\delta}\|f\|_1.
\end{equation}

The next theorem is the discrete version of the principle of
Bakry, Coulhon, Ledoux, and Saloff-Coste
\cite[Theorem~9.1]{BCLS1995} together with~\cite[Theorems~3.1--3.4,
Proposition~3.5 and Corollary~6.3]{BCLS1995}. The corresponding argument for
heat semigroup Besov functionals appears in
\cite[Section~4.1]{AlonsoRuizBaudoin2021} and a metric ball functional version can be found in~\cite[Theorem~4.5]{Baudoin2024KS}.

\begin{theorem}
\label{thm:abstract-GN}\label{thm:main-GN}
Assume $\mathrm{P}(p,\alpha)$ with $1\le p<\infty$ and $\alpha>0$, as well as~\eqref{eq:positive-atoms} and~\eqref{eq:abstract-UC}. Let $a,b\in(0,\infty]$ and $\theta\in(0,1]$ satisfy
\begin{equation}\label{eq:abstract-GN-parameters}
 \frac1a=\theta\left(\frac1p-\frac\alpha\delta\right)+\frac{1-\theta}{b}.
\end{equation}
If $p\alpha=\delta$ and $p>1$, assume in addition that $a<\infty$.
Then, for every $f\in\ell^p(V,m)$,
\begin{equation}\label{eq:abstract-GN}
 \|f\|_a
 \le C \|f\|_b^{1-\theta} \cN_{p,1}(f)^\theta.
\end{equation}
When $\theta=1$, the factor $\|f\|_b^{1-\theta}$ is omitted. The constant may depend on $p,\alpha,\delta,a,b,\theta$, on $m_*^{-1}$, and on $C_{\rm U},C_{\rm P}$, but is independent of $f$.
If $({\rm WM}_{p,\alpha})$ holds, then also
\begin{equation}\label{eq:abstract-GN-WM}
 \|f\|_a
 \le C\|f\|_b^{1-\theta}
 \left(\inf_{s\ge1}\cN_{p,s}(f)\right)^\theta\,,
\end{equation}
with dependence on $C_{\rm WM}$ in place of $C_{\rm P}$.
\end{theorem}

\begin{proof}
First consider nonnegative finitely supported functions $f$ such that $\cN_{p,1}(f)<\infty$. This class is stable
under truncations. By applying~\cite[Lemma~7.1]{BCLS1995} to the kernel $k_1(x,y)m(y)$ and the measure $m$, the seminorm $\cN_{p,1}$ satisfies the hypotheses $(H_p)$ and $(H_\infty^+)$ of that paper, with constants depending only on $p$. These are the truncation inequalities used below to pass from a weak interpolation estimate to strong inequalities.

For $t>0$, set
\begin{equation}\label{eq:Mt-reference-scale}
 \mathsf M_t f:=
 \begin{cases}
  f,&0<t<1,\\
  K_{n_t}f,&t\ge1,\qquad n_t:=\lceil t^{1/\alpha}\rceil.
 \end{cases}
\end{equation}
We claim that
\begin{align}
 \|\mathsf M_t f\|_\infty
 &\le Ct^{-\delta/\alpha}\|f\|_1,\label{eq:Mt-smoothing}\\
 \|f-\mathsf M_t f\|_p
 &\le Ct\cN_{p,1}(f),\label{eq:Mt-approximation}
\end{align}
with constants independent of $t$ and $f$. For $0<t<1$, the first estimate follows from $\|f\|_\infty\le m_*^{-1}\|f\|_1$ and the second is immediate. For $t\ge1$, ultracontractivity, Lemma~\ref{lem:abstract-pseudo-Poincare}, and $\mathrm{P}(p,\alpha)$ give
\begin{align*}
 \|K_{n_t}f\|_\infty
 &\le C_{\rm U}n_t^{-\delta}\|f\|_1
 \le C_{\rm U}t^{-\delta/\alpha}\|f\|_1,\\
 \|f-K_{n_t}f\|_p
 &\le\cN_{p,n_t}(f)
 \le C_{\rm P}n_t^\alpha\cN_{p,1}(f)
 \le Ct\cN_{p,1}(f),
\end{align*}
since $t^{1/\alpha}\le n_t\le2t^{1/\alpha}$. This proves the claim. Put $\nu=\delta/\alpha$ and $\tau=1+1/\nu$.
For $f\ne0$ and $\lambda>0$, choose $t$ so that $Ct^{-\nu}\|f\|_1=\lambda/3$.
Chebyshev's inequality and the two estimates above give
\begin{equation*}
 m(\{f\ge\lambda\})
 \le C\lambda^{-p\tau}\cN_{p,1}(f)^p\|f\|_1^{p/\nu}.
\end{equation*}
This is the weak interpolation estimate of~\cite[Theorem~9.1]{BCLS1995}, with $R=\infty$.
Its Sobolev parameter is $1/q=1/p-1/\nu$.
The truncation hypotheses upgrade this estimate to the strong inequalities by~\cite[Theorems~3.1--3.3 and Proposition~3.5]{BCLS1995}, with interpolation relation~\eqref{eq:abstract-GN-parameters}. Consequently,
\begin{equation}\label{eq:abstract-GN-fixed-scale}
 \|f\|_a
 \le C\|f\|_b^{1-\theta}\cN_{p,1}(f)^\theta
\end{equation}
for the asserted parameters when $p\alpha\ne\delta$.
When $p\alpha=\delta$, the weak inequality furnished by the same
theorem, together with~\cite[Theorem~3.3 and
Proposition~3.5]{BCLS1995}, gives
\eqref{eq:abstract-GN-fixed-scale} for finite $a$.
For $p=1$, \cite[Corollary~6.3(1)]{BCLS1995} also gives
$\|f\|_\infty\le C\cN_{1,1}(f)$, and interpolation yields the
remaining critical inequalities.

To extend the estimate to arbitrary functions, let $f\in\ell^p(V,m)$ with $\cN_{p,1}(f)<\infty$, and set
\[
 f_j=(|f|-j^{-1})_+\wedge j,
 \qquad j\ge1.
\]
The bound $m_*>0$ implies
\[
 \#\operatorname{supp}f_j
 \le m_*^{-1}j^p\|f\|_p^p<\infty.
\]
Moreover, $f_j\uparrow|f|$ and, by normal contraction,
$\cN_{p,1}(f_j)\le\cN_{p,1}(f)$. Applying
\eqref{eq:abstract-GN-fixed-scale} to $f_j$ and passing to the limit
therefore proves the same inequality for $f$ whenever its right-hand
side is finite. Here one uses monotone convergence for finite
exponents and $\|f_j\|_\infty\uparrow\|f\|_\infty$ for the exponent
$\infty$; when $\theta=1$, the factor involving $\|f_j\|_b$ is omitted.
This proves~\eqref{eq:abstract-GN}.

Now assume $({\rm WM}_{p,\alpha})$ and fix an integer $s\ge1$. The kernel family $(k_{sn})_{n\ge1}$ satisfies $\mathrm{P}(p,\alpha)$ with constant $C_{\rm WM}$, since
\begin{equation*}
 n^{-\alpha}\cN_{p,sn}(f)\le C_{\rm WM}\cN_{p,s}(f),
 \qquad n\ge1\,.
\end{equation*}
Moreover, $k_{sn}(x,y)\le C_{\rm U}s^{-\delta}n^{-\delta}\le C_{\rm U}n^{-\delta}$. Applying~\eqref{eq:abstract-GN} to this family gives
\begin{equation*}
 \|f\|_a\le C\|f\|_b^{1-\theta}\cN_{p,s}(f)^\theta
\end{equation*}
with $C$ independent of $s$. Taking the infimum over $s\ge1$ proves~\eqref{eq:abstract-GN-WM}.
\end{proof}

We point out the following  special cases which are of special interest.

\begin{corollary}
\label{cor:abstract-GN-specializations}
Under the hypotheses of Theorem~\ref{thm:abstract-GN}, the following statements hold:
\begin{enumerate}[wide=0em,label={\rm(\roman*)}]
 \item (Sobolev inequality) If $p\alpha<\delta$ and
 $q=p\delta/(\delta-p\alpha)$, then
 \begin{equation}\label{eq:abstract-Sobolev}
  \|f\|_q\le C\cN_{p,1}(f).
 \end{equation}
 \item (Nash inequality) If $p>1$, then
 \begin{equation}\label{eq:abstract-Nash}
  \|f\|_p
  \le C\cN_{p,1}(f)^{\eta_p}\|f\|_1^{1-\eta_p},
  \qquad
  \eta_p:=\frac{(p-1)\delta}{p(\alpha+\delta)-\delta}.
 \end{equation}
 \item If $p\alpha=\delta$ and $p>1$, then, for
 $0<b<a<\infty$,
 \begin{equation}\label{eq:abstract-critical-GN}
  \|f\|_a\le C\cN_{p,1}(f)^{1-b/a}\|f\|_b^{b/a}.
 \end{equation}
 \item If $p\alpha>\delta$, or if $p\alpha=\delta$ and $p=1$, then,
 for $1\le b<\infty$,
 \begin{equation}\label{eq:abstract-supercritical-GN}
  \|f\|_\infty
  \le C\cN_{p,1}(f)^{\theta_b}\|f\|_b^{1-\theta_b},
  \qquad
  \theta_b:=\frac{p\delta}
  {p\delta+b(p\alpha-\delta)}.
 \end{equation}
 If $\operatorname{supp}f\subset\Omega$ and $m(\Omega)<\infty$, then
 \begin{equation}\label{eq:abstract-supported-GN}
  \|f\|_\infty
  \le C\cN_{p,1}(f)m(\Omega)^{\alpha/\delta-1/p}.
 \end{equation}
\end{enumerate}
If $({\rm WM}_{p,\alpha})$ holds, each occurrence of $\cN_{p,1}(f)$ may be replaced by $\inf_{s\ge1}\cN_{p,s}(f)$.
\end{corollary}

\begin{proof}
Use respectively $\theta=1$; $a=p$ and $b=1$; $p\alpha=\delta$ and
$\theta=1-b/a$; and $a=\infty$ in
Theorem~\ref{thm:abstract-GN}.  The supported estimate follows from
\eqref{eq:abstract-supercritical-GN} with $b=1$ and
$\|f\|_1\le\|f\|_\infty m(\Omega)$.
\end{proof}

At the critical exponent, the truncation argument also gives exponential integrability and a logarithmic estimate for functions of finite support.

\begin{corollary}\label{cor:abstract-Trudinger}
Assume the hypotheses of Theorem~\ref{thm:abstract-GN}, with $p>1$ and $p\alpha=\delta$. There exist constants $c,C>0$ such that, for every $\Omega\subset V$ with $0<m(\Omega)<\infty$ and every nonzero $f\in\ell^p(V,m)$ supported in $\Omega$ with $\cN_{p,1}(f)<\infty$,
\begin{equation}\label{eq:abstract-Trudinger}
 \sum_{x\in\Omega}
 \exp\!\left[c\left(\frac{|f(x)|}{\cN_{p,1}(f)}\right)^{p/(p-1)}\right]m(x)
 \le C m(\Omega)\,.
\end{equation}
Consequently,
\begin{equation}\label{eq:abstract-critical-logarithmic}
 \|f\|_\infty
 \le C\cN_{p,1}(f)
 \left[\log\!\left(\frac{e\,m(\Omega)}{m_*}\right)\right]^{1-1/p}\,.
\end{equation}
The latter inequality also holds for $f=0$. Under $({\rm WM}_{p,\alpha})$, each occurrence of $\cN_{p,1}(f)$ may be replaced by $\inf_{s\ge1}\cN_{p,s}(f)$. The constants depend only on $p,\alpha,\delta,m_*^{-1},C_{\rm U},C_{\rm P}$, with $C_{\rm WM}$ in place of $C_{\rm P}$ for the weak monotonicity conclusion.
\end{corollary}

\begin{proof}
Since $m_*>0$, the set $\Omega$ is finite. Inequality~\eqref{eq:abstract-critical-GN}, with $a=2p$ and $b=p$, gives $\cN_{p,1}(f)>0$ for nonzero $f$. By~\cite[Lemma~7.1]{BCLS1995}, this seminorm satisfies the uniform truncation condition $(H_p)$ of that paper. Thus~\eqref{eq:abstract-critical-GN} and~\cite[Corollary~6.4]{BCLS1995}, applied to $|f|$, give~\eqref{eq:abstract-Trudinger}. At a point where $|f|$ attains its maximum,
\begin{equation*}
 m_*\exp\!\left[c\left(\frac{\|f\|_\infty}{\cN_{p,1}(f)}\right)^{p/(p-1)}\right]
 \le C m(\Omega)\,.
\end{equation*}
Taking logarithms proves~\eqref{eq:abstract-critical-logarithmic}, after changing the constant.

Under $({\rm WM}_{p,\alpha})$, apply the preceding argument to $(k_{sn})_{n\ge1}$ whenever $\cN_{p,s}(f)<\infty$, as in the proof of Theorem~\ref{thm:abstract-GN}. The constants are independent of $s$. Moreover,~\eqref{eq:abstract-GN-WM}, with $a=2p$ and $b=p$, gives $\inf_{s\ge1}\cN_{p,s}(f)>0$ for nonzero $f$. Taking a sequence of scales approaching this infimum and applying Fatou's lemma proves the assertion.
\end{proof}

As we explain, for $p=1$, the Sobolev inequality has an isoperimetric interpretation. Similarly to~\cite[Section~6.3]{AlonsoEtAl2020I}, we associate a perimeter to the kernel family $k=(k_n)_{n\ge1}$.

\begin{definition}\label{def:kernel-alpha-perimeter}
Let $\alpha>0$ and $A\subset V$ with $m(A)<\infty$. Its $\alpha$-perimeter associated with $k$ is
\begin{equation}\label{eq:kernel-alpha-perimeter}
 \mathcal P_{\alpha,k}(A):=\|\1_A\|_{1,\alpha}
 =\sup_{n\ge1}n^{-\alpha}
 \sum_{\substack{x\in A\\y\notin A}}
 \bigl(k_n(x,y)+k_n(y,x)\bigr)m(x)m(y)\,.
\end{equation}
We say that $A$ has finite $\alpha$-perimeter if $\mathcal P_{\alpha,k}(A)<\infty$.
\end{definition}

If the kernels are symmetric, the sum in~\eqref{eq:kernel-alpha-perimeter} is twice $\sum_{x\in A,\,y\notin A}k_n(x,y)m(x)m(y)$, the boundary size associated with $K_n$ in~\cite[equation~(1.8)]{Kaimanovich1992}. Under $\mathrm{P}(1,\alpha)$, one has
\begin{equation}\label{eq:kernel-perimeter-comparison}
 \cN_{1,1}(\1_A)\le\mathcal P_{\alpha,k}(A)
 \le C_{\rm P}\cN_{1,1}(\1_A)\,.
\end{equation}

The next consequence is the discrete counterpart of~\cite[Proposition~6.5 and Corollary~6.6]{AlonsoEtAl2020I}. For a single symmetric kernel, the equivalence between $L^1$ Sobolev and isoperimetric inequalities is contained in~\cite[Theorems~4.1 and~4.3]{Kaimanovich1992}.

\begin{corollary}\label{cor:kernel-isoperimetry}
Assume the hypotheses of Theorem~\ref{thm:abstract-GN} with $p=1$ and $0<\alpha\le\delta$. Set $q=\delta/(\delta-\alpha)$ if $\alpha<\delta$, and $q=\infty$ if $\alpha=\delta$. Then
\begin{equation}\label{eq:kernel-L1-Sobolev}
 \|f\|_q\le C\|f\|_{1,\alpha},
 \qquad f\in\ell^1(V,m)\,.
\end{equation}
In particular, every $A\subset V$ with $0<m(A)<\infty$ satisfies
\begin{equation}\label{eq:kernel-isoperimetry}
 m(A)^{1-\alpha/\delta}\le C\mathcal P_{\alpha,k}(A)\,.
\end{equation}
When $\alpha=\delta$, the inequality is the uniform lower bound $\mathcal P_{\alpha,k}(A)\ge c$ for nonempty sets of finite measure. More generally, assuming only $\mathrm{P}(1,\alpha)$ in addition to the kernel normalization, inequalities~\eqref{eq:kernel-L1-Sobolev} and~\eqref{eq:kernel-isoperimetry} are equivalent up to constants.
\end{corollary}

\begin{proof}
Corollary~\ref{cor:abstract-GN-specializations}, parts~{\rm(i)} and~{\rm(iv)}, gives~\eqref{eq:kernel-L1-Sobolev}, since $\cN_{1,1}(f)\le\|f\|_{1,\alpha}$. Apply this inequality to $f=\1_A$ to obtain~\eqref{eq:kernel-isoperimetry}.

For the converse, assume only $\mathrm{P}(1,\alpha)$ and~\eqref{eq:kernel-isoperimetry}. We use the classical coarea argument of Federer--Fleming and Maz'ya relating $L^1$ Sobolev and isoperimetric inequalities; see~\cite[equation~(2.6)]{Kaimanovich1992} and~\cite[Chapter~8]{Ledoux1996}. For $0\le f\in\ell^1(V,m)$, Tonelli's theorem and the identity
\begin{equation*}
 |f(x)-f(y)|
 =\int_0^\infty
 |\1_{\{f>t\}}(x)-\1_{\{f>t\}}(y)|\,dt
\end{equation*}
give, for every $n\ge1$,
\begin{equation*}
 \cN_{1,n}(f)
 =\int_0^\infty\cN_{1,n}(\1_{\{f>t\}})\,dt\,.
\end{equation*}
Taking the supremum for the lower bound and using~\eqref{eq:kernel-perimeter-comparison} for the upper bound yields
\begin{equation}\label{eq:kernel-coarea}
 \|f\|_{1,\alpha}
 \le\int_0^\infty\mathcal P_{\alpha,k}(\{f>t\})\,dt
 \le C_{\rm P}\cN_{1,1}(f)
 \le C_{\rm P}\|f\|_{1,\alpha}\,.
\end{equation}
Since $m(\{f>t\})<\infty$ for $t>0$, Minkowski's inequality and~\eqref{eq:kernel-isoperimetry} now imply
\begin{align*}
 \|f\|_q
 &\le\int_0^\infty\|\1_{\{f>t\}}\|_q\,dt\\
 &\le C\int_0^\infty\mathcal P_{\alpha,k}(\{f>t\})\,dt
 \le CC_{\rm P}\|f\|_{1,\alpha}\,.
\end{align*}
The same argument applies when $q=\infty$, with $\|\1_{\{f>t\}}\|_\infty$ equal to $1$ for nonempty superlevel sets and $0$ otherwise. For general real $f$, we apply the estimate to $|f|$ and use normal contraction.
\end{proof}

\section{Weighted graphs and heat kernels}\label{sec:graphs}

In this section, we give the definitions and label the assumptions we are working with in the framework of weighted graphs.

\subsection{Weighted graphs and Markov operators}
Let $G=(V,E)$ be a connected, locally finite, countable graph. We equip the graph with symmetric conductances $(c_{xy})_{x,y\in V}$ and take the vertex measure to be
\begin{equation}\label{E:def_conductance_measure}
    m(x):=\sum_{y\in V}c_{xy},
\end{equation}
where the conductances are nonnegative and $0<m(x)<\infty$ for every $x\in V$. For $x\ne y$, one has $c_{xy}>0$ if and only if $\{x,y\}\in E$; the diagonal conductances $c_{xx}\ge0$ are arbitrary, so loops are allowed. The \emph{graph distance} between $x,y\in V$ is
\[
d(x,y):=\min\bigl\{n\ge0:\ \exists\,x_0,\ldots,x_n\in V,\ 
x_0=x,\ x_n=y,\ \{x_{i-1},x_i\}\in E
\text{ for }1\le i\le n\bigr\}.
\]
Thus $d(x,y)$ is the minimum number of edges in a path joining
$x$ to $y$, with $d(x,x)=0$. Each edge has length one,
independently of its conductance.

We write
\[
 \Vol(x,r):=m(B(x,r)),
 \qquad
 B(x,r):=\{y\in V:d(x,y)\le r\}.
\]

\medskip

On this weighted graph, there is a naturally associated Markov chain $X:=\{X_n\}_{n\geq 0}$, whose law and expectation when started from $x\in V$ will be denoted by $\mathbb P^x$ and $\mathbb E^x$, respectively. The corresponding transition probability kernel and Markov operator are given by
\begin{equation}\label{eq:transition}
    \pi(x,y):=\frac{c_{xy}}{m(x)}
    \qquad\text{and}\qquad
    Pf(x):=\sum_{y\in V}\pi(x,y)f(y).
\end{equation}
We use the latter to define the positive normalized Laplacian
\begin{equation}\label{eq:laplacian}
    L:=I-P,
\end{equation}
which is a bounded operator on every $\ell^p(V,m)$ with norm at most $2$, and its associated Dirichlet form on $\ell^2(V,m)$ is given by
\begin{equation}\label{eq:def_graph_DF}
    \begin{aligned}
    \cE(f,g)&:=\langle Lf,g\rangle_{\ell^2}
    =\frac12\sum_{x,y\in V}c_{xy}\bigl(f(x)-f(y)\bigr)\bigl(g(x)-g(y)\bigr),\\
    \cF &:=\ell^2(V,m).
    \end{aligned}
\end{equation}
Since $\cE(f,f)\le2\|f\|_2^2$, the form norm is equivalent to the $\ell^2$ norm, and finitely supported functions form a core. Thus the form is regular for the discrete topology on $V$; see for instance ~\cite[Lemma~1.16]{KellerLenzWojciechowski2021}.

Moreover, $\cE(f,f)=\frac12\|\nabla_2f\|_2^2$, where for a general $1\leq p<\infty$ the discrete $p$-gradient is defined as
\begin{equation}\label{eq:p-gradient}
    |\nabla_p f|(x)
    :=\bigg(\sum_{y\in V}\pi(x,y)|f(x)-f(y)|^p\bigg)^{1/p}
\end{equation}
and thus
\begin{equation}\label{eq:gradient-edge-identity}
    \|\nabla_p f\|_p^p
    =\sum_{x,y\in V}c_{xy}|f(x)-f(y)|^p.
\end{equation}
For $A\subset V$, its weighted edge perimeter is
\begin{equation}\label{eq:graph-perimeter}
    \Per(A):=\sum_{\substack{x\in A\\y\notin A}}c_{xy}\,.
\end{equation}
Let $\pi_n(x,y)$, $n\geq 0$, denote the $n$-step transition probabilities of the Markov chain $X$. The $n$-step heat kernel is defined as the density of $\pi_n$ with respect to $m$, i.e.
\begin{equation}\label{eq:heat-kernel}
    p_n(x,y):=\frac{\pi_n(x,y)}{m(y)},
    \qquad n\geq 0.
\end{equation}
In particular, $\{p_n\}_{n\geq 0}$ is a sequence of nonnegative symmetric kernels satisfying~\eqref{eq:averaging-normalization} as well as
\begin{equation}\label{eq:kernel-properties}
\begin{aligned}
    p_0(x,y)&=\frac{1}{m(y)}\1_{\{x\}}(y),\\
    p_{n+k}(x,y)&=\sum_{z\in V}p_n(x,z)p_k(z,y)m(z),\qquad k\geq 0, n\geq 1.
\end{aligned}
\end{equation}
The associated averaging operators $P^n$ form a semigroup $P^{n+k}=P^n P^k$, are positive contractions on $\ell^p(V,m)$ for every $1\leq p\leq\infty$ and they are self-adjoint on $\ell^2(V,m)$.

\vspace{0.4cm}

For $\alpha>0$ and $A\subset V$ with $m(A)<\infty$, the heat kernel $\alpha$-perimeter is $\mathcal P_{\alpha,(p_n)}(A)$ in the notation of Definition~\ref{def:kernel-alpha-perimeter}. Since $p_1(x,y)m(x)m(y)=c_{xy}$, one has $\cN_{1,1}(\1_A)=2\Per(A)$. Hence, whenever the heat kernel functionals satisfy $\mathrm{P}(1,\alpha)$,~\eqref{eq:kernel-perimeter-comparison} gives
\begin{equation}\label{eq:heat-kernel-perimeter-comparison}
 2\Per(A)\le\mathcal P_{\alpha,(p_n)}(A)
 \le2C_{\rm P}\Per(A)\,.
\end{equation}
Thus the heat kernel $\alpha$-perimeter and the weighted edge perimeter are comparable under $\mathrm{P}(1,\alpha)$. In particular, under the hypotheses of Corollary~\ref{cor:kernel-isoperimetry} with $k_n=p_n$, its conclusion implies $m(A)^{1-\alpha/\delta}\le C\Per(A)$ for $0<m(A)<\infty$.

\subsection{Controlled weights, uniform diagonal lower bounds and volume growth}
We introduce the weight and volume conditions used in the sequel.

\begin{definition}\label{def:weight-assumptions}
A weighted graph $(V,E)$ with transition probability kernel $\pi(x,y)$ is said to satisfy the \emph{uniform edge control property}, denoted~\eqref{eq:p0}, if there exists $\varepsilon_0>0$ such that
\begin{equation}\label{eq:p0}
    \pi(x,y)\geq\varepsilon_0
    \qquad\text{for }x\sim y,\;x\neq y.\tag{${\rm p}_0$}
\end{equation}
The graph satisfies the \emph{diagonal lower bound}~\eqref{eq:LB} if there exists $\varepsilon_{\rm L}>0$ such that
\begin{equation}\label{eq:LB}
    \pi(x,x)\geq\varepsilon_{\rm L},
    \qquad x\in V.
    \tag{${\rm LB}$}
\end{equation}
We say that the graph has \emph{controlled weights} when both conditions hold.
\end{definition}

\begin{remark}\label{rem:p_0&LB}
Since $\pi(x,y)=p_1(x,y)m(y)$, controlled weights are equivalent to the existence of $\varepsilon>0$ such that
\begin{equation*}
 p_1(x,y)m(y)\ge\varepsilon,
 \qquad d(x,y)\le1\, ,
\end{equation*}
where one may take $\varepsilon=\min\{\varepsilon_0,\varepsilon_{\rm L}\}$.

\end{remark}
\begin{remark}
We note that condition~\eqref{eq:LB} excludes $-1$ from the spectrum of $P$ on $\ell^2(V,m)$ and therefore implies discrete analyticity of the semigroup $\{P^n\}_{n\geq 1}$ on $\ell^2(V,m)$; see
\cite[Section~III, Proposition~3, and Section~V, Remarque~2]{CoulhonSaloffCoste1990},
\cite[Lemma~1.3]{Dungey2006}, and
\cite[Section~1]{Feneuil2018}.
\end{remark}
\begin{remark}\label{rem:controlled-weight-comparison}
If~\eqref{eq:p0} holds and $\Vol(x,1)\simeq1$ uniformly in $x$, then $m(x)\simeq1$ and $c_{xy}\simeq1$ on off-diagonal edges. Indeed, the off-diagonal degree is at most $\varepsilon_0^{-1}$, and reversibility gives $m(y)\le\varepsilon_0^{-1}m(x)$ for $x\sim y$. Thus
\begin{equation*}
 \Vol(x,1)\le(1+\varepsilon_0^{-2})m(x)\,.
\end{equation*}
The lower bound on $m(x)$ follows, while $m(x)\le\Vol(x,1)$ gives the upper bound. Finally, $\varepsilon_0m(x)\le c_{xy}\le m(x)$ for $x\sim y$, $x\ne y$.
\end{remark}

\begin{definition}\label{def:volume-growth}
Let $d_f>0$. The graph satisfies the \emph{volume growth condition}~\eqref{eq:volume-growth} if there exist constants $0<c_V\le C_V<\infty$ such that
\begin{equation}\label{eq:volume-growth}
 c_Vr^{d_f}\le\Vol(x,r)\le C_Vr^{d_f},
 \qquad x\in V,\quad r\ge1\,.
 \tag{$\mathrm{V}_{d_f}$}
\end{equation}
\end{definition}
\medskip

\subsection{Sub-Gaussian heat kernel estimates}
The definition of the discrete sub-Gaussian heat kernel estimates below is standard and we refer the reader to~\cite{GrigoryanTelcs2001,BCK2005,Barlow2017} for further details.

\begin{definition}\label{def:sub-gaussian}
Fix $d_f>0$ and $d_w\geq2$. A weighted graph satisfies the condition~$({\rm SG}_{d_f,d_w})$ if there exist positive constants $c_i$ such that
\begin{equation}\label{eq:SG-upper}
 p_n(x,y)
 \leq c_1n^{-d_f/d_w}
 \exp\left[-c_2
  \left(\frac{d(x,y)}{n^{1/d_w}}\right)^{d_w/(d_w-1)}\right],
\end{equation}
for every $x,y\in V$ and $n\geq1$, and
\begin{equation}\label{eq:SG-lower}
 p_n(x,y)+p_{n+1}(x,y)
 \geq c_3n^{-d_f/d_w}
 \exp\left[-c_4\left(\frac{d(x,y)}{n^{1/d_w}}\right)^{d_w/(d_w-1)}\right]
\end{equation}
for any $x,y\in V$ with $d(x,y)\leq n$.
\end{definition}

\begin{proposition}\label{prop:volume-growth}
Assume~\eqref{eq:LB} and $({\rm SG}_{d_f,d_w})$. Then~\eqref{eq:volume-growth} holds. In particular, the measure is doubling.
\end{proposition}

\begin{proof}

Fix $x\in V$ and $r\geq1$, and set $n=\lceil r^{d_w}\rceil$.  If
$y\in B(x,r)$, then $d(x,y)\leq r$ and the exponential factor in
\eqref{eq:SG-lower} is bounded below uniformly.  Hence
\[
 2=\sum_{y\in V}\bigl(p_n(x,y)+p_{n+1}(x,y)\bigr)m(y)
 \geq c r^{-d_f}\Vol(x,r),
\]
which proves $\Vol(x,r)\leq Cr^{d_f}$. In particular, we have
\begin{equation}\label{eq:uniform-atom-upper-bound}
    m(x)\le\Vol(x,1)\le C, \qquad x\in V.
\end{equation}

Next we note that~\eqref{eq:LB} and~\eqref{eq:SG-upper} imply
\[
 \varepsilon_{\rm L}\leq\pi(x,x)=p_1(x,x)m(x)\leq C m(x)
\]
whence
\begin{equation}\label{eq:uniform-atom-lower-bound}
 m(x)\geq c \qquad \forall\,x\in V.
\end{equation}

For the lower volume bound, choose $K\ge2$ and, for $r\ge2K$, put $n=\lfloor(r/K)^{d_w}\rfloor$. The heat kernel upper bound and the upper volume bound give
\begin{align*}
 \sum_{d(x,y)>r}p_n(x,y)m(y)
 &\le C\sum_{j\ge0}n^{-d_f/d_w}(2^{j+1}r)^{d_f}
 \exp\left[-c\left(\frac{2^jr}{n^{1/d_w}}\right)^{d_w/(d_w-1)}\right]\\
 &\le C\sum_{j\ge0}(2^jK)^{d_f}
 \exp\left[-c(2^jK)^{d_w/(d_w-1)}\right]\,.
\end{align*}
Fix $K$ large enough that this sum is at most $1/2$. Conservativeness now yields
\begin{equation*}
 \frac12\le\sum_{y\in B(x,r)}p_n(x,y)m(y)
 \le Cn^{-d_f/d_w}\Vol(x,r)\,,
\end{equation*}
and hence $\Vol(x,r)\ge cr^{d_f}$. For $1\le r<2K$, the same conclusion follows from~\eqref{eq:uniform-atom-lower-bound}.
\end{proof}

\section{Heat kernel functional weak monotonicity}\label{sec:HK-Besov}

We prove weak monotonicity for the heat kernel functionals under~\eqref{eq:LB}, with exponent $1/\min\{p,2\}$. Weak Bakry--\'Emery estimates improve this exponent for $1\le p<2$.

\subsection{Weak monotonicity under the diagonal lower bound}
\label{sec:WM-p-ge-two}

The sequence of kernels $\{p_n\}_{n\geq0}$ fits into the framework
of Section~\ref{subsec:WM}, with $k_n(x,y)=p_n(x,y)$.
Following Definition~\ref{def:averaging-functional}, we set
\begin{equation}\label{eq:Npn}
 \cN_{p,n}^{\rm H}(f)
 :=\left(
   \sum_{x,y\in V}|f(x)-f(y)|^p
   p_n(x,y)m(x)m(y)
  \right)^{1/p}
\end{equation}
for $1\leq p<\infty$, $n\geq1$, and $f\in\ell^p(V,m)$.
We write $\|f\|_{p,\alpha}=\sup_{n\ge1}n^{-\alpha}\cN_{p,n}^{\rm H}(f)$ for $\alpha\ge0$, as in Definition~\ref{def:averaging-functional}. This seminorm is the graph analogue of the heat kernel based Besov seminorm introduced in \cite{AlonsoEtAl2020I}.

\vspace{0.4cm}

Conservativeness and symmetry give, for every $n\ge1$,
\begin{equation}\label{eq:basic-upper}
 \cN_{p,n}^{\rm H}(f)^p
 \le2^{p-1}\sum_{x,y\in V}\bigl(|f(x)|^p+|f(y)|^p\bigr)
 p_n(x,y)m(x)m(y)
 =2^p\|f\|_p^p.
\end{equation}
Moreover, $\cN_{p,1}^{\rm H}(f)=\|\nabla_pf\|_p$. Hence $\|\nabla_pf\|_p\le\|f\|_{p,\alpha}\le2\|f\|_p$ for every $\alpha\ge0$. In particular, we always have
\[
\left\{ f \in \ell^p(V,m) : \| f \|_{p,\alpha} <+\infty \right\}=\ell^p(V,m),
\]
which is one of the main differences between the theory we develop here and \cite{AlonsoEtAl2020I}.

\vspace{0.4cm}

A first useful observation is the following.

\begin{proposition}\label{prop:l2-identity}
Let $1\le p\le2$. For every $f\in\ell^p(V,m)$ and integers $n,s\ge1$,
\begin{equation}\label{eq:heat-small-p-multiples}
 \cN_{p,ns}^{\rm H}(f)^p\le n\cN_{p,s}^{\rm H}(f)^p\,.
\end{equation}
Consequently, for every $\alpha\ge1/p$,
\begin{equation}\label{eq:small-p-identity}
 \|f\|_{p,\alpha}=\|\nabla_pf\|_p\,.
\end{equation}
In particular, for $p=2$ and $\alpha\ge1/2$,
\begin{equation}\label{eq:l2-identity}
 \|f\|_{2,\alpha}^2=2\cE(f,f)\,.
\end{equation}
\end{proposition}

\begin{proof}
For $p=2$, we  apply the spectral theorem to the self-adjoint contraction $P^s$ and the inequality $1-\lambda^n\le n(1-\lambda)$ on $[-1,1]$ gives~\eqref{eq:heat-small-p-multiples}.

Let $1\le p<2$. For real $a,b$, a change of variables gives
\begin{equation}\label{eq:phase-quadratic-identity}
 \int_0^\infty |e^{ita}-e^{itb}|^2\,\frac{dt}{t^{1+p}}
 =c_p|a-b|^p\,,\qquad 0<c_p<\infty\,.
\end{equation}
Let $f\in\ell^p(V,m)$. For each $t>0$, the function $e^{itf}-1$ belongs to $\ell^2(V,m)$, because
\begin{equation*}
 |e^{itf}-1|^2\le2^{2-p}t^p|f|^p\,.
\end{equation*}
Apply the $p=2$ extimate \eqref{eq:heat-small-p-multiples} to this function and then integrate against $t^{-1-p}\,dt$. Tonelli's theorem proves~\eqref{eq:heat-small-p-multiples}. Taking $s=1$ gives~\eqref{eq:small-p-identity}, since $n^{1/p-\alpha}\le1$ and equality holds at $n=1$. The last identity follows from $\|\nabla_2f\|_2^2=2\cE(f,f)$.
\end{proof}

The following lemma compares the heat-kernel functionals at different times and will be used repeatedly in the sequel.

\begin{lemma}
\label{lem:time-subadditivity}
\label{lem:heat-earlier-time}
Let $1\leq p<\infty$ and $f\in\ell^p(V,m)$. Suppose that
\eqref{eq:LB} holds with constant $\varepsilon_{\rm L}$.
Then, for all integers $1\leq r\leq s$,
\begin{equation}\label{eq:heat-earlier-time-comparison}
 \cN_{p,r}^{\rm H}(f)
 \leq 2\varepsilon_{\rm L}^{-1/p}
       \cN_{p,s}^{\rm H}(f).
\end{equation}
\end{lemma}

\begin{proof}
We first show that, for all integers $1\leq k\leq s$,
\begin{equation}\label{eq:heat-even-time-comparison}
 \cN_{p,2k}^{\rm H}(f)
 \leq 2\cN_{p,s}^{\rm H}(f).
\end{equation}
By the Chapman-Kolmogorov identity, symmetry, and conservativeness of the heat kernel, one has
\[
\begin{aligned}
 &\sum_{x,y\in V}|f(x)-f(y)|^p p_{2k}(x,y)m(x)m(y)\\
 &\quad=
 \sum_{x,y,z,\tilde z\in V}
 |f(x)-f(y)|^p
 p_k(x,\tilde z)p_k(y,\tilde z)p_{s-k}(z,\tilde z)
 m(x)m(y)m(z)m(\tilde z).
\end{aligned}
\]
Next we apply Minkowski's inequality to
\[
 f(x)-f(y)=(f(x)-f(z))+(f(z)-f(y))
\]
with respect to the measure on $V^4$ appearing in the last sum. Then it follows again from the Chapman-Kolmogorov identity, symmetry, and conservativeness  that each term on the resulting right-hand side equals $\cN_{p,s}^{\rm H}(f)$. This proves \eqref{eq:heat-even-time-comparison}.

Now let $1\leq r\leq s$. If $r$ is even, applying \eqref{eq:heat-even-time-comparison} with $k=r/2$ gives \eqref{eq:heat-earlier-time-comparison}. If $r$ is odd, condition \eqref{eq:LB} implies
\[
 p_{r+1}(x,y)\geq \varepsilon_{\rm L}p_r(x,y).
\]
Since $(r+1)/2\leq s$, we may apply
\eqref{eq:heat-even-time-comparison} with $k=(r+1)/2$ to obtain
\[
 \cN_{p,r}^{\rm H}(f)
 \leq \varepsilon_{\rm L}^{-1/p}\cN_{p,r+1}^{\rm H}(f)
 \leq 2\varepsilon_{\rm L}^{-1/p}\cN_{p,s}^{\rm H}(f).
\]
This completes the proof.
\end{proof}

Taking $r=1$ in~\eqref{eq:heat-earlier-time-comparison} gives, under~\eqref{eq:LB},
\begin{equation}\label{eq:heat-infimum-gradient}
 \frac{\varepsilon_{\rm L}^{1/p}}2\|\nabla_pf\|_p
 \le\inf_{s\ge1}\cN_{p,s}^{\rm H}(f)
 \le\|\nabla_pf\|_p,
 \qquad 1\le p<\infty\,.
\end{equation}
Thus the infimum $\inf_{s\ge1}\cN_{p,s}(f)$ is equivalent to the edge gradient norm for heat kernels under~\eqref{eq:LB}.

\vspace{0.6cm}

For $p\geq2$, the linear growth at integer multiples in
Lemma~\ref{lem:time-subadditivity} improves to square-root growth.

\begin{lemma}\label{lem:heat-p-ge-two-multiples}
Let $2\leq p<\infty$. There exists $A_p>0$ such that
\begin{equation}\label{eq:heat-p-ge-two-multiples}
 \cN_{p,ns}^{\rm H}(f)
 \leq A_p\sqrt n\,\cN_{p,s}^{\rm H}(f)
\end{equation}
for all integers $n,s\geq1$ and all $f\in\ell^p(V,m)$.
\end{lemma}

\begin{proof}
For $p=2$, the assertion is~\eqref{eq:heat-small-p-multiples} in Proposition~\ref{prop:l2-identity}.

For $p>2$, we use the moment estimate
from~\cite[Theorem~4.5]{NaorPeresSchrammSheffield2006}:
\begin{equation}\label{eq:moment_estimate_Naor}
 \left(\mathbb E|g(Z_n)-g(Z_0)|^p\right)^{1/p}
 \leq16\sqrt{pn}\,
      \left(\mathbb E|g(Z_1)-g(Z_0)|^p\right)^{1/p},
\end{equation}
where $(Z_j)_{j\geq0}$ is a finite stationary reversible
Markov chain and $g$ is a real-valued function on its state space. Fix $n,s\geq1$ and first suppose that $f$ has finite nonempty
support. Since the graph is locally finite, the set
\[
 D_f:=\{x\in V:d(x,\operatorname{supp}f)\leq ns+1\}
\]
is also finite. On $D_f$, we consider the Markov chain with  transition kernel
\[
 \pi_{D_f}(x,y):=
 \begin{cases}
  \pi(x,y),&x\neq y,\\[2mm]
  \displaystyle\pi(x,x)+\sum_{z\notin D_f}\pi(x,z),
    &x=y.
 \end{cases}
\]
Note that $\pi_{D_f}$ is reversible with respect to the probability measure $m|_{D_f}/m(D_f)$. 
Let $(Z_j^f)_{j\geq0}$ be the corresponding chain, started with this stationary distribution.

For every $1\leq j\leq ns$, the time-$j$ functional of the
restricted chain, computed with the unnormalized measure
$m|_{D_f}$ and applied to $f|_{D_f}$, equals
$\cN_{p,j}^{\rm H}(f)$. Indeed, any path contributing a nonzero
term has an endpoint in $\operatorname{supp}f$ and length
at most $ns$. Every vertex of such a path lies within distance
$ns$ of $\operatorname{supp}f$, where the restriction leaves
all transitions unchanged. Consequently,
\[
 \mathbb E|f(Z_j^f)-f(Z_0^f)|^p
 =\frac{\cN_{p,j}^{\rm H}(f)^p}{m(D_f)},
 \qquad 1\leq j\leq ns.
\]

Apply \eqref{eq:moment_estimate_Naor} to the stationary reversible
chain $(Z_{js}^f)_{j\geq0}$ and the function $f|_{D_f}$.
Multiplying by $m(D_f)^{1/p}$ gives
\eqref{eq:heat-p-ge-two-multiples} for finitely supported $f$;
the case $f=0$ is immediate.

Finally, finitely supported functions are dense in $\ell^p(V,m)$,
and Minkowski's inequality and~\eqref{eq:basic-upper} give
\[
 \bigl|\cN_{p,j}^{\rm H}(f)-\cN_{p,j}^{\rm H}(g)\bigr|
 \leq\cN_{p,j}^{\rm H}(f-g)
 \leq2\|f-g\|_p.
\]
Passing to the limit proves the result for every
$f\in\ell^p(V,m)$.
\end{proof}

Combining these time comparisons yields weak monotonicity for every finite~$p\ge1$.

\begin{theorem}
\label{thm:HK-WM-p-ge-2}\label{cor:universal-heat-WM}
Let $1\le p<\infty$. Under~\eqref{eq:LB}, there exists $C_{p,\varepsilon_{\rm L}}>0$ such that, for every $f\in\ell^p(V,m)$ and all integers $r\ge s\ge1$,
\begin{equation}\label{eq:heat-WM-p-ge-two}
 r^{-1/\min\{p,2\}}\cN_{p,r}^{\rm H}(f)
 \le C_{p,\varepsilon_{\rm L}}s^{-1/\min\{p,2\}}\cN_{p,s}^{\rm H}(f)\,.
\end{equation}
Without assuming~\eqref{eq:LB}, the heat kernel functionals satisfy $\mathrm{P}(p,\alpha)$ for every $\alpha\ge1/\min\{p,2\}$; equivalently, for every $f\in\ell^p(V,m)$,
\begin{equation}\label{eq:heat-p-ge-two-bridge}
 \|\nabla_pf\|_p\le\|f\|_{p,\alpha}\le C_p\|\nabla_pf\|_p\,.
\end{equation}
\end{theorem}

\begin{proof}
For integers $r\ge s\ge1$, put $n=\lceil r/s\rceil$. Then $r\le ns$ and $n\le2r/s$. By~\eqref{eq:heat-earlier-time-comparison}, followed by~\eqref{eq:heat-small-p-multiples} for $p\le2$ and~\eqref{eq:heat-p-ge-two-multiples} for $p\ge2$,
\begin{align*}
 \cN_{p,r}^{\rm H}(f)
 &\le2\varepsilon_{\rm L}^{-1/p}\cN_{p,ns}^{\rm H}(f)\\
 &\le C_{p,\varepsilon_{\rm L}}n^{1/\min\{p,2\}}\cN_{p,s}^{\rm H}(f)
 \le C_{p,\varepsilon_{\rm L}}(r/s)^{1/\min\{p,2\}}\cN_{p,s}^{\rm H}(f)\,.
\end{align*}
This proves~\eqref{eq:heat-WM-p-ge-two}. 

For~\eqref{eq:heat-p-ge-two-bridge}, we apply ~\eqref{eq:heat-small-p-multiples} for $p\le2$ and~\eqref{eq:heat-p-ge-two-multiples} for $p\ge2$ at $ns$ with $s=1$. Then 
\[\cN_{p,n}^{\rm H}(f)
 \le C_{p}n^{1/\min\{p,2\}}\cN_{p,1}^{\rm H}(f).
 \]
  Since $\cN_{p,1}^{\rm H}(f)=\|\nabla_pf\|_p$ and $n^{1/\min\{p,2\}-\alpha}\le1$, taking the supremum over $n\ge1$ gives the upper bound without~\eqref{eq:LB}, the lower bound follows from $n=1$.
\end{proof}

\subsection{Weak monotonicity from weak Bakry--\'Emery estimates}
\label{sec:WM-wBE}

We introduce the following discrete-time version of the weak Bakry--\'Emery condition from~\cite[Definition~3.1]{AlonsoEtAl2021III}.
\begin{definition}\label{def:wBE}
Let $d_w\ge2$ and $0<\kappa<d_w$. A weighted graph $(V,E)$ satisfies the \emph{weak Bakry-\'Emery condition}~\eqref{eq:wBE} if there exists $C>0$ such that
\begin{equation}\label{eq:wBE}
 |P^nf(x)-P^nf(y)|
 \le C\min\left\{1,
 \left(\frac{d(x,y)}{n^{1/d_w}}\right)^\kappa\right\}
 \|f\|_\infty\tag{${\rm wBE}(\kappa)$}
\end{equation}
for every bounded $f$, $n\ge1$, and $x,y\in V$.
\end{definition}

For $1\le p<\infty$, set
\begin{equation}\label{eq:beta-p}
 \beta_p(\kappa):=\frac1p+\left(1-\frac2p\right)\frac\kappa{d_w}\,.
\end{equation}

The next result proves that the weak Bakry-\'Emery estimate implies a weak monotonicity of the heat kernel functionals in the range $1 \le p \le 2$.

\begin{theorem}
\label{thm:HK-WM}\label{thm:WM-wBE}
Assume~\eqref{eq:LB}, $({\rm SG}_{d_f,d_w})$, and
$({\rm wBE}(\kappa))$, where $0<\kappa<d_w$.
For every $1\le p\le2$, $f\in\ell^p(V,m)$, and integers $r\ge s\ge1$,
\begin{equation}\label{eq:HK-WM}
 r^{-\beta_p(\kappa)}\cN_{p,r}^{\rm H}(f)
 \le C_p s^{-\beta_p(\kappa)}\cN_{p,s}^{\rm H}(f).
\end{equation}
Thus the heat kernel functionals satisfy
$({\rm WM}_{p,\beta_p(\kappa)})$.
\end{theorem}

We first compare the heat kernel $p_n(x,y)$ at two finite times. The lemma below is a discrete analogue of~\cite[Lemma~2.3]{AlonsoEtAl2021III}, see also~\cite[Section~2.2]{AlonsoRuizBaudoin2021}.
The constants in the comparison depend only on the parameters and constants in the hypotheses.

\begin{lemma}\label{lem:HK-dilations}
Assume~\eqref{eq:LB} and $({\rm SG}_{d_f,d_w})$. For every $\lambda>0$, there exist a constant $C>0$ and an integer $A\ge2$ such that
\begin{equation}\label{eq:weighted-HK-dilation}
 \left[1+\left(\frac{d(x,y)}{n^{1/d_w}}\right)^\lambda\right]
 p_n(x,y)\le Cp_{An}(x,y),
 \qquad \forall\; n\in\mathbb{N}.
\end{equation}
\end{lemma}

\begin{proof}
The assumption~\eqref{eq:LB} and the semigroup property give
\[
p_n(x,y)\ge p_{n-1}(x,y)\pi(y,y)\ge \varepsilon_{\rm L}p_{n-1}(x,y).
\]  Applying
\eqref{eq:SG-lower} at time $n-1$, we obtain
\begin{equation}\label{eq:exact-HK-lower}
 p_n(x,y)\ge \frac{\varepsilon_{\rm L}}{1+\varepsilon_{\rm L}}(p_{n-1}(x,y)+p_n(x,y))\ge cn^{-d_f/d_w}
 \exp\left[-C\left(\frac{d(x,y)}{n^{1/d_w}}\right)^{\frac{d_w}{d_w-1}}\right]
\end{equation}
for $n\ge2$ and $d(x,y)\le n-1$.

If $p_n(x,y)>0$, then $d(x,y)\le n$. Comparing~\eqref{eq:SG-upper} at time $n$
with~\eqref{eq:exact-HK-lower} at time $An$ gives
\[
 \frac{p_n(x,y)}{p_{An}(x,y)}
 \le CA^{d_f/d_w}
 \exp\left[-\bigl(c-C'A^{-1/(d_w-1)}\bigr) \left(\frac{d(x,y)}{n^{1/d_w}}\right)^{\frac{d_w}{d_w-1}}\right].
\]
We choose $A$ large enough so that $c-C'A^{-1/(d_w-1)}$ is positive, then~\eqref{eq:weighted-HK-dilation} follows since the exponential absorbs the polynomial.
\end{proof}

The next step is a telescopic identity that involves operators acting on functions $H\colon V\times V\to\mathbb{R}$ in weighted function spaces. To define these spaces, for each $j\ge1$, let $\nu_j$ be the measure on
$V\times V$ given by
\begin{equation}\label{eq:nu-j}
    \nu_j(\{(x,y)\}):=p_j(x,y)m(x)m(y).
\end{equation}
Define also $Df(x,y):=f(x)-f(y)$. With these notations, we have
\[\|Df\|_{L^p(\nu_j)}=\cN_{p,j}^{\rm H}(f).
\]

For functions $H$ on $V\times V$, we choose representatives that vanish outside the support of the relevant $\nu_j$.
Fix integers $n,s\ge1$. For $0\le k\le n-1$, define
\[
T_{n,s}^kH(u,v)
:=
\sum_{a,b\in V}
p_{ks}(u,a)p_{(n-1-k)s}(v,b)
H(a,b)m(a)m(b),
\]
and set
\begin{equation}\label{eq:bridge-T}
T_{n,s}H
:=
\sum_{k=0}^{n-1}T_{n,s}^kH.
\end{equation}
\begin{remark}

We observe that:
\begin{itemize}
 \item If $p_s(u,v)>0$, the sum in $T_{n,s}^kH$ has nonzero
 coefficients only for pairs $(a,b)$ with $p_{ns}(a,b)>0$.
 Thus the operators are well defined on equivalence classes.
 \item Jensen's inequality and the semigroup property show that
 each $T_{n,s}^k$ is a contraction from $L^q(\nu_{ns})$ to
 $L^q(\nu_s)$, $1\le q\le\infty$.  In particular,
 $\|T_{n,s}\|_{L^q(\nu_{ns})\to L^q(\nu_s)}\le n$.
 \item A function $H$ on $V\times V$ is antisymmetric if $H(y,x)=-H(x,y)$. Interchanging $u,v$ and replacing $k$ by $n-1-k$ shows  that $T_{n,s}$ commutes with interchanging the two variables. In particular, it preserves antisymmetry.

\end{itemize}
\end{remark}

The telescoping identity is given in the following lemma.
\begin{lemma}\label{lem:bridge-identity}
Let $n,s\ge1$ be integers, $1\le p<\infty$, and $q$ its conjugate exponent.  Let $T_{n,s}$ be defined by~\eqref{eq:bridge-T}.
For every $f\in\ell^p(V,m)$ and $H\in L^q(\nu_{ns})$, we have
\begin{align}\label{eq:bridge-identity}
 & \int_{V \times V} Df\, H \,d\nu_{ns}=\int_{V \times V} Df\, T_{n,s}H \,d\nu_{s}
\end{align}
and both sides converge absolutely.
\end{lemma}

\begin{proof}
For fixed $a,b\in V$, set
\[
 A_k(a,b):=\sum_{z\in V}
 p_{ks}(a,z)f(z)p_{(n-k)s}(z,b)m(z),
 \qquad 0\le k\le n.
\]
These sums are finite, since the graph is locally finite and the
walk has finite range.  For $0\le k\le n-1$, it follows from the semigroup identity that
\begin{align*}
 \sum_{u,v\in V}p_{ks}(a,u)p_s(u,v)p_{(n-1-k)s}(v,b)
      (f(u)-f(v))m(u)m(v)=A_k(a,b)-A_{k+1}(a,b).
\end{align*}
We sum over $k$ for the above equality, then
\begin{align*}
 &\sum_{k=0}^{n-1}\sum_{u,v\in V}p_{ks}(a,u)p_s(u,v)p_{(n-1-k)s}(v,b)
      (f(u)-f(v))m(u)m(v)\\
 &\qquad= A_0(a,b)-A_n(a,b)=(f(a)-f(b))p_{ns}(a,b).
\end{align*}
Multiplying both sides of the above equation by $H(a,b)m(a)m(b)$ and summing over $a,b\in V$, we obtain \eqref{eq:bridge-identity} by the symmetry $p_{ks}(a,u)=p_{ks}(u,a)$.

For completeness, H\"older's inequality and the bound
$\|T_{n,s}\|_{L^q(\nu_{ns})\to L^q(\nu_s)}\le n$ give
\begin{align*}
 &\int_{V \times V} |Df| \,|T_{n,s}H|\,d\nu_s
 \le \|Df\|_{L^p(\nu_s)}\|T_{n,s}H\|_{L^q(\nu_s)}
 \le2n\|f\|_p\|H\|_{L^q(\nu_{ns})}.
\end{align*}
This proves absolute convergence of the expanded sum and justifies
the rearrangements, including when $p=1$ and $q=\infty$.
The left-hand side of~\eqref{eq:bridge-identity} is absolutely
convergent by the same H\"older estimate with $\nu_{ns}$.
\end{proof}

In the following three lemmas, we prove an $L^2$ bound and a weighted $L^\infty$ bound for $T_{n,s}$, and then deduce, by interpolation, a weighted $L^q$ bound for $2\le q\le\infty$.
\begin{lemma}\label{lem:L2-bridge}
For every antisymmetric $H\in L^2(\nu_{ns})$,
\begin{equation}\label{eq:L2-bridge}
 \|T_{n,s}H\|_{L^2(\nu_s)}
 \le\sqrt n\,\|H\|_{L^2(\nu_{ns})}.
\end{equation}
\end{lemma}

\begin{proof}
Let $h\in L^2(\nu_s)$ be antisymmetric, and set
\[
    h_j:=T_{j,s}^*h\in L^2(\nu_{js}),\qquad j\ge1.
\]
The adjoint is well defined by recalling that
$\|T_{j,s}\|_{L^2(\nu_{js})\to L^2(\nu_s)}\le j$. Since $T_{j,s}$ commutes with
interchanging the two variables, its adjoint preserves antisymmetry,
so each $h_j$ is antisymmetric.
Expanding~\eqref{eq:bridge-T} gives $h_1=h$ and
\begin{equation}\label{eq:adjoint-kernel-recursion}
    \pi_{(j+1)s}(x,y)h_{j+1}(x,y)
    =\sum_{z\in V}\pi_{js}(x,z)\pi_s(z,y)
      \bigl(h_j(x,z)+h(z,y)\bigr).
\end{equation}
We choose representatives that vanish where the corresponding kernel vanishes.
Define
\[
    a(x):=\sum_{y\in V}\pi_s(x,y)h(x,y),
    \qquad
    a_j(x):=\sum_{y\in V}\pi_{js}(x,y)h_j(x,y).
\]
Jensen's inequality implies that
$a,a_j\in\ell^2(V,m)$. Summing
\eqref{eq:adjoint-kernel-recursion} over $y$ gives
\[
    a_{j+1}=a_j+P^{js}a.
\]
Since $a_1=a$, induction gives $a_j=\sum_{k=0}^{j-1}P^{ks}a$.
Moreover, reversibility and antisymmetry yield
\begin{equation}\label{eq:a_j}
    \sum_{x\in V}m(x)\pi_{js}(x,z)h_j(x,z)
    =-m(z)a_j(z).
\end{equation}

In view of \eqref{eq:adjoint-kernel-recursion} and the fact $\sum_{z\in V}\pi_{js}(x,z)\pi_s(z,y) =\pi_{(j+1)s}(x,y)$, we have from Jensen's inequality that
\begin{align*}
   \|h_{j+1}\|_{L^2(\nu_{(j+1)s})}^2&=\sum_{x,y\in V} \Bigg(\sum_{z\in V}
      \bigl(h_j(x,z)+h(z,y)\bigr)\frac{\pi_{js}(x,z)\pi_s(z,y)}{\pi_{(j+1)s}(x,y)}\Bigg)^2 \pi_{(j+1)s}(x,y)m(x)
    \\ &\le   \sum_{x,z,y\in V}
    \bigl(h_j(x,z)+h(z,y)\bigr)^2\pi_{js}(x,z)\pi_s(z,y) m(x).
\end{align*}
Expanding the square and using conservativeness together with the identity\eqref{eq:a_j} give
\begin{equation}\label{eq:adjoint-energy-recursion}
    \|h_{j+1}\|_{L^2(\nu_{(j+1)s})}^2
    \le
    \|h_j\|_{L^2(\nu_{js})}^2
    +\|h\|_{L^2(\nu_s)}^2
    -2\langle a_j,a\rangle_2.
\end{equation}

We observe that $\langle a_j,a\rangle_2\ge0$.
Indeed, $P^s$ is a self-adjoint contraction, so consecutive even and odd powers satisfy
\[
 \langle a,(P^{2rs}+P^{(2r+1)s})a\rangle_2
 =\langle P^{rs}a,(I+P^s)P^{rs}a\rangle_2\ge0, \quad \forall r\ge 0.
\]
If an unpaired even power of $P^s$ remains, it has nonnegative quadratic form.

Induction in~\eqref{eq:adjoint-energy-recursion} therefore yields
\[
    \|T_{n,s}^*h\|_{L^2(\nu_{ns})}^2
    \le n\|h\|_{L^2(\nu_s)}^2.
\]
For antisymmetric $H$, the function $T_{n,s}H$ is antisymmetric.
We thus prove \eqref{eq:L2-bridge} by duality on the antisymmetric subspaces.

\end{proof}

Next we will prove a weighted $L^\infty$ estimate for $T_{n,s}$, where the spatial weight is given by
\begin{equation}\label{eq:bridge-weight}
 w_s(u,v):=1+\left(\frac{d(u,v)}{s^{1/d_w}}\right)^\kappa.
\end{equation}

\begin{lemma}\label{lem:Linfinity-bridge}
Assume $({\rm wBE}(\kappa))$. For every antisymmetric $H\in L^\infty(\nu_{ns})$,
\begin{equation}\label{eq:Linfinity-bridge}
 \|w_s^{-1}T_{n,s}H\|_{L^\infty(\nu_s)}
 \le Cn^{1-\kappa/d_w}\|H\|_{L^\infty(\nu_{ns})}.
\end{equation}
\end{lemma}

\begin{proof}
We choose an antisymmetric representative of $H$ which vanishes where
$p_{ns}=0$ and satisfies
$|H(a,b)|\le\|H\|_{L^\infty(\nu_{ns})}$ for every $a,b\in V$.
Without loss of generality, we assume $n\ge 2$. In fact, the estimate is trivial for $n=1$ since $T_{1,s}H=H$.
For any fixed $v\in V$, set
\[
 F_{k,v}(a):=\sum_b p_{(n-1-k)s}(v,b)H(a,b)m(b).
\]
Then conservativeness gives $\|F_{k,v}\|_\infty\le\|H\|_{L^\infty(\nu_{ns})}$, and
$T_{n,s}H(u,v)=\sum_{k=0}^{n-1}P^{ks}F_{k,v}(u)$.
We further observe from~\eqref{eq:bridge-T} that $T_{n,s}H(v,v)=0$ since $H$ is antisymmetric.
Consequently,
\[
 T_{n,s}H(u,v)
 =\sum_{k=0}^{n-1}
   \bigl(P^{ks}F_{k,v}(u)-P^{ks}F_{k,v}(v)\bigr).
\]
For $k=0$, the absolute value of the summand is at most
$2\|H\|_\infty$.  For $k\ge1$, we apply $({\rm wBE}(\kappa))$ at
time $ks$ to the bounded function $F_{k,v}$. Hence
\[
 |P^{ks}F_{k,v}(u)-P^{ks}F_{k,v}(v)|
 \le C\left(\frac{d(u,v)}{s^{1/d_w}}\right)^\kappa
       k^{-\kappa/d_w}\|H\|_{L^\infty(\nu_{ns})}.
\]
Summing these estimates, we obtain
\begin{align*}
 |T_{n,s}H(u,v)|
 &\le\left[2+C\left(\frac{d(u,v)}{s^{1/d_w}}\right)^\kappa
          \sum_{k=1}^{n-1}k^{-\kappa/d_w}\right]\|H\|_{L^\infty(\nu_{ns})}\\
 &\le Cn^{1-\kappa/d_w}w_s(u,v)\|H\|_{L^\infty(\nu_{ns})}.
\end{align*}
For the last inequality, we have used $n^{1-\kappa/d_w}\ge1$ and
\(\sum_{k=1}^{n-1}k^{-\kappa/d_w}
 \le C n^{1-\kappa/d_w}\) for $n\ge2$.
 Dividing by $w_s(u,v)$ and taking the essential
supremum with respect to $\nu_s$ completes the proof.
\end{proof}

\begin{lemma}\label{lem:weighted-interpolation}
Assume $({\rm wBE}(\kappa))$. For $2\le q\le\infty$ and antisymmetric $H\in L^q(\nu_{ns})$,
\begin{equation}\label{eq:weighted-interpolation}
 \|w_s^{2/q-1}T_{n,s}H\|_{L^q(\nu_s)}
 \le C_q n^{\frac1q+(1-\frac2q)(1-\frac\kappa{d_w})}
       \|H\|_{L^q(\nu_{ns})}.
\end{equation}
\end{lemma}

\begin{proof}
Since the $L^2$ and weighted $L^{\infty}$ estimates for $T_{n,s}$ are established only for antisymmetric functions, we introduce the antisymmetrization operator $\mathsf A H(u,v)=(H(u,v)-H(v,u))/2$. By symmetry of $\nu_{ns}$, this is a contraction on every
$L^q(\nu_{ns})$, and $\mathsf A H=H$ whenever $H$ is antisymmetric.

Consider the operator $w_s^{-1}T_{n,s}\mathsf A$. It follows from Lemma \ref{lem:L2-bridge} that
\begin{align*}
 \|w_s^{-1}T_{n,s}\mathsf A H\|_{L^2(w_s^2\nu_s)}^2
 &=\int |w_s^{-1}T_{n,s}\mathsf A H|^2w_s^2\,d\nu_s\\
 &=\|T_{n,s}\mathsf A H\|_{L^2(\nu_s)}^2
 \le n\|H\|_{L^2(\nu_{ns})}^2.
\end{align*}
For the $L^\infty$ norm, since $w_s^2\nu_s$ and $\nu_s$ have the same null sets,  Lemma~\ref{lem:Linfinity-bridge} gives
\[
\|w_s^{-1}T_{n,s}\mathsf A H\|_{L^{\infty}(w_s^2\nu_s)}\le Cn^{1-\kappa/d_w}\|\mathsf AH\|_{L^{\infty}(\nu_{ns})} \le Cn^{1-\kappa/d_w}\|H\|_{L^{\infty}(\nu_{ns})}.
\]

Applying the Riesz--Thorin theorem to the
operator $w_s^{-1}T_{n,s}\mathsf A$ (see \cite[Chapter~V, Section~1]{SteinWeiss1971}), we obtain
\[
 \|w_s^{2/q-1}T_{n,s}\mathsf AH\|_{L^q(\nu_s)}=\|w_s^{-1}T_{n,s}\mathsf AH\|_{L^q(w_s^2\nu_s)}
 \le C_q n^{\frac1q+(1-\frac2q)(1-\frac\kappa{d_w})}
       \|H\|_{L^q(\nu_{ns})}.
\]
The assertion follows from the fact that $\mathsf A H=H$ for antisymmetric $H$.
\end{proof}

\begin{proof}[Proof of Theorem~\ref{thm:HK-WM}]
Let $q$ be conjugate to $p$, so that $2\le q\le\infty$.
Choose $A$ in Lemma~\ref{lem:HK-dilations} with $\lambda=\kappa$.
Since $p(1-2/q)=2-p\in[0,1]$ and $w_s\ge1$, the kernel comparison in \eqref{eq:weighted-HK-dilation}
and time subadditivity in Lemma \ref{lem:time-subadditivity} give
\begin{align}
 \|w_s^{1-2/q}Df\|_{L^p(\nu_s)}^p
 &=\int |Df|^p w_s^{2-p}\,d\nu_s\notag
 \le\int |Df|^p w_s\,d\nu_s
 \\ &\le C\bigl(\cN_{p,As}^{\rm H}(f)\bigr)^p
 \le CA^p\bigl(\cN_{p,s}^{\rm H}(f)\bigr)^p.
 \label{eq:weighted-D-bound}
\end{align}
For antisymmetric $H\in L^q(\nu_{ns})$, the telescoping identity \eqref{eq:bridge-identity},
H\"older's inequality, and Lemma~\ref{lem:weighted-interpolation} yield
\begin{align*}
 \left|\int_{V \times V} Df\,H\,d\nu_{ns}\right|=\left|\int_{V \times V}  Df\,T_{n,s}Hd\nu_s\right|
 &\le\|w_s^{1-2/q}Df\|_{L^p(\nu_s)}
       \|w_s^{2/q-1}T_{n,s}H\|_{L^q(\nu_s)}\\
 &\le C_p n^{\beta_p(\kappa)}
       \cN_{p,s}^{\rm H}(f)\|H\|_{L^q(\nu_{ns})},
\end{align*}
where the last line follows from
\(
 \frac1q+\left(1-\frac2q\right)
              \left(1-\frac\kappa{d_w}\right)
 =\frac1p+\left(1-\frac2p\right)\frac\kappa{d_w}
 =\beta_p(\kappa).
\)

Since $Df$ is antisymmetric, its $L^p(\nu_{ns})$ norm equals the supremum of $|\int DfH\,d\nu_{ns}|$ over antisymmetric $H$ with $\|H\|_{L^q(\nu_{ns})}\le1$.
This remains true for $p=1$ by taking $H=\operatorname{sgn}(Df)$.
Hence
\begin{equation}\label{eq:exact-multiples}
 \cN_{p,ns}^{\rm H}(f)
 \le C_p n^{\beta_p(\kappa)}\cN_{p,s}^{\rm H}(f).
\end{equation}

For general $r\ge s$, take $n=\lceil r/s\rceil$. Then $r\le ns$ and $n\le2r/s$. By~\eqref{eq:heat-earlier-time-comparison} and~\eqref{eq:exact-multiples},
\begin{equation*}
 \cN_{p,r}^{\rm H}(f)
 \le2\varepsilon_{\rm L}^{-1/p}\cN_{p,ns}^{\rm H}(f)
 \le C_p(r/s)^{\beta_p(\kappa)}\cN_{p,s}^{\rm H}(f)\,.
\end{equation*}
Since $\beta_p(\kappa)>0$, this proves~\eqref{eq:HK-WM}.
\end{proof}

For the property $\mathrm{P}(p,\beta_p(\kappa))$, we only need a weaker form of Bakry-\'Emery.

\begin{theorem}\label{thm:HK-P-wBE}
Let $d_w\ge2$ and $0<\kappa<d_w$. Assume that there exists $C>0$ such that
\begin{equation}\label{eq:edge-smoothing}
 |P^ng(x)-P^ng(y)|\le Cn^{-\kappa/d_w}\|g\|_\infty
\end{equation}
for every bounded $g$,  $n\ge1$, and $x\sim y$. Then, for $1\le p\le2$, the heat kernel functionals satisfy $\mathrm{P}(p,\beta_p(\kappa))$.
 In particular, this conclusion holds under $({\rm wBE}(\kappa))$, without assuming~\eqref{eq:LB} or $({\rm SG}_{d_f,d_w})$.
\end{theorem}

\begin{proof}
Consider the operator $T_{n,1}$ defined in~\eqref{eq:bridge-T}. For antisymmetric $H\in L^\infty(\nu_n)$, the cancellation used in the proof of Lemma~\ref{lem:Linfinity-bridge} and~\eqref{eq:edge-smoothing} give
\begin{align*}
 \|T_{n,1}H\|_{L^\infty(\nu_1)}
 &\le\left(2+C\sum_{k=1}^{n-1}k^{-\kappa/d_w}\right)\|H\|_{L^\infty(\nu_n)}\\
 &\le Cn^{1-\kappa/d_w}\|H\|_{L^\infty(\nu_n)}\,.
\end{align*}
Indeed, the off-diagonal support of $\nu_1$ consists of adjacent pairs, and $T_{n,1}H$ vanishes on the diagonal. Lemma~\ref{lem:L2-bridge} gives the corresponding $L^2$ bound with constant $\sqrt n$. Apply the Riesz--Thorin theorem to $T_{n,1}\mathsf A$, where $\mathsf A$ is the antisymmetrization operator from Lemma~\ref{lem:weighted-interpolation}. For $2\le q\le\infty$ and antisymmetric $H$,
\begin{equation*}
 \|T_{n,1}H\|_{L^q(\nu_1)}
 \le C_qn^{\frac1q+(1-\frac2q)(1-\frac\kappa{d_w})}
 \|H\|_{L^q(\nu_n)}\,.
\end{equation*}
Let $q$ be conjugate to $p$ and let $f\in\ell^p(V,m)$. The telescoping identity~\eqref{eq:bridge-identity} and H\"older's inequality yield
\begin{equation*}
 \left|\int_{V \times V} Df\,H\,d\nu_n\right|
 \le C_pn^{\beta_p(\kappa)}\|\nabla_pf\|_p\|H\|_{L^q(\nu_n)}\,.
\end{equation*}
For $1<p\le 2$, taking the supremum over antisymmetric $H$ with norm at most one gives
\begin{equation*}
 \cN_{p,n}^{\rm H}(f)\le C_pn^{\beta_p(\kappa)}\|\nabla_pf\|_p\,,
\end{equation*}
For $p=1$, the result follows by taking $H=\operatorname{sgn}(Df)$.
\end{proof}

\section{Weak Bakry-\'Emery estimates}
\label{sec:wBE}
In this section, we establish weak Bakry-\'Emery estimates under sub-Gaussian heat kernel bounds. We first show that controlled weights imply $({\rm wBE}(\kappa))$ for some positive exponent $\kappa$, using parabolic Harnack regularity. We then obtain the explicit exponent $\kappa=d_w-d_f$ in the strongly recurrent case.

\subsection{Sub-Gaussian graphs}\label{sec:wBE-general}

We  show that controlled weights and sub-Gaussian estimates actually always imply $({\rm wBE}(\kappa))$ for some positive exponent. The argument uses parabolic Harnack regularity  but does not identify the optimal value of $\kappa$.

\begin{proposition}\label{prop:wBE-sub-gaussian}
Assume controlled weights and $({\rm SG}_{d_f,d_w})$. Then $({\rm wBE}(\kappa))$ holds for some $\kappa\in(0,1]$, depending only on the parameters and constants in the hypotheses.
\end{proposition}
\begin{proof}
The main idea, which is classical, is to apply the parabolic Harnack inequality to $u(k,z)=P^kf(z)$ and deduce from it the needed oscillation estimates. All time variables are integers.

\emph{Step 1: Parabolic Harnack inequality.}
We first recall explicitly the parabolic Harnack inequality. By Proposition~\ref{prop:volume-growth}, condition~\eqref{eq:volume-growth} holds. Hence $({\rm SG}_{d_f,d_w})$ gives the volume-normalized sub-Gaussian estimates of~\cite{GrigoryanTelcs2002}. Since the weights are controlled, $\pi(x,y)\ge\min\{\varepsilon_0,\varepsilon_{\rm L}\}$ whenever $d(x,y)\le1$, which is condition $(p_0)$ of~\cite{GrigoryanTelcs2002}; hence these estimates imply the parabolic Harnack inequality by~\cite[Theorem~3.1]{GrigoryanTelcs2002}.

More precisely, for $x\in V$, $R>0$ and $t\in\mathbb Z$, set $N_R:=\lceil R^{d_w}\rceil$ and
\begin{equation*}
 Q_R(x,t):=\bigl([t-4N_R,t]\cap\mathbb Z\bigr)\times B(x,2R+1)\,.
\end{equation*}
We call $u\colon Q_R(x,t)\to\R$ \emph{caloric in $Q_R(x,t)$} if, for every integer $t-4N_R\le k<t$ and every $z\in B(x,2R)$,
\begin{equation}\label{eq:wBE-heat-equation}
 u(k+1,z)=\sum_{w\in B(x,2R+1)}\pi(z,w)u(k,w)\,.
\end{equation}
Note that $z\in B(x,2R)$ and $\pi(z,w)>0$ imply $w\in B(x,2R+1)$, so the sum contains every nonzero transition from~$z$. We also note that the set of caloric functions is a linear space, contains constant functions and is closed by restriction to $Q_\rho(x,t)$, $0<\rho\le R$, since $P 1=1$.

For $R\ge8$, we have $N_R\simeq R^{d_w}$ and $N_R\ge R^2\ge2R$. Applying~\cite[Definition~2.2]{GrigoryanTelcs2002}, whose balls are open, with radius $\lfloor R\rfloor+\tfrac12$, $N=N_R$ and time origin $t-4N_R$, we obtain that every nonnegative caloric function~$u$ in $Q_R(x,t)$ satisfies
\begin{equation}\label{eq:wBE-PHI-GT}
 \sup_{\substack{t-3N_R\le k<t-2N_R\\z\in B(x,R)}}u(k,z)
 \le C_{\rm H}\inf_{\substack{t-N_R\le k<t\\z\in B(x,R)}}
 \bigl(u(k,z)+u(k+1,z)\bigr)\,,
\end{equation}
where $C_{\rm H}\ge1$ depends only on the parameters and constants in the hypotheses. By~\eqref{eq:LB} and~\eqref{eq:wBE-heat-equation}, we have
\begin{equation*}
 u(k+1,z)\ge\varepsilon_{\rm L}u(k,z),
\end{equation*}
and therefore
\begin{equation*}
 u(k,z)+u(k+1,z)\le(1+\varepsilon_{\rm L}^{-1})u(k+1,z)\,.
\end{equation*}
Thus, with $C_1:=(1+\varepsilon_{\rm L}^{-1})C_{\rm H}$,
\begin{equation}\label{eq:wBE-parabolic-Harnack}
 \sup_{Q_R^-(x,t)}u\le C_1\inf_{Q_R^+(x,t)}u\,,
\end{equation}
where
\begin{align*}
 Q_R^-(x,t)&:=\bigl([t-3N_R,t-2N_R)\cap\mathbb Z\bigr)\times B(x,R),\\
 Q_R^+(x,t)&:=\bigl([t-N_R+1,t]\cap\mathbb Z\bigr)\times B(x,R)\,.
\end{align*}

\emph{Step 2: Oscillation decay.}
We adapt the argument of~\cite[Proposition~3.2]{BarlowHambly2009} to the time scale~$R^{d_w}$. Set
\begin{equation*}
 q:=1-\frac1{2C_1}\in\Bigl[\frac12,1\Bigr),
 \qquad
 \kappa:=\frac{\log(1/q)}{\log4}\in\Bigl(0,\frac12\Bigr],
\end{equation*}
so that $q=4^{-\kappa}$, and write $\operatorname{osc}_Qu:=\sup_Qu-\inf_Qu$. We claim that every real-valued caloric function~$u$ in $Q_R(x,t)$, with $R\ge8$, satisfies
\begin{equation}\label{eq:wBE-oscillation-decay}
 \operatorname{osc}_{Q_{R/4}(x,t)}u\le q\operatorname{osc}_{Q_R(x,t)}u\,.
\end{equation}
First, $Q_{R/4}(x,t)\subset Q_R^+(x,t)$: the spatial inclusion follows from $R/2+1\le R$, while $d_w\ge2$ and $R\ge8$ give
\begin{equation*}
 4N_{R/4}\le4^{1-d_w}R^{d_w}+4\le\tfrac14R^{d_w}+4\le R^{d_w}-1\le N_R-1\,.
\end{equation*}
If $\operatorname{osc}_{Q_R(x,t)}u=0$, the claim is immediate. Otherwise, subtracting the infimum and dividing by the oscillation, we may assume that $\inf_{Q_R(x,t)}u=0$ and $\sup_{Q_R(x,t)}u=1$. Both $u$ and $1-u$ are then nonnegative and caloric. Replacing $u$ by $1-u$ if necessary, we may also assume that $\sup_{Q_R^-(x,t)}u\ge1/2$. By~\eqref{eq:wBE-parabolic-Harnack},
\begin{equation*}
 \inf_{Q_{R/4}(x,t)}u\ge\inf_{Q_R^+(x,t)}u\ge\frac1{2C_1}\,.
\end{equation*}
Since $\sup_{Q_{R/4}(x,t)}u\le1$, this proves~\eqref{eq:wBE-oscillation-decay}. Iterating at the radii $R,R/4,\ldots,4^{1-j}R$ gives
\begin{equation}\label{eq:wBE-iterated-oscillation}
 \operatorname{osc}_{Q_{4^{-j}R}(x,t)}u\le4^{-\kappa j}\operatorname{osc}_{Q_R(x,t)}u
\end{equation}
for every integer $j\ge1$ such that $4^{1-j}R\ge8$; the same estimate is immediate for $j=0$.

\emph{Step 3: H\"older estimate for $P^nf$.}
Let $f$ be bounded, $n\ge1$, and $x\ne y$. The bound $|P^nf(x)-P^nf(y)|\le2\|f\|_\infty$ proves~\eqref{eq:wBE} whenever $n<8^{d_w+1}$ or $d(x,y)>(n/8)^{1/d_w}$, since in these cases
\begin{equation*}
 \min\left\{1,\left(\frac{d(x,y)}{n^{1/d_w}}\right)^\kappa\right\}
 \ge8^{-\kappa(1+1/d_w)}\,.
\end{equation*}
We may therefore assume that $n\ge8^{d_w+1}$ and $1\le d(x,y)\le R:=(n/8)^{1/d_w}$. Then $R\ge8$ and
\begin{equation*}
 4N_R=4\lceil n/8\rceil\le n/2+4\le n\,.
\end{equation*}
Consequently, $u(k,z):=P^kf(z)$ is defined and caloric throughout $Q_R(x,n)$, and contractivity gives $\operatorname{osc}_{Q_R(x,n)}u\le2\|f\|_\infty$.

Set $r:=\max\{8,d(x,y)\}\le R$ and choose the largest integer $j\ge0$ for which $4^{-j}R\ge r$. Then $r\le4^{-j}R<4r$, and~\eqref{eq:wBE-iterated-oscillation} applies. Both $(n,x)$ and $(n,y)$ lie in $Q_{4^{-j}R}(x,n)$, so
\begin{align}\label{eq:wBE-main-case}
 |P^nf(x)-P^nf(y)|
 &\le\operatorname{osc}_{Q_{4^{-j}R}(x,n)}u
 \le2\cdot4^{-\kappa j}\|f\|_\infty\notag\\
 &\le2\left(\frac{4r}{R}\right)^\kappa\|f\|_\infty
 \le C\left(\frac{d(x,y)}{n^{1/d_w}}\right)^\kappa\|f\|_\infty\,,
\end{align}
where we used $r\le8d(x,y)$. Since $d(x,y)\le R<n^{1/d_w}$, this is~\eqref{eq:wBE} in the remaining case.
\end{proof}

\begin{remark}
Condition~\eqref{eq:LB} cannot simply be dropped to prove the weak Bakry-\'Emery estimate. Indeed, the simple random walk on~$\mathbb Z$ satisfies $({\rm SG}_{1,2})$, with the lower bound formulated for $p_n+p_{n+1}$, but $f(j)=(-1)^j$ satisfies $P^nf=(-1)^nf$. Thus $|P^nf(0)-P^nf(1)|=2$ for every~$n$, which excludes $({\rm wBE}(\kappa))$ for every $\kappa>0$.
\end{remark}

\subsection{Strongly recurrent graphs}\label{sec:wBE-recurrent}

We now prove the discrete counterpart of~\cite[Theorem~3.7]{AlonsoEtAl2021III}. A graph is strongly recurrent here if it satisfies $({\rm SG}_{d_f,d_w})$ with
\begin{equation}\label{eq:strong-recurrence}
 1\le d_f<d_w.
\end{equation}

\begin{theorem}\label{thm:wBE-strong-recurrence}
Assume~\eqref{eq:LB} and $({\rm SG}_{d_f,d_w})$ with $1\le d_f<d_w$. Then $({\rm wBE}(d_w-d_f))$ holds. Consequently, for every $1\le p\le2$, the heat kernel functionals satisfy $({\rm WM}_{p,\beta_p(d_w-d_f)})$, where
\begin{equation}\label{eq:strong-recurrent-beta}
 \beta_p(d_w-d_f)
 =\frac1p+\left(1-\frac2p\right)\frac{d_w-d_f}{d_w}\,.
\end{equation}
\end{theorem}

The proof combines a local Poisson estimate with the temporal bound $\|LP^n\|_{\ell^\infty\to\ell^\infty}\le Cn^{-1}$. We first establish these two estimates.

To prove~\eqref{eq:wBE} we will make use of the \emph{effective resistance distance} associated with the weighted graph. To introduce the latter, we first extend the definition of the Dirichlet form $(\cE,\cF)$ in~\eqref{eq:def_graph_DF} to $(\cE,\cF_e)$, where
\begin{equation}\label{eq:finite-energy-class}
 \cF_{\rm e}:=  \Big\{u\colon V\to\R\quad\text{ s.t. }\quad \cE(u,u)=\frac12\sum_{x,y\in V}c_{xy}|u(x)-u(y)|^2<\infty\Big\}.
\end{equation}
The effective resistance distance between two disjoint nonempty sets $A,B\subset V$ is defined as
\begin{equation}\label{eq:resistance-definition}
 \cR(A,B)^{-1} :=\inf\{\cE(v,v):v\in\cF_{\rm e},\ v|_A=1,\ v|_B=0\}
\end{equation}
and for any two distinct vertices $x,y\in V$,
\begin{equation}\label{eq:point-resistance-dual}
 \cR(x,y) =\sup_{\substack{v\in\cF_{\rm e}\\ \cE(v,v)>0}} \frac{|v(x)-v(y)|^2}{\cE(v,v)}.
\end{equation}
\begin{proposition}
\label{prop:resistance-growth}
Assume~\eqref{eq:LB} and $({\rm SG}_{d_f,d_w})$ with  $1\le d_f<d_w$. Then, there are $c,C>0$ such that
    \begin{equation}\label{eq:point-resistance-growth}
     c d(x,y)^{d_w-d_f}\le\cR(x,y)\le C d(x,y)^{d_w-d_f},
     \qquad x,y\in V,\quad x\ne y,                            
    \end{equation}
 and
 \begin{equation}\label{eq:annular-resistance}
     \cR(x,B(x,R)^c)\simeq R^{d_w-d_f},
    \hspace{7em} x\in V,\quad R\ge1.                               
 \end{equation}
\end{proposition}
\begin{proof}
Proposition~\ref{prop:volume-growth} gives~\eqref{eq:volume-growth} and uniform two-sided bounds on $m(x)$.
Thus the graph satisfies volume doubling, the volume growth condition $({\rm VG}(d_w^-))$ and the heat kernel estimates $({\rm HK}(d_w))$ of~\cite{BCK2005}; the results of~\cite{BCK2005} used below do not require their condition $(p_0)$.
By~\cite[Proposition~4.2(1)]{BCK2005}, it satisfies the corresponding scaled Poincar\'e inequality.
Since $d_f<d_w$,~\cite[Lemma~2.3(b)]{BCK2005} therefore gives
\begin{equation*}
 \cR(x,y)\le C\frac{d(x,y)^{d_w}}{\Vol(x,d(x,y))}
 \le Cd(x,y)^{d_w-d_f},\qquad x\ne y.
\end{equation*}
Moreover,~\cite[Proposition~4.2(2)]{BCK2005} yields
\begin{equation*}
 \cR(x,B(x,R)^c)\ge c\frac{R^{d_w}}{\Vol(x,R)}
 \ge cR^{d_w-d_f},\qquad R\ge1.
\end{equation*}
If $d(x,y)\ge2$, monotonicity of effective resistance gives
\begin{equation*}
 \cR(x,y)\ge\cR\bigl(x,B(x,d(x,y)/2)^c\bigr)
 \ge cd(x,y)^{d_w-d_f}.
\end{equation*}
For $d(x,y)=1$, test~\eqref{eq:point-resistance-dual} with $\1_{\{x\}}$ and use
\begin{equation*}
 \cE(\1_{\{x\}},\1_{\{x\}})\le m(x)\le C.
\end{equation*}
Finally, connectedness and infinitude give a vertex $z\notin B(x,R)$ with $d(x,z)\le R+1$. Hence
\begin{equation*}
 \cR(x,B(x,R)^c)\le\cR(x,z)
 \le C(R+1)^{d_w-d_f}\le CR^{d_w-d_f}.
\end{equation*}
\end{proof}

\begin{proposition}\label{prop:temporal-regularity}
Assume~\eqref{eq:LB} and $({\rm SG}_{d_f,d_w})$. There exists $C>0$ such that
\begin{equation}\label{eq:infinity-Ritt}
 \|LP^n\|_{\infty\to\infty}\le\frac Cn,
 \qquad n\ge1\,.
\end{equation}
\end{proposition}

\begin{proof}
By Proposition~\ref{prop:volume-growth}, condition~\eqref{eq:volume-growth} holds. Together with~\eqref{eq:LB} and~\eqref{eq:SG-upper}, this verifies the hypotheses (LB), $(\mathrm D_\beta)$ and $(\mathrm{UE}_\beta)$ of~\cite[Proposition~2.6]{Feneuil2015SG} with $\beta\equiv d_w$ and quasidistance $d^{d_w}$. That proposition, obtained by the method of~\cite[Theorem~1.1]{Dungey2006}, applied with derivative order~$1$ at time~$n+1$, gives, for every $M>0$,
\begin{equation}\label{eq:temporal-kernel}
 |p_{n+1}(z,w)-p_n(z,w)|
 \le C_M n^{-1-d_f/d_w}
 \left(1+\frac{d(z,w)^{d_w}}n\right)^{-M},
 \qquad n\ge1\,.
\end{equation}
Choose $M>d_f/d_w$ and put $r=n^{1/d_w}$. Summing over $B(z,2r)$ and the annuli $B(z,2^{j+1}r)\setminus B(z,2^jr)$, $j\ge1$, yields
\begin{equation*}
 \sum_{w\in V}|p_{n+1}(z,w)-p_n(z,w)|m(w)
 \le\frac Cn\left(1+\sum_{j\ge1}2^{j(d_f-Md_w)}\right)
 \le\frac Cn\,.
\end{equation*}
Since $LP^n=P^n-P^{n+1}$, this proves~\eqref{eq:infinity-Ritt}.
\end{proof}

The key to proving the weak Bakry--\'Emery estimate is a local Poisson estimate
based on local Green's functions, defined as follows:
For a finite, nonempty, and proper subset $D\subset V$, let
\begin{equation*}
    P_D f:=\1_D P(\1_D f),
\end{equation*}
and denote by $p_k^D(z,w)$, $z,w\in D$, the kernel of $P_D^k$
restricted to $D$. Further, let
\begin{equation}\label{eq:killed-Green}
    G_D(z,w):=\sum_{k=0}^{\infty}p_k^D(z,w),
    \qquad z,w\in D,
\end{equation}
and note that this is a convergent series, see e.g.~\cite[Theorem~1.31]{Barlow2017}. For each $w\in D$, $G_D(\cdot,w)$ may be extended to $V$ by setting $G_D(z,w)=0$ for $z\in D^c$.

\medskip

Let
\begin{equation*}
 \tau_D:=\inf\{k\ge0:X_k\notin D\},
 \qquad T_w:=\inf\{k\ge0:X_k=w\}.
\end{equation*}
For $w\in D$ and $z\in V$, the standard Green identities are
\begin{equation}\label{eq:Green-resistance}
 G_D(w,w)=\cR(w,D^c)
\end{equation}
and
\begin{equation}\label{eq:Green-hitting}
 \mathbb P^z(T_w<\tau_D)=\frac{G_D(z,w)}{G_D(w,w)};
\end{equation}
see~\cite[Theorems~1.31, 1.34(c), and~2.8]{Barlow2017}. We note that Gao and Yang~\cite[Proposition~4.3 and Section~5]{GaoYang2026} obtain weak Bakry--\'Emery estimates from Poisson regularity on strongly local Dirichlet spaces. Here we use a discrete Green representation and the estimate of Proposition~\ref{prop:temporal-regularity}.

\begin{lemma}
\label{lem:local-Poisson}
Assume~\eqref{eq:LB} and $({\rm SG}_{d_f,d_w})$ with $1\le d_f<d_w$. Let $R\ge1$ and $x,y\in V$ with $r=d(x,y)\le R$. For any $u\in \ell^\infty(V)$, 
\begin{equation}\label{eq:local-Poisson}
 |u(x)-u(y)|
 \le Cr^{d_w-d_f}\left(R^{d_f-d_w}\|u\|_\infty
 +R^{d_f}\|Lu\|_\infty\right).
\end{equation}
\end{lemma}

\begin{proof}
The assertion is immediate if $x=y$. Otherwise, set $D=B(y,4R)$ and identify $D^c$ with a single vertex $\partial$, as in~\cite[Definition~1.5(2)]{Barlow2017}. Denote the effective resistance on this finite network by $\widehat{\cR}_D$. The variational definition gives, for $a,b\in D$,
\begin{equation*}
 \widehat{\cR}_D(a,\partial)=\cR(a,D^c),
 \qquad
 \widehat{\cR}_D(a,b)\le\cR(a,b)\,.
\end{equation*}
By~\cite[Proposition~2.55]{Barlow2017},
\begin{equation*}
 2G_D(a,w)=\widehat{\cR}_D(a,\partial)
 +\widehat{\cR}_D(w,\partial)-\widehat{\cR}_D(a,w),
 \qquad a,w\in D\,.
\end{equation*}
Since effective resistance is a metric~\cite[Theorem~2.64]{Barlow2017}, the triangle inequality and Proposition~\ref{prop:resistance-growth} yield
\begin{equation}\label{eq:Green-difference}
 |G_D(x,w)-G_D(y,w)|
 \le\widehat{\cR}_D(x,y)
 \le\cR(x,y)
 \le Cr^{d_w-d_f},
 \qquad w\in D\,.
\end{equation}
Moreover, by~\eqref{eq:Green-resistance} and Proposition~\ref{prop:resistance-growth},
\begin{equation*}
 G_D(y,y)=\widehat{\cR}_D(y,\partial)
 =\cR(y,D^c)\ge cR^{d_w-d_f}\,.
\end{equation*}
Combining these estimates with~\eqref{eq:Green-hitting}, we obtain
\begin{equation}\label{eq:miss-y}
 \mathbb P^x(\tau_D<T_y)
 =\frac{G_D(y,y)-G_D(x,y)}{G_D(y,y)}
 \le C(r/R)^{d_w-d_f}\,.
\end{equation}

Let $H_Du(z)=\mathbb E^z[u(X_{\tau_D})]$ for $z\in D$, and extend $H_Du$ with $u$ on $D^c$. Then $H_Du$ is harmonic in $D$ and $\|H_Du\|_\infty\le\|u\|_\infty$. The strong Markov property at $T_y\wedge\tau_D$ gives
\begin{equation*}
 H_Du(x)-H_Du(y)
 =\mathbb E^x\!\left[
 \bigl(u(X_{\tau_D})-H_Du(y)\bigr)
 \1_{\{\tau_D<T_y\}}\right]\,.
\end{equation*}
Hence, by~\eqref{eq:miss-y},
\begin{equation}\label{eq:harmonic-oscillation}
 |H_Du(x)-H_Du(y)|
 \le2\|u\|_\infty\mathbb P^x(\tau_D<T_y)
 \le C(r/R)^{d_w-d_f}\|u\|_\infty\,.
\end{equation}

Since $D$ is finite, the convergent Green series~\eqref{eq:killed-Green} gives $(I-P_D)^{-1}=\sum_{k\ge0}P_D^k$ on functions supported in $D$. As $u-H_Du$ vanishes on $D^c$ and satisfies $(I-P_D)(u-H_Du)=Lu$ in $D$, we have
\begin{equation*}
 u(z)=H_Du(z)+\sum_{w\in D}G_D(z,w)Lu(w)m(w),
 \qquad z\in D\,.
\end{equation*}
Using~\eqref{eq:Green-difference}, \eqref{eq:harmonic-oscillation}, and $m(D)\le CR^{d_f}$ from Proposition~\ref{prop:volume-growth}, we conclude that
\begin{align*}
 |u(x)-u(y)|
 &\le |H_Du(x)-H_Du(y)|
 +\sum_{w\in D}|G_D(x,w)-G_D(y,w)|\,|Lu(w)|m(w)\\
 &\le Cr^{d_w-d_f}\left(
 R^{d_f-d_w}\|u\|_\infty+R^{d_f}\|Lu\|_\infty\right)\,.
\end{align*}
\end{proof}

We now have all the ingredients to conclude the proof of Theorem~\ref{thm:wBE-strong-recurrence}.

\begin{proof}[Proof of Theorem~\ref{thm:wBE-strong-recurrence}]
The Markov property and Proposition~\ref{prop:temporal-regularity} give, respectively,
\[
 \|P^n f\|_\infty\le\|f\|_\infty,
 \qquad\text{and}\qquad
 \|LP^nf\|_\infty\le Cn^{-1}\|f\|_\infty.
\]
If $r=d(x,y)\le n^{1/d_w}$, Lemma~\ref{lem:local-Poisson} implies
\begin{align*}
 |P^nf(x)-P^nf(y)|
 &\le Cr^{d_w-d_f}\left(n^{(d_f-d_w)/d_w}
 +n^{d_f/d_w-1}\right)\|f\|_\infty\\
 &\le C\left(\frac r{n^{1/d_w}}\right)^{d_w-d_f}
 \|f\|_\infty.
\end{align*}

If $r>n^{1/d_w}$, the $L^{\infty}$ contractivity of $P^n$ yields $|P^nf(x)-P^nf(y)|\le 2\|f\|_\infty$, and the two estimates prove $({\rm wBE}(d_w-d_f))$.
The weak monotonicity assertion follows from Theorem~\ref{thm:HK-WM}.
\end{proof}

\section{Weak monotonicity of ball averaged functionals}\label{sec:BL}

In this section we study ball averaged functionals, the discrete counterparts of the Korevaar-Schoen functionals studied in \cite{Baudoin2024KS}. Under a uniform lower bound on the masses of vertices and~\eqref{eq:volume-growth}, we prove weak monotonicity with exponent~$1+(d_f-1)/p$ for every~$1\le p<\infty$. 

\subsection{Weak monotonicity from volume growth}\label{sec:BL-volume}

For every real $r\ge1$, define the ball-averaging kernel by
\begin{equation}\label{eq:BL-kernel}
 b_r(x,y):=\frac{\1_{B(x,r)}(y)}{\Vol(x,r)}\,.
\end{equation}
Thus $\sum_yb_r(x,y)m(y)=1$.

\begin{definition}[Ball averaged functionals]
\label{def:BL-functional}
For $1\le p<\infty$, $r\ge1$, and $f\in\ell^p(V,m)$, set, with values
allowed in $[0,\infty]$,
\begin{equation}\label{eq:BL-functional}
 \cB_{p,r}(f)
 :=\left(\sum_{x\in V}\frac{m(x)}{\Vol(x,r)}
 \sum_{y\in B(x,r)}|f(x)-f(y)|^pm(y)\right)^{1/p}.
\end{equation}

\end{definition}

Under~\eqref{eq:volume-growth}, one has
\begin{equation*}
 \sum_{x\in V}b_r(x,y)m(x)
 \le\frac{m(B(y,r))}{c_Vr^{d_f}}
 \le\frac{C_V}{c_V}\,.
\end{equation*}
Together with the row normalization, Minkowski's inequality therefore gives
\begin{equation}\label{eq:ell_p}
 \cB_{p,r}(f)\le\bigl(1+(C_V/c_V)^{1/p}\bigr)\|f\|_p
 \le C\|f\|_p\,,\qquad r\ge1.
\end{equation}

In the following proof, we combine the discretization method of~\cite{Kanai1985} with the geodesic path-counting argument of~\cite[proof of Theorem~1.1]{Chen2019}.

\begin{theorem}
\label{thm:BL-WM}
Assume $m_*:=\inf_{x\in V}m(x)>0$ and~\eqref{eq:volume-growth} for some $d_f>0$.
For $1\le p<\infty$, set
\begin{equation}\label{eq:sigma-p}
 \sigma_p:=1+\frac{d_f-1}{p}.
\end{equation}
Then there exists a constant $C>0$ depending on
$p,d_f,m_*^{-1},c_V,C_V$ such that
\begin{equation}\label{eq:BL-WM}
 r^{-\sigma_p}\cB_{p,r}(f)
 \le Cs^{-\sigma_p}\cB_{p,s}(f),
 \qquad r\ge s\ge1,
\end{equation}
for every $f\in\ell^p(V,m)$.  Thus the integer subfamily of the ball averaged functionals satisfies $({\rm WM}_{p,\sigma_p})$.  
\end{theorem}

\begin{proof}
We divide the proof into four steps.
\medskip

\emph{Step 1: vertex masses and pair sums.}
The assumptions give $m_*\le m(x)\le\Vol(x,1)\le C_V$. Hence $m(x)\simeq1$, and the graph has uniformly bounded degree. Condition~\eqref{eq:volume-growth} gives
\begin{equation}\label{eq:BL-pair-sum}
 \cB_{p,t}(f)^p
 \simeq t^{-d_f}\sum_{d(x,y)\le t}|f(x)-f(y)|^pm(x)m(y),
 \qquad t\ge1\,.
\end{equation}

\medskip

\emph{Step 2: a partition at scale $s$ and its associated graph.} 
Following the discretization approach of \cite{Kanai1985} (see also \cite{BarlowBassKumagai2006}), we construct a partition of $V$ at scale $s$ and an auxiliary graph roughly isometric to $G=(V,E)$.
It suffices first to consider $s\ge24$.  Put
$\varrho=\lfloor s/8\rfloor\simeq s$ and choose a maximal $\varrho$-separated family $(x_i)_{i\in I}$.  A nearest-center partition may be chosen so that
\begin{equation}\label{eq:coarse-cells}
 V=\bigsqcup_{i\in I}A_i,
 \qquad
 B(x_i,\varrho/3)\subset A_i\subset B(x_i,\varrho).
\end{equation}
Here the symbol $\bigsqcup$ means disjoint union.
Condition~\eqref{eq:volume-growth} gives
\begin{equation}\label{eq:cell-mass}
m(A_i)\simeq s^{d_f}.
\end{equation}

We define an auxiliary graph $G_c$ with vertex set $I$, which joins distinct indices $i,j\in I$ by an edge whenever $d(x_i,x_j)\le4\varrho$.  Then $G_c$ is connected.
Indeed, choose vertices $(z_\ell)_{\ell=0}^M\subset V$ along a geodesic in $G$ connecting $x_i$ and $x_j$ such that $z_0=x_i$, $z_M=x_j$ and 
\[
d(z_\ell,z_{\ell+1})\le \varrho, \quad M\le1+d(x_i,x_j)/\varrho.
\]
For each $\ell$, let $i_{\ell}\in I$ be the unique index such that $z_{\ell}\in A_{i_{\ell}}$. Then $d(x_{i_\ell},x_{i_{\ell+1}})\le3\varrho$ so either $x_{i_\ell}=x_{i_{\ell+1}}$ or they are  adjacent in \(G_{\mathrm c}\).  Hence the graph $G_c$ is connected. Denote its metric by $d_c$, the preceding construction gives that
\begin{equation}\label{eq:coarse-metric}
 d_{\rm c}(i,j)
 \le C\left(1+d(x_i,x_j)/\varrho\right),
 \qquad
 d(x_i,x_j)\le4\varrho\,d_{\rm c}(i,j).
\end{equation}

For every $i\in I$ and $R>0$, consider $B_c(i,R):=\{j\in I:d_c(i,j)\le R\}$.  For any $j\in B_c(i,R)$, one has $d(x_i,x_j)\le 4\varrho R$ and hence
\(
\bigsqcup_{j\in B_c(i,R)}B(x_j,\varrho/3) \subset B(x_i,4\varrho R+\varrho/3)
\), which yields
\begin{equation}\label{eq:coarse-volume}
 \#B_{\rm c}(i,R)\le C(1+R)^{d_f},
 \qquad i\in I,\quad R\ge1.
\end{equation}
In particular, taking \(R=1\) shows that $G_c$ has uniformly bounded degree.

\medskip 

\emph{Step 3: decomposition of $f$ and the estimate for $s\ge24$.}
Define
\begin{equation}\label{eq:cell-decomposition}
 f_{A_i}:=\frac1{m(A_i)}\sum_{x\in A_i}f(x)m(x),
 \qquad
 g:=\sum_i  f_{A_i}\1_{A_i},
 \qquad
 h:=f-g.
\end{equation}
Jensen's inequality, $\operatorname{diam}(A_i)\le2\varrho\le s$, and
\eqref{eq:BL-pair-sum} yield
\begin{align}
 \|h\|_p^p
 &\le\sum_i\frac1{m(A_i)}\sum_{x,y\in A_i}
 |f(x)-f(y)|^pm(x)m(y)\notag\\
 &\le Cs^{-d_f}\sum_{d(x,y)\le s}|f(x)-f(y)|^pm(x)m(y)
 \le C\cB_{p,s}(f)^p.                                \label{eq:h-bound}
\end{align}
Consequently, together with \eqref{eq:ell_p} we have
\begin{equation}\label{eq:BL-Lp-bound}
 \cB_{p,r}(h)^p\le C\|h\|_p^p
 \le C\cB_{p,s}(f)^p.
\end{equation}

Put $R=r/s\ge1$.  Equations \eqref{eq:BL-pair-sum},
\eqref{eq:cell-mass}, and the comparison of metrics \eqref{eq:coarse-metric} give
\begin{align}\label{eq:g-coarse-start}
 \cB_{p,r}(g)^p
&\le Cr^{-d_f}
\sum_{i,j}\!
\sum_{\substack{x\in A_i,\ y\in A_j\\d(x,y)\le r}}
|f_{A_i}-f_{A_j}|^p m(x)m(y)
\\ & \le Cs^{d_f}R^{-d_f}
 \sum_i\sum_{d_{\rm c}(i,j)\le CR}|f_{A_i}-f_{A_j}|^p.
\end{align}
We claim
\begin{equation}\label{eq:coarse-path-count}
 \sum_i\sum_{d_{\rm c}(i,j)\le CR}|f_{A_i}-f_{A_j}|^p
 \le CR^{p-1+2d_f}\sum_{i\sim_{\rm c}j}|f_{A_i}-f_{A_j}|^p.
\end{equation}
This can be shown using the geodesic path-counting idea in the proof of \cite[Theorem 1.1]{Chen2019}. Indeed, for each pair \(i,j\in I\) such that $d_{\mathrm c}(i,j)\le CR$, choose a shortest path with subsequent vertices $i_0=i,i_1,\ldots,i_L=j$, then $L\le CR$ and H\"older’s inequality give
\[
\begin{aligned}
|f_{A_i}-f_{A_j}|^p
&=\left|\sum_{k=0}^{L-1}
(f_{A_{i_k}}-f_{A_{i_{k+1}}})\right|^p\le CR^{p-1}\sum_{k=0}^{L-1}
|f_{A_{i_k}}-f_{A_{i_{k+1}}}|^p.
\end{aligned}
\]
Fix an edge with endpoint \(u\). Every chosen path containing that edge has both its endpoints in \(B_{\mathrm c}(u,CR)\). Since this ball contains at most \(CR^{d_f}\) vertices, there are at most
\(\bigl(\#B_{\mathrm c}(u,CR)\bigr)^2\le CR^{2d_f}\)
possible ordered endpoint pairs. Thus each edge is counted at most \(CR^{2d_f}\) times and the claim follows. 

For adjacent indices $i,j\in I$, Jensen's inequality gives
\[
 |f_{A_i}-f_{A_j}|^p
 \le\frac1{m(A_i)m(A_j)}\sum_{x\in A_i}\sum_{y\in A_j}
 |f(x)-f(y)|^pm(x)m(y).
\]
Here $d(x,y)\le6\varrho\le s$.  Recall from Step 2 that the graph $G_c$ has bounded degree. Then \eqref{eq:BL-pair-sum} and \eqref{eq:cell-mass} imply
\begin{equation}\label{eq:coarse-edge-bound}
 s^{d_f}\sum_{i\sim_{\rm c}j}|f_{A_i}-f_{A_j}|^p
 \le C\cB_{p,s}(f)^p.
\end{equation}
Combining \eqref{eq:g-coarse-start}, \eqref{eq:coarse-path-count} and \eqref{eq:coarse-edge-bound}, we
obtain
\begin{equation}\label{eq:g-final}
 \cB_{p,r}(g)^p
 \le C R^{p+d_f-1}\cB_{p,s}(f)^p.
\end{equation}

Since $f=g+h$, Minkowski's inequality, \eqref{eq:BL-Lp-bound},
and \eqref{eq:g-final} yield
\[
 \cB_{p,r}(f)^p
 \le C\bigl(\cB_{p,r}(g)^p+\cB_{p,r}(h)^p\bigr)
 \le C(r/s)^{p+d_f-1}\cB_{p,s}(f)^p,
\]
where we used $r/s\ge1$. Since $p\sigma_p=p+d_f-1$,
we obtain \eqref{eq:BL-WM} for $s\ge24$.
\medskip

\emph{Step 4: the estimate for $1\le s<24$.} Applying H\"older's inequality
along geodesics in $G$ and counting the paths
through each edge, as in Step~3, gives
\begin{equation}\label{eq:BL-bounded-scale-path}
 \sum_{d(x,y)\le r}|f(x)-f(y)|^pm(x)m(y)
 \le Cr^{p-1+2d_f}
 \sum_{x\sim y}|f(x)-f(y)|^p\,.
\end{equation}
The last sum is comparable to $\cB_{p,1}(f)^p$, since $m(x)\simeq1$ and $\Vol(x,1)\simeq1$. Multiplying by $r^{-d_f}$ and using \eqref{eq:BL-pair-sum} from Step 1, we obtain
\[
 \cB_{p,r}(f)^p
 \le Cr^{p+d_f-1}\cB_{p,1}(f)^p,
\]
and hence
\(\cB_{p,r}(f)\le Cr^{\sigma_p}\cB_{p,1}(f)\).
For $1\le s<24$, the sum in \eqref{eq:BL-pair-sum} at scale $s$ contains all adjacent pairs. Consequently,
\(\cB_{p,1}(f)\le C\cB_{p,s}(f)\) and 
it follows that
\[
 \cB_{p,r}(f) \le Cr^{\sigma_p}\cB_{p,s}(f)
 \le C 24^{\sigma_p}(r/s)^{\sigma_p}\cB_{p,s}(f).
\]
This proves \eqref{eq:BL-WM} for all $1\le s\le r$.
\end{proof}

\subsection{Comparison with heat kernel functionals}\label{sec:BL-heat}

Under~\eqref{eq:LB} and $({\rm SG}_{d_f,d_w})$, the heat kernel and ball averaged functionals are uniformly comparable. This transfers weak monotonicity to the heat kernel functionals for every finite $p\ge1$.

\begin{theorem}\label{thm:ball-heat-comparison}
Assume~\eqref{eq:LB} and $({\rm SG}_{d_f,d_w})$. Let $1\le p<\infty$ and $\sigma_p=1+(d_f-1)/p$. There is a constant $C\ge1$ such that
\begin{equation}\label{eq:ball-heat-comparison}
 C^{-1}\cB_{p,n^{1/d_w}}(f)
 \le\cN_{p,n}^{\rm H}(f)
 \le C\cB_{p,n^{1/d_w}}(f)
\end{equation}
for every $n\ge1$ and $f\in\ell^p(V,m)$. Moreover, for integers $r\ge s\ge1$
\begin{equation}\label{eq:heat-WM-from-balls}
 r^{-\sigma_p/d_w}\cN_{p,r}^{\rm H}(f)
 \le C s^{-\sigma_p/d_w}\cN_{p,s}^{\rm H}(f).
\end{equation}
 Thus the heat kernel functionals satisfy $({\rm WM}_{p,\sigma_p/d_w})$. 
\end{theorem}

\begin{proof}
Put $R=n^{1/d_w}$. Proposition~\ref{prop:volume-growth} gives~\eqref{eq:volume-growth} and the uniform lower bound on vertex masses required by Theorem~\ref{thm:BL-WM}. Split the sum defining $\cN_{p,n}^{\rm H}(f)^p$ into $d(x,y)\le R$ and the annuli $2^jR<d(x,y)\le2^{j+1}R$, $j\ge0$. By~\eqref{eq:SG-upper} and~\eqref{eq:BL-pair-sum},
\begin{align*}
 \cN_{p,n}^{\rm H}(f)^p
 &\le C\cB_{p,R}(f)^p
 +C\sum_{j\ge0}2^{(j+1)d_f}
 e^{-c2^{j d_w/(d_w-1)}}\cB_{p,2^{j+1}R}(f)^p\\
 &\le C\left(1+\sum_{j\ge0}2^{(j+1)(d_f+p\sigma_p)}
 e^{-c2^{j d_w/(d_w-1)}}\right)\cB_{p,R}(f)^p
 \le C\cB_{p,R}(f)^p\,,
\end{align*}
where the second inequality uses Theorem~\ref{thm:BL-WM}.

For the reverse inequality, $R\le2n-1$ for every $n\ge1$. The exact lower bound~\eqref{eq:exact-HK-lower} at time $2n$ therefore gives $p_{2n}(x,y)\ge c n^{-d_f/d_w}$ whenever $d(x,y)\le R$. By~\eqref{eq:BL-pair-sum} and time subadditivity,
\begin{equation*}
 \cB_{p,R}(f)\le C\cN_{p,2n}^{\rm H}(f)
 \le2C\cN_{p,n}^{\rm H}(f)\,.
\end{equation*}
This proves~\eqref{eq:ball-heat-comparison}. Finally, for integers $r\ge s\ge1$, Theorem~\ref{thm:BL-WM} yields
\begin{equation*}
 \cN_{p,r}^{\rm H}(f)
 \le C\cB_{p,r^{1/d_w}}(f)
 \le C(r/s)^{\sigma_p/d_w}\cB_{p,s^{1/d_w}}(f)
 \le C(r/s)^{\sigma_p/d_w}\cN_{p,s}^{\rm H}(f)\,,
\end{equation*}
which proves~\eqref{eq:heat-WM-from-balls}.
\end{proof}

Taking $n=1$ in~\eqref{eq:ball-heat-comparison} gives
\begin{equation}\label{eq:ball-gradient-comparison}
 \cB_{p,1}(f)\simeq\|\nabla_pf\|_p
 \qquad\text{under~\eqref{eq:LB} and $({\rm SG}_{d_f,d_w})$}\,.
\end{equation}
The same comparison holds under~\eqref{eq:p0} and~\eqref{eq:volume-growth}, by Remark~\ref{rem:controlled-weight-comparison} and~\eqref{eq:BL-pair-sum}.

\begin{remark}\label{R:comparison_ball_vs_path}
Assume~\eqref{eq:LB} and $({\rm SG}_{d_f,d_w})$ with $1\le d_f<d_w$. Theorem~\ref{thm:ball-heat-comparison} gives the heat kernel weak monotonicity with exponent $\sigma_p/d_w$. If, in addition, the weak Bakry--\'Emery estimate $({\rm wBE}(d_w-d_f))$ holds, then Theorem~\ref{thm:HK-WM} yields the exponent $\beta_p(d_w-d_f)$ for $1\le p\le 2$.  Note that
\begin{equation}\label{eq:prescribed-strict-improvement}
 \frac{\sigma_p}{d_w}-\beta_p(d_w-d_f)
 =\frac{\bigl(1-(d_w-d_f)\bigr)(1-1/p)}{d_w}.
\end{equation}
Consequently, when $d_w-d_f<1$ and $1<p\le 2$, the weak Bakry--\'Emery estimate gives a strictly smaller exponent. The two exponents coincide when $d_w-d_f=1$ or $p=1$. For $p\ge2$, Theorems~\ref{thm:HK-WM-p-ge-2} and~\ref{thm:ball-heat-comparison} give weak monotonicity with exponent $\min\{1/2,\sigma_p/d_w\}$.
\end{remark}

\vspace{0.6cm}

Assume either~\eqref{eq:p0} and~\eqref{eq:volume-growth}, or~\eqref{eq:LB} and $({\rm SG}_{d_f,d_w})$. Then Theorem~\ref{thm:BL-WM} applies with $\sigma_1=d_f$. The $\alpha$-perimeter of Definition~\ref{def:kernel-alpha-perimeter}, applied to the integer family $k_n=b_n$, therefore satisfies
\begin{equation}\label{eq:ball-kernel-perimeter-comparison}
 \mathcal P_{d_f,(b_n)}(A)
 =\sup_{n\ge1}n^{-d_f}\cB_{1,n}(\1_A)
 \simeq\cB_{1,1}(\1_A)
 \simeq\Per(A)
\end{equation}
for $m(A)<\infty$. The first comparison follows from $\mathrm{P}(1,d_f)$, which is implied by Theorem~\ref{thm:BL-WM}; the second follows from the scale-one gradient comparison above and $\|\nabla_1\1_A\|_1=2\Per(A)$.

\section{Critical exponents}\label{sec:critical-exponents}

We study the infima of the exponents for which the heat kernel and ball averaged functionals are bounded by the $\ell^p$-edge gradient functional. Those are the critical exponents for the property $\mathrm{P}(p,\alpha)$ and play the same role in our framework as the exponents respectively considered in \cite{AlonsoEtAl2020I} and \cite{Baudoin2024KS}.

\subsection{Definitions and comparison}

\begin{definition}\label{def:critical-exponents}
Let $1\le p<\infty$. The critical exponents for the heat kernel and ball averaged functionals are, respectively, defined by
\begin{equation}\label{eq:critical-heat-definition}
 \alpha_p^{\rm H}:=\inf\left\{\alpha\ge0:
 (\cN_{p,n}^{\rm H})_{n\ge1}\text{ satisfies }\mathrm{P}(p,\alpha)\right\}
\end{equation}
and
\begin{equation}\label{eq:critical-ball-definition}
 \alpha_p^{\rm B}:=\inf\left\{\alpha\ge0:
 \sup_{r\ge1}r^{-\alpha}\cB_{p,r}(f)
 \le C\|\nabla_pf\|_p\ \text{for every }f\in\ell^p(V,m)\right\}\,.
\end{equation}
The constants may depend on $p$, $\alpha$ and of course the graph, but not on $f$ or the scale; we set $\inf\varnothing=+\infty$. An exponent is \emph{admissible} if the corresponding property or estimate holds. The critical exponent is \emph{attained} if it is admissible.
\end{definition}

Admissibility is preserved when $\alpha$ increases. Thus every exponent strictly above a finite critical exponent is admissible; attainment requires a separate argument. 

Note that the infimum is the same if restrict the $\inf$ to finitely supported functions since if $f_j\to f$ in $\ell^p(V,m)$, then
\begin{equation*}
 \|\nabla_p(f_j-f)\|_p\le2\|f_j-f\|_p\longrightarrow0\,,
\end{equation*}
and Fatou's lemma applies to the pair sums.

These exponents differ from those defined by nontriviality or density of continuous Besov classes in~\cite[Section~5.2]{AlonsoEtAl2020I} and~\cite[Section~4.1]{Baudoin2024KS} but are the critical exponents for $\mathrm{P}(p,\alpha)$. Indeed, for every $\alpha\ge0$,
\begin{equation*}
 \{f\in\ell^p(V,m):\|f\|_{p,\alpha}<\infty\}=\ell^p(V,m)\,.
\end{equation*}
For heat kernel functionals, $\cN_{p,1}^{\rm H}(f)=\|\nabla_pf\|_p$ gives
\begin{equation*}
 ({\rm WM}_{p,\alpha})\Longrightarrow\mathrm{P}(p,\alpha)
 \quad\Longleftrightarrow\quad
 \|f\|_{p,\alpha}\le C\|\nabla_pf\|_p\,.
\end{equation*}
For ball averaged functionals, $\cB_{p,r}=\cB_{p,\lfloor r\rfloor}$ for $r\ge1$, so $\mathrm{P}(p,\alpha)$ for the integer family also controls every real radius. Under either~\eqref{eq:p0} and~\eqref{eq:volume-growth}, or~\eqref{eq:LB} and $({\rm SG}_{d_f,d_w})$, one has $\cB_{p,1}(f)\simeq\|\nabla_pf\|_p$. Under these hypotheses, ball admissibility is therefore equivalent to $\mathrm{P}(p,\alpha)$ for $(\cB_{p,n})_{n\ge1}$.

\begin{proposition}\label{prop:critical-comparison}
Assume~\eqref{eq:LB} and $({\rm SG}_{d_f,d_w})$. For every $1\le p<\infty$, an exponent $\alpha\ge0$ is admissible in~\eqref{eq:critical-heat-definition} if and only if $d_w\alpha$ is admissible in~\eqref{eq:critical-ball-definition}. Consequently,
\begin{equation}\label{eq:critical-comparison}
 \alpha_p^{\rm B}=d_w\alpha_p^{\rm H}\,,
\end{equation}
and one critical exponent is attained if and only if the other is attained.
\end{proposition}

\begin{proof}
For $r\ge1$, choose $n=\lceil r^{d_w}\rceil$. Then
\begin{equation*}
 r\le n^{1/d_w}\le2^{1/d_w}r\,,
 \qquad
 \cB_{p,r}(f)\le C\cB_{p,n^{1/d_w}}(f)
\end{equation*}
by~\eqref{eq:BL-pair-sum}. Together with Theorem~\ref{thm:ball-heat-comparison}, this gives, for every $\alpha\ge0$,
\begin{equation*}
 \sup_{r\ge1}r^{-d_w\alpha}\cB_{p,r}(f)
 \simeq\sup_{n\ge1}n^{-\alpha}\cN_{p,n}^{\rm H}(f)
 =\|f\|_{p,\alpha}\,,
\end{equation*}
which proves both assertions.
\end{proof}
\subsection{General bounds}

We first compare the exponents for different values of $p$, without heat kernel assumptions.

\begin{proposition}\label{prop:critical-p-comparison}
Let $1\le q<p<\infty$. For either the heat kernel or the ball averaged functionals, admissibility of $\alpha$ at $p$ implies admissibility of $p\alpha/q$ at $q$. Moreover,
\begin{equation}\label{eq:P-p-comparison}
 \mathrm{P}(p,\alpha)\Longrightarrow\mathrm{P}(q,p\alpha/q),
\end{equation}
and
\begin{equation}\label{eq:WM-p-comparison}
 ({\rm WM}_{p,\alpha})\Longrightarrow({\rm WM}_{q,p\alpha/q})\,.
\end{equation}
In particular,
\begin{equation}\label{eq:critical-p-comparison}
 q\alpha_q^{\rm H}\le p\alpha_p^{\rm H},
 \qquad q\alpha_q^{\rm B}\le p\alpha_p^{\rm B}\,.
\end{equation}
\end{proposition}

\begin{proof}
As in Proposition~\ref{prop:l2-identity}, we use the identity
\begin{equation*}
 \int_0^\infty |e^{ita}-e^{itb}|^p\,\frac{dt}{t^{1+q}}
 =c_{p,q}|a-b|^q\,,\qquad a,b\in\R\,,\quad 0<c_{p,q}<\infty\,.
\end{equation*}
The integral converges because $0<q<p$. For $f\in\ell^q(V,m)$,
\begin{equation*}
 |e^{itf}-1|^p\le2^{p-q}t^q|f|^q\,,
 \qquad e^{itf}-1\in\ell^p(V,m)\,.
\end{equation*}
Extend the assumed estimate to complex-valued functions by applying it to real and imaginary parts. In the notation of Definition~\ref{def:averaging-functional}, Tonelli's theorem gives
\begin{align*}
 \int_0^\infty\cN_{p,n}(e^{itf}-1)^p\,\frac{dt}{t^{1+q}}
 &=c_{p,q}\cN_{q,n}(f)^q\,,\\
 \int_0^\infty\|\nabla_p(e^{itf}-1)\|_p^p\,\frac{dt}{t^{1+q}}
 &=c_{p,q}\|\nabla_qf\|_q^q\,.
\end{align*}
Integrating the $p$-th power of the admissibility estimate therefore yields
\begin{equation*}
 \cN_{q,n}(f)^q\le Cn^{p\alpha}\|\nabla_qf\|_q^q\,.
\end{equation*}
The same calculation with $\cN_{p,1}$ or $\cN_{p,s}$ on the right proves~\eqref{eq:P-p-comparison} or~\eqref{eq:WM-p-comparison}, respectively. For ball averaged functionals, replace $k_n$ by $b_r$. Taking infima gives~\eqref{eq:critical-p-comparison}.
\end{proof}

The following cutoff criterion will be used to prove the lower bounds. See also the cutoff argument in~\cite[proof of Proposition~1.3]{Chen2019}.

\begin{proposition}\label{prop:critical-Poincare}
Let $1\le p<\infty$. Assume that for some $c_V,d_f>0$, one has $\Vol(x,r)\ge c_Vr^{d_f}$ for every $x\in V$ and $r\ge1$. Suppose that for some $\sigma\ge0$ there exist $R_j\to\infty$ and finitely supported functions $f_j$ satisfying
\begin{equation}\label{eq:critical-cutoffs}
 m(\operatorname{supp}f_j)\le CR_j^{d_f},\qquad
 \|f_j\|_p\ge cR_j^{d_f/p},\qquad
 \|\nabla_pf_j\|_p\le CR_j^{d_f/p-\sigma}\,.
\end{equation}
Then $\alpha_p^{\rm B}\ge\sigma$.
\end{proposition}

\begin{proof}
Choose $D\ge1$ so that $\Vol(x,DR_j)\ge2m(\operatorname{supp}f_j)$ for all $x$ and all sufficiently large $j$. For $x\in\operatorname{supp}f_j$, at least half the mass of $B(x,DR_j)$ lies outside $\operatorname{supp}f_j$. Thus every ball admissible exponent $\alpha$ satisfies
\begin{equation*}
 c R_j^{d_f}\le\tfrac12\|f_j\|_p^p
 \le\cB_{p,DR_j}(f_j)^p
 \le C R_j^{p\alpha}\|\nabla_pf_j\|_p^p
 \le C R_j^{d_f+p(\alpha-\sigma)}\,.
\end{equation*}
Letting $j\to\infty$ gives $\alpha\ge\sigma$.
\end{proof}
\begin{theorem}\label{thm:critical-bounds}
The following holds:
\begin{enumerate}[wide=0em,label={\rm(\roman*)}]
 \item Assume~\eqref{eq:p0} and~\eqref{eq:volume-growth}. For every $1\le p<\infty$,
\begin{equation}\label{eq:critical-ball-bounds}
 1\le\alpha_p^{\rm B}\le1+\frac{d_f-1}{p}\,.
\end{equation}
\item Assume~\eqref{eq:LB} and $({\rm SG}_{d_f,d_w})$, then for every $1\le p<\infty$,
\begin{equation}\label{eq:critical-heat-bounds}
 \max\left\{\frac1{d_w},\frac1{\max\{p,2\}}\right\}
 \le\alpha_p^{\rm H}
 \le\min\left\{\frac1{\min\{p,2\}},\frac{d_f+p-1}{pd_w}\right\}\,.
\end{equation}
Every exponent in the upper bounds of $\alpha_p^{\rm B}$ and $\alpha_p^{\rm H}$ is admissible. In particular,
\begin{equation}\label{eq:critical-two-and-limit}
 \alpha_2^{\rm H}=\frac12,\qquad
 \alpha_2^{\rm B}=\frac{d_w}{2},\qquad
 \lim_{p\to\infty}\alpha_p^{\rm H}=\frac1{d_w},\qquad
 \lim_{p\to\infty}\alpha_p^{\rm B}=1\,,
\end{equation}
and both critical exponents at $p=2$ are attained.
\item Assume~\eqref{eq:LB}, $({\rm SG}_{d_f,d_w})$ and $({\rm wBE}(\kappa))$, then
\begin{equation}\label{eq:critical-wBE-bound}
 \alpha_p^{\rm H}\le\beta_p(\kappa)
 =\frac1p+\left(1-\frac2p\right)\frac\kappa{d_w}\,,
 \qquad 1\le p\le2\,,
\end{equation}
with $\beta_p(\kappa)$ admissible, whereas $\alpha_p^{\rm H}\ge\beta_p(\kappa)$ for $2\le p<\infty$.
\end{enumerate}
\end{theorem}

\begin{proof}
\begin{enumerate}[wide=0em,label={\rm (\roman*)}]
\item The upper bound for $\alpha_p^{\rm B}$ and its admissibility follow from Theorems~\ref{thm:BL-WM}. For the lower bound, we fix $o\in V$ and put $f_R(x)=(1-d(o,x)/R)_+$, $R\ge2$. By~\eqref{eq:volume-growth} and~\eqref{eq:gradient-edge-identity},
\begin{equation*}
 \|\nabla_pf_R\|_p^p\le CR^{-p}m(B(o,R+1))\le CR^{d_f-p}\,.
\end{equation*}
Choose $A>2$ such that $c_VA^{d_f}\ge2C_V$. Then
\begin{equation*}
 m(B(x,AR)\setminus B(o,R))\ge\frac{c_V}{2C_V}\Vol(x,AR)\,,
 \qquad x\in V\,.
\end{equation*}
Since $f_R\ge1/2$ on $B(o,R/2)$,
\begin{align*}
 \cB_{p,AR}(f_R)^p&\ge 
 \sum_{x\in B(o,R/2)}\frac{m(x)}{V(x,AR)} \sum_{y\in B(x,AR)\setminus B(o,R)}|f_R(x)-f_R(y)|^pm(y)
 \\&\ge2^{-p-1}\frac{c_V}{C_V}m(B(o,R/2))
 \ge cR^{d_f}\,.
\end{align*}
Every admissible exponent $\alpha$ satisfies $R^{d_f}\le CR^{p\alpha+d_f-p}$ for $R\ge2$.
Hence, $\alpha_p^{\rm B}\ge1$. 

\item The upper bound for $\alpha_p^{\rm H}$ and its admissibility are consequence of Theorems~\ref{thm:HK-WM-p-ge-2} and \ref{thm:ball-heat-comparison}. It remains to show the lower bound. Put $n=\lceil R^{d_w}\rceil$. By summing~\eqref{eq:SG-upper} over annuli, we can choose a fixed $K>2$ such that
\begin{equation*}
g_R:=P^n\1_{B(o,KR)}\ge\frac34\quad\text{on }B(o,R),
 \qquad
 g_R\le\frac14\quad\text{on }V\setminus B(o,2KR)\,.
\end{equation*}
The function $h_R=(g_R-1/4)_+$ satisfies
\begin{equation*}
 0\le h_R\le1\,,\qquad
 h_R\ge\frac12\ \text{on }B(o,R),\qquad
 \operatorname{supp}h_R\subset B(o,2KR)\,.
\end{equation*}
Since the mapping $t\mapsto(t-1/4)_+$ is $1$-Lipschitz, we have
$\|\nabla_2 h_R\|_2^2\le \|\nabla_2 g_R\|_2^2$. 
On the other hand, the condition \eqref{eq:LB} implies $\operatorname{spec}(P)\subset[-1+2\varepsilon_{\rm L},1]$ and note that
\begin{equation*}
 \sup_{-1+2\varepsilon_{\rm L}\le\lambda\le1}
 (1-\lambda)|\lambda|^{2n}\le\frac Cn\,.
\end{equation*}
Then the spectral calculus yields
\begin{equation*}
 \|\nabla_2 h_R\|_2^2\le \|\nabla_2g_R\|_2^2=2\|(I-P)^{1/2}P^n\1_{B(o,KR)}\|_2^2
 \le\frac Cn\|\1_{B(o,KR)}\|_2^2
 \le CR^{d_f-d_w}\,.
\end{equation*}
Next applying H\"older's inequality for $p\le2$, and  $|h_R(x)-h_R(y)|\le1$ for $p\ge2$, we have 
\begin{equation*}
 \|\nabla_ph_R\|_p^p\le
 \begin{cases}
 CR^{d_f-pd_w/2},&1\le p\le2,\\
 CR^{d_f-d_w},&2\le p<\infty.
 \end{cases}
\end{equation*}
Choose a fixed $D>2K$ sufficiently large and set $T=\lceil(DR)^{d_w}\rceil$. The heat kernel upper bound~\eqref{eq:SG-upper} gives
\begin{equation*}
 P^T\1_{\operatorname{supp}h_R}
 \le CT^{-d_f/d_w}m(B(o,2KR))\le\frac12\,.
\end{equation*}
By the conservativeness of $P^T$, it follows that
\begin{align*}
 \cN_{p,T}^{\rm H}(h_R)^p
 &\ge\sum_{x\in B(o,R)}|h_R(x)|^p P^T\1_{(\operatorname{supp}h_R)^c}(x)m(x)  
 \\ &\ge\frac12\sum_{x\in B(o,R)}|h_R(x)|^pm(x)
 \ge cR^{d_f}\,.
\end{align*}
Hence any admissible heat kernel exponent $\alpha$ satisfies 
\begin{equation*}
 R^{d_f}\le CR^{p\alpha d_w}
 \begin{cases}
 R^{d_f-pd_w/2},&1\le p\le2,\\
 R^{d_f-d_w},&2\le p<\infty,
 \end{cases}
 \quad \text{ for }R\ge 2
 \,.
\end{equation*}
As a consequence, $\alpha_p^{\rm H}\ge\frac1{\max\{p,2\}}$.

The identities and limits in~\eqref{eq:critical-two-and-limit} now follow from~\eqref{eq:critical-heat-bounds} and~\eqref{eq:critical-comparison}.

\item For $1\le p\le 2$, Theorem~\ref{thm:HK-WM} gives~\eqref{eq:critical-wBE-bound} and the admissibility of $\beta_p(\kappa)$. For $p\ge 2$, we follow the strategy in the proof of (ii). If $({\rm wBE}(\kappa))$ holds, then
\begin{equation*}
 \sup_{x\sim y}|h_R(x)-h_R(y)|
 \le\sup_{x\sim y}|g_R(x)-g_R(y)|\le CR^{-\kappa}\,,
\end{equation*}
and hence
\begin{equation*}
 \|\nabla_ph_R\|_p^p
 \le CR^{-(p-2)\kappa}\|\nabla_2h_R\|_2^2
 \le CR^{d_f-d_w-(p-2)\kappa}\,.
\end{equation*}
Recall that $\cN_{p,T}^{\rm H}(h_R)^p\ge cR^{d_f}$. With the sharper bound above, any admissible $\alpha$ satisfies 
\begin{equation*}
 R^{d_f}\le CR^{p\alpha d_w+d_f-d_w-(p-2)\kappa}.
\end{equation*}
We conclude that $\alpha_{p}^{\rm H}\ge\frac{d_w+(p-2)\kappa}{pd_w}=\beta_p(\kappa)$.
\end{enumerate}
\end{proof}

Under~\eqref{eq:LB} and $({\rm SG}_{d_f,d_w})$ with $1\le d_f<d_w$, set $\kappa=d_w-d_f$. Theorems~\ref{thm:wBE-strong-recurrence} and~\ref{thm:critical-bounds} give
\begin{equation}\label{eq:critical-strongly-recurrent-small-p}
 \frac12\le\alpha_p^{\rm H}\le\beta_p(\kappa),
 \qquad 1\le p\le2\,,
\end{equation}
where $\beta_p(\kappa)$ is admissible. The bounds for $p\ge2$ are
\begin{equation}\label{eq:critical-strongly-recurrent-large-p}
 \max\left\{\frac1{d_w},\beta_p(\kappa)\right\}
 \le\alpha_p^{\rm H}
 \le\min\left\{\frac12,\frac{d_f+p-1}{pd_w}\right\},
 \qquad 2\le p<\infty.
\end{equation}
The corresponding ball bounds are obtained by multiplying by $d_w$. At $p=1$, equality in the upper bound follows for graphs with arbitrarily large sets of uniformly bounded perimeter, by Proposition~\ref{prop:critical-perimeter} below.

\begin{corollary}\label{cor:critical-dw-df-one}
Assume~\eqref{eq:LB} and $({\rm SG}_{d_f,d_w})$ with $d_w=d_f+1$. For every $2\le p<\infty$,
\begin{equation}\label{eq:critical-dw-df-one}
 \alpha_p^{\rm H}=\frac{d_f+p-1}{pd_w},
 \qquad \alpha_p^{\rm B}=\frac{d_f+p-1}{p}\,.
\end{equation}
Both critical exponents are attained, and both families satisfy weak monotonicity at their critical exponents.
\end{corollary}

\begin{proof}
Since $d_w\ge2$, one has $d_f=d_w-1\ge1$, so Theorem~\ref{thm:wBE-strong-recurrence} gives $({\rm wBE}(1))$. For $p\ge2$, Theorem~\ref{thm:critical-bounds} yields
\begin{equation*}
 \frac{d_f+p-1}{pd_w}=\beta_p(1)
 \le\alpha_p^{\rm H}\le\frac{d_f+p-1}{pd_w}\,.
\end{equation*}
Proposition~\ref{prop:critical-comparison} gives $\alpha_p^{\rm B}=d_w\alpha_p^{\rm H}$. Weak monotonicity and attainment follow from Theorems~\ref{thm:ball-heat-comparison} and~\ref{thm:BL-WM}.
\end{proof}

\subsection{Poincar\'e inequalities and sharpness}

\begin{remark}
Fix $1\le p<\infty$ and assume~\eqref{eq:p0} and~\eqref{eq:volume-growth}. A local Poincar\'e inequality of the form
\begin{equation}\label{eq:critical-local-Poincare}
 \sum_{z\in B}|f(z)-f_B|^pm(z)
 \le Cr^{p\sigma}\sum_{z\in B(x,Ar)}|\nabla_pf(z)|^pm(z)
\end{equation}
holding for every $B=B(x,r)$, $r\ge1$, and every $f\in\ell^p(V,m)$ with fixed $A\ge1$ and $\sigma\ge0$, makes $\sigma$ admissible for the ball averaged functionals. Here $f_B$ is the average of $f$ on $B$. Indeed, cover $V$ by balls $B(x_i,r)$ about a maximal $r$-separated set. Every pair at distance at most $r$ lies in some $B(x_i,3r)$. Using~\eqref{eq:BL-pair-sum} and the bound $m(B(x_i,3r))\le Cr^{d_f}$, we obtain
\begin{align*}
 \cB_{p,r}(f)^p
 &\le C\sum_i\sum_{z\in B(x_i,3r)}|f(z)-f_{B(x_i,3r)}|^pm(z)\\
 &\le Cr^{p\sigma}\|\nabla_pf\|_p^p\,.
\end{align*}
The last inequality uses~\eqref{eq:critical-local-Poincare} and bounded overlap of $B(x_i,3Ar)$ by~\eqref{eq:volume-growth}. If there also exist cutoffs satisfying~\eqref{eq:critical-cutoffs}, Proposition~\ref{prop:critical-Poincare} gives $\alpha_p^{\rm B}=\sigma$ with attainment. See~\cite[Section~5.2]{Baudoin2024KS} and~\cite{Chen2019}.
\end{remark}

At $p=1$, coarea reduces admissibility to estimates for sets.

\begin{proposition}\label{prop:critical-perimeter}
An exponent $\alpha\ge0$ is admissible for the heat kernel functionals at $p=1$ if and only if, uniformly over $n\ge1$ and finite sets $A\subset V$,
\begin{equation}\label{eq:critical-perimeter-test}
 \sum_{x\in A}P^n\1_{A^c}(x)m(x)
 \le Cn^\alpha\Per(A)\,.
\end{equation}
For the ball averaged functionals, admissibility is equivalent to
\begin{equation*}
 \cB_{1,r}(\1_A)\le Cr^\alpha\Per(A)\,,
 \qquad r\ge1\,,\quad A\subset V\ \text{finite}\,.
\end{equation*}
Under~\eqref{eq:LB} and $({\rm SG}_{d_f,d_w})$, if finite sets $A_j$ satisfy $m(A_j)\to\infty$ and $\sup_j\Per(A_j)<\infty$, then
\begin{equation}\label{eq:critical-bounded-perimeter}
 \alpha_1^{\rm H}=\frac{d_f}{d_w},\qquad
 \alpha_1^{\rm B}=d_f\,,
\end{equation}
and both critical exponents are attained.
\end{proposition}

\begin{proof}
Symmetry and~\eqref{eq:graph-perimeter} give
\begin{equation*}
 \cN_{1,n}^{\rm H}(\1_A)=2\sum_{x\in A}P^n\1_{A^c}(x)m(x)\,,
 \qquad \|\nabla_1\1_A\|_1=2\Per(A)\,.
\end{equation*}
For finitely supported $f$, use the layer-cake identity underlying~\eqref{eq:kernel-coarea},
\begin{equation*}
 |f(x)-f(y)|=\int_{\R}|\1_{\{f>t\}}(x)-\1_{\{f>t\}}(y)|\,dt\,.
\end{equation*}
Apply the set estimate to $\{f>t\}$, replacing cofinite sets by their finite complements, and integrate, as in the proof of Corollary~\ref{cor:kernel-isoperimetry}. Approximation extends the resulting functional estimate to $\ell^1(V,m)$. The same argument applies to ball averaged functionals without symmetry.

For the last assertion, take $R_j=m(A_j)^{1/d_f}$ and $f_j=\1_{A_j}$. Since $\|f_j\|_1=R_j^{d_f}$ and $\|\nabla_1f_j\|_1=2\Per(A_j)$ is bounded, Proposition~\ref{prop:critical-Poincare} applies with $\sigma=d_f$. The upper bound $d_f$ for the ball exponent is admissible by Theorem~\ref{thm:ball-heat-comparison} and Proposition~\ref{prop:critical-comparison}, which also gives the value and attainment of the heat exponent.
\end{proof}

\subsection{Stability under tensor products}

For a weighted graph $G=(V,E)$ and an integer $N\ge1$, define the product $G^{\otimes N}$ on $V^N$ by
\begin{equation}\label{eq:tensor-product-kernel}
 m_N(\mathbf x):=\prod_{i=1}^Nm(x_i),
 \qquad
 \pi_N(\mathbf x,\mathbf y):=\prod_{i=1}^N\pi(x_i,y_i),
 \qquad
 c_N(\mathbf x,\mathbf y):=\prod_{i=1}^Nc_{x_i y_i}\,.
\end{equation}
Note that for the walk several coordinates may therefore change in one step. Norms on $V^N$ are taken with respect to $m_N$; the gradient $\nabla_p^{(N)}$ and perimeter $\operatorname{Per}_{G^{\otimes N}}$ are defined by~\eqref{eq:p-gradient} and~\eqref{eq:graph-perimeter}, with $\pi_N$ and $c_N$ in place of $\pi$ and $c$.

Under~\eqref{eq:LB}, each coordinate path can wait at its endpoint, so
\begin{equation}\label{eq:product-distance-balls}
 d_N(\mathbf x,\mathbf y)=\max_{1\le i\le N}d(x_i,y_i),
 \qquad
 B_N(\mathbf x,r)=\prod_{i=1}^NB(x_i,r)\,.
\end{equation}

\begin{proposition}\label{prop:critical-tensor}
Let $1\le p<\infty$. Assume controlled weights and~\eqref{eq:volume-growth}. For every integer $N\ge1$, the graphs $G$ and $G^{\otimes N}$ have the same admissible exponents for the heat kernel functionals, and also for the ball averaged functionals. In particular, both critical exponents and their attainment are preserved by finite tensor powers.
\end{proposition}

\begin{proof}

We may assume $N\ge2$. For $F\in\ell^p(V^N,m_N)$, $1\le i\le N$ and $\mathbf w\in V^{N-1}$, let
\begin{equation*}
 F_{i,\mathbf w}(u):=F(w_1,\ldots,w_{i-1},u,w_i,\ldots,w_{N-1})\,,\qquad u\in V.
\end{equation*}
Since $\sum_{\mathbf w}m_{N-1}(\mathbf w)\|F_{i,\mathbf w}\|_p^p=\|F\|_p^p$, it belongs to $\ell^p(V,m)$. For $\mathbf x,\mathbf y\in V^N$, put $\mathbf z^{(i)}:=(y_1,\ldots,y_i,x_{i+1},\ldots,x_N)$, $0\le i\le N$, so that
\begin{equation*}
 F(\mathbf x)-F(\mathbf y)=\sum_{i=1}^N\bigl(F(\mathbf z^{(i-1)})-F(\mathbf z^{(i)})\bigr)\,.
\end{equation*}
The $i$th difference depends neither on $x_j$, $j<i$, nor on $y_j$, $j>i$, and the heat kernel of $G^{\otimes N}$ is $\prod_{j=1}^Np_n(x_j,y_j)$. Minkowski's inequality and summation over these variables, using $\sum_{x\in V}p_n(x,y)m(x)=\sum_{y\in V}p_n(x,y)m(y)=1$, therefore give
\begin{equation}\label{eq:tensor-telescoping}
 \cN_{p,n}^{{\rm H},G^{\otimes N}}(F)
 \le\sum_{i=1}^N\Bigl(\sum_{\mathbf w\in V^{N-1}}m_{N-1}(\mathbf w)\,\cN_{p,n}^{\rm H}(F_{i,\mathbf w})^p\Bigr)^{1/p}\,.
\end{equation}
On the other hand, keeping in $\|\nabla_p^{(N)}F\|_p^p=\sum_{\mathbf x,\mathbf y\in V^N}c_N(\mathbf x,\mathbf y)|F(\mathbf x)-F(\mathbf y)|^p$ only the pairs with $x_j=y_j$ for all $j\ne i$, and using $c_{x_jx_j}\ge\varepsilon_{\rm L}m(x_j)$ from~\eqref{eq:LB}, we obtain
\begin{equation}\label{eq:tensor-gradient-sections}
 \sum_{\mathbf w\in V^{N-1}}m_{N-1}(\mathbf w)\|\nabla_pF_{i,\mathbf w}\|_p^p
 \le\varepsilon_{\rm L}^{-(N-1)}\|\nabla_p^{(N)}F\|_p^p\,,
 \qquad 1\le i\le N\,.
\end{equation}
If $\cN_{p,n}^{\rm H}(f)\le Cn^\alpha\|\nabla_pf\|_p$ for all $f\in\ell^p(V,m)$ and $n\ge1$, then applying this estimate to each section in~\eqref{eq:tensor-telescoping} and using~\eqref{eq:tensor-gradient-sections} yields
\begin{equation*}
 \cN_{p,n}^{{\rm H},G^{\otimes N}}(F)
 \le NC\varepsilon_{\rm L}^{-(N-1)/p}n^\alpha\|\nabla_p^{(N)}F\|_p\,,
\end{equation*}
so $\alpha$ is admissible for the heat kernel functionals on $G^{\otimes N}$.

For the ball averaged functionals, the method is similar. The ball-averaging kernel of $G^{\otimes N}$ is $\prod_{j=1}^Nb_r(x_j,y_j)$ by~\eqref{eq:product-distance-balls}, with $b_r$ as in~\eqref{eq:BL-kernel}, and the same argument applies with one change: $b_r$ is not symmetric, so each sum over $x_j$, $j<i$, is only bounded by $\sum_{x\in V}b_r(x,y)m(x)\le C_V/c_V$; see the estimate preceding~\eqref{eq:ell_p}. Thus~\eqref{eq:tensor-telescoping} holds for $\cB_{p,r}^{G^{\otimes N}}$ and $\cB_{p,r}$ with the additional factor $(C_V/c_V)^{(N-1)/p}$ on the right-hand side, and admissibility on $G^{\otimes N}$ follows as before.

Conversely, fix $o\in V$, put $h_R(x)=(1-d(o,x)/R)_+$, and define
\begin{equation*}
 F_R(x_1,\ldots,x_N)=f(x_1)\prod_{j=2}^Nh_R(x_j)\,.
\end{equation*}
Condition~\eqref{eq:volume-growth}, time subadditivity and~\eqref{eq:BL-pair-sum} give
\begin{equation*}
 \frac{\|\nabla_ph_R\|_p}{\|h_R\|_p}\le\frac CR\,,\qquad
 \frac{\cN_{p,n}^{\rm H}(h_R)}{\|h_R\|_p}\le\frac{Cn}{R}\,,\qquad
 \frac{\cB_{p,r}(h_R)}{\|h_R\|_p}\le\frac{Cr}{R}\,.
\end{equation*}
Coordinate telescoping as before and the reverse triangle inequality yield
\begin{align*}
 \frac{\cN_{p,n}^{{\rm H},G^{\otimes N}}(F_R)}{\|h_R\|_p^{N-1}}
 &\ge\cN_{p,n}^{\rm H}(f)-C_N\|f\|_p\frac nR,\\
 \frac{\|\nabla_p^{(N)}F_R\|_p}{\|h_R\|_p^{N-1}}
 &\le\|\nabla_pf\|_p+C_N\|f\|_p\frac1R\,.
\end{align*}
If $\alpha$ is admissible on the product, then
\begin{equation*}
 \cN_{p,n}^{\rm H}(f)-C_N\|f\|_p\frac nR
 \le Cn^\alpha\left(\|\nabla_pf\|_p+C_N\|f\|_p\frac1R\right)\,.
\end{equation*}
Letting $R\to\infty$ for each fixed $n$ gives the factor estimate, with $C$ independent of $n$. For ball averaged functionals, the same argument uses the third cutoff estimate and the uniform column bounds underlying~\eqref{eq:ell_p}; the constants remain independent of the fixed radius $r$.
\end{proof}

\subsection{Stability under rough isometries}

We recall that a rough isometry in the sense of~\cite{Kanai1985} is a map $\Phi:V\to V'$ for which there exist $A\ge1$ and $B\ge0$ such that
\begin{equation*}
 A^{-1}d(x,y)-B\le d'(\Phi x,\Phi y)\le Ad(x,y)+B\,,
 \qquad x,y\in V\,,
\end{equation*}
and every vertex of $V'$ lies within a fixed distance of $\Phi(V)$.

\begin{proposition}\label{prop:critical-rough-isometry}
Let $G,G'$ be roughly isometric graphs satisfying~\eqref{eq:volume-growth} with the same $d_f$, and let $1\le p<\infty$.
\begin{itemize}
\item If both graphs satisfy~\eqref{eq:p0}, their admissible exponents for the ball averaged functionals coincide.
\item If both satisfy~\eqref{eq:LB} and $({\rm SG}_{d_f,d_w})$ with the same $d_w$, their admissible exponents coincide for both the ball averaged and heat kernel functionals, without assuming~\eqref{eq:p0}. 
\end{itemize}
In each case, the corresponding critical exponents and their attainment are preserved.
\end{proposition}

\begin{proof}
Under~\eqref{eq:p0}, Remark~\ref{rem:controlled-weight-comparison} gives $m(x)\simeq1$ and $c_{xy}\simeq1$ on off-diagonal edges. Under~\eqref{eq:LB} and $({\rm SG}_{d_f,d_w})$, use Proposition~\ref{prop:volume-growth} and~\eqref{eq:ball-gradient-comparison} instead. In either case, both graphs have bounded degrees, $m(x)\simeq1$, and
\begin{equation*}
 \|\nabla_pf\|_p^p\simeq\sum_{x\sim y}|f(x)-f(y)|^p\,.
\end{equation*}
Let $\Psi:V'\to V$ be a coarse inverse of $\Phi$. Both maps have uniformly bounded fibres, and $\sup_xd(x,\Psi\Phi x)<\infty$. For $f\in\ell^p(V,m)$, put $g=f\circ\Psi$. Summing differences along paths of bounded length gives
\begin{equation*}
 \|\nabla_pg\|_{p,G'}\le C\|\nabla_pf\|_{p,G},\qquad
 \|f-g\circ\Phi\|_{p,G}\le C\|\nabla_pf\|_{p,G}\,.
\end{equation*}
The bounded fibres and~\eqref{eq:BL-pair-sum} also give, for fixed $A'\ge1$, $B'\ge0$,
\begin{equation*}
 \cB_{p,r}^G(g\circ\Phi)\le C\cB_{p,A'r+B'}^{G'}(g)\,.
\end{equation*}
If $\alpha\ge0$ is admissible on $G'$, then~\eqref{eq:ell_p} yields
\begin{align*}
 \cB_{p,r}^G(f)
 &\le C\|f-g\circ\Phi\|_{p,G}
      +C\cB_{p,A'r+B'}^{G'}(g)\\
 &\le C\bigl(1+(A'r+B')^\alpha\bigr)\|\nabla_pf\|_{p,G}
 \le Cr^\alpha\|\nabla_pf\|_{p,G}\,.
\end{align*}
Interchanging $G$ and $G'$ proves equality of the admissible exponents for the ball averaged functionals. Proposition~\ref{prop:critical-comparison} gives the heat kernel assertion.
\end{proof}

\section{Applications and examples}\label{sec:examples}

We compute in this section the critical exponents and recover classical functional inequalities on several spaces that fit into the discussed framework: Cayley graphs of polynomial growth, strongly recurrent graphs, three fractal examples and products. Condition~\eqref{eq:LB} is assumed throughout the section. 

\subsection{Cayley graphs of polynomial growth}
\label{sec:examples-cayley}\label{sec:examples-lattices}

Let $\Gamma$ be an infinite group with identity $e$ generated by a finite symmetric set $S\subset\Gamma\setminus\{e\}$. Write $|x|_S$ for the length of a word $x$ with respect to $S$, and equip its Cayley graph with the counting measure and the distance function $d(x,y)=|x^{-1}y|_S$. Further, consider the Markov operator
\begin{equation*}
 Pf(x):=\frac12f(x)+\frac1{2|S|}\sum_{a\in S}f(xa),
\end{equation*}
which corresponds to $c_{xx}=1/2$, $c_{x,xa}=1/(2|S|)$ for $a\in S$, and $m(x)=1$.

\begin{theorem}\label{thm:cayley-WM}
Assume $(\mathrm V_D)$ for some $D\ge1$. For every $1\le p<\infty$, there exists $C_p>0$ such that for any $f\in\ell^p(\Gamma)$,
\begin{equation*}
 r^{-1}\cB_{p,r}(f)
 \le C_p s^{-1}\cB_{p,s}(f),
 \qquad r\ge s\ge1\,,
\end{equation*}
and
\begin{equation*}
 n^{-1/2}\cN_{p,n}^{\rm H}(f)
 \le C_p k^{-1/2}\cN_{p,k}^{\rm H}(f)
\end{equation*}
for all integers $n\ge k\ge1$. 
Moreover,
\begin{equation}\label{eq:critical-lattices}
 \alpha_p^{\rm H}=\frac12,\qquad
 \alpha_p^{\rm B}=1,\qquad 1\le p<\infty\, ,
\end{equation}
and both critical exponents are attained.
\end{theorem}

\begin{proof}
Since $\pi(x,x)=1/2$ and $\pi(x,xa)=1/(2|S|)$ for $a\in S$, the weights are controlled. As $m\equiv1$, we have $p_n(x,y)=\mu^{(n)}(x^{-1}y)$, where $\mu^{(n)}$ is the $n$-fold convolution power of the symmetric probability density $\mu:=\bigl(\1_{\{e\}}+|S|^{-1}\1_S\bigr)/2$, whose support $S\cup\{e\}$ is finite and generates $\Gamma$. The Gaussian estimates~(13) and~(15) in~\cite[Theorem~5.1]{HebischSaloffCoste1993} therefore give~\eqref{eq:SG-upper} with $(d_f,d_w)=(D,2)$ for all $x,y\in\Gamma$, and~\eqref{eq:SG-lower} for $d(x,y)<n/C''$, where $C''$ is the constant in~(15). If $n/C''\le d(x,y)\le n$, following a geodesic from $x$ to $y$ and then holding at $y$ gives
\begin{equation*}
 p_n(x,y)\ge(2|S|)^{-d(x,y)}\,2^{-(n-d(x,y))}\ge(2|S|)^{-n}
 \ge\exp\Bigl(-(C'')^2\log(2|S|)\,\frac{d(x,y)^2}{n}\Bigr)\,.
\end{equation*}
Since $n^{-D/2}\le1$, this proves~\eqref{eq:SG-lower} in the remaining range, with $c_4\ge(C'')^2\log(2|S|)$. Thus $({\rm SG}_{D,2})$ holds. Next, let $y=xa$ with $a\in S$, and put $g=z^{-1}x$, so that $|g|_S=d(x,z)$. By symmetry of $\mu^{(n)}$, $p_n(x,z)-p_n(y,z)=\mu^{(n)}(g)-\mu^{(n)}(ga)$, and the gradient estimate in~\cite[Theorem~5.1, equation~(14)]{HebischSaloffCoste1993} gives
\begin{equation*}
 |p_n(x,z)-p_n(y,z)|
 \le Cn^{-(D+1)/2}
       \exp\left(-c\frac{d(x,z)^2}{n}\right),
 \qquad x\sim y,\quad n\ge1\,.
\end{equation*}
Splitting $\Gamma$ into $B(x,\sqrt n)$ and dyadic annuli, and using $(\mathrm V_D)$, yields
\begin{equation*}
 \sum_{z\in\Gamma}\exp\left(-c\frac{d(x,z)^2}{n}\right)
 \le Cn^{D/2}\left(1+\sum_{j\ge0}2^{(j+1)D}e^{-c4^j}\right)
 \le Cn^{D/2}\,.
\end{equation*}
Hence, for every bounded $g$ and adjacent $x,y$,
\begin{equation*}
 |P^ng(x)-P^ng(y)|
 \le\|g\|_\infty\sum_{z\in\Gamma}|p_n(x,z)-p_n(y,z)|
 \le Cn^{-1/2}\|g\|_\infty\,.
\end{equation*}
Summing along a shortest path and using the $\ell^\infty$ contraction of $P^n$, we obtain
\begin{equation*}
 |P^ng(x)-P^ng(y)|
 \le C\min\left\{1,\frac{d(x,y)}{\sqrt n}\right\}\|g\|_\infty,
 \qquad x,y\in\Gamma,\quad n\ge1\,.
\end{equation*}
This is $({\rm wBE}(1))$ with $d_w=2$, for which
\begin{equation*}
 \beta_p(1)=\frac1p+\frac12\left(1-\frac2p\right)=\frac12\,.
\end{equation*}
Theorem~\ref{thm:HK-WM} for $1\le p\le2$ and Theorem~\ref{thm:HK-WM-p-ge-2} for $p\ge2$ therefore give
\begin{equation*}
 \cN_{p,n}^{\rm H}(f)
 \le C_p\sqrt{\frac nk}\,\cN_{p,k}^{\rm H}(f),
 \qquad n\ge k\ge1\,.
\end{equation*}
For real $r\ge s\ge1$, put $n=\lceil r^2\rceil$ and $k=\lfloor s^2\rfloor$. By volume growth and Theorem~\ref{thm:ball-heat-comparison},
\begin{equation*}
 \cB_{p,r}(f)
 \le C_p\cN_{p,n}^{\rm H}(f)
 \le C_p\sqrt{\frac nk}\,\cN_{p,k}^{\rm H}(f)
 \le C_p\frac rs\cB_{p,s}(f)\,.
\end{equation*}
Taking the reference scale $s=1$ one and using
\begin{equation*}
 \cB_{p,1}(f)\simeq
 \cN_{p,1}^{\rm H}(f)=\|\nabla_pf\|_p
\end{equation*}
gives admissibility of the ball exponent $1$ and the heat kernel exponent $1/2$. Theorem~\ref{thm:critical-bounds} gives the reverse bounds
\begin{equation*}
 \alpha_p^{\rm B}\ge1,\qquad
 \alpha_p^{\rm H}\ge\frac1{d_w}=\frac12\,.
\end{equation*}
\end{proof}

With $\alpha=1/2$ and $\delta=D/2$, Theorem~\ref{thm:abstract-GN} gives
\begin{equation*}
 \|f\|_a
 \le C\|f\|_b^{1-\theta}\|\nabla_pf\|_p^\theta,
 \qquad
 \frac1a
 =\theta\left(\frac1p-\frac1D\right)
  +\frac{1-\theta}{b}\,,
\end{equation*}
for $f\in\ell^p(\Gamma)$ and the parameters specified there. In particular, Corollary~\ref{cor:abstract-GN-specializations} yields
\begin{align*}
 \|f\|_{pD/(D-p)}
 &\le C_p\|\nabla_pf\|_p,
 &&1\le p<D\,,\\
 \|f\|_2^{2+4/D}
 &\le C\cE(f,f)\|f\|_1^{4/D},
 &&f\in\ell^1(\Gamma)\,.
\end{align*}
When $D=1$, the endpoint Sobolev inequality is $\|f\|_\infty\le C\|\nabla_1f\|_1$ for $f\in\ell^1(\Gamma)$. Applying the $p=1$ inequality to $f=\1_A$ and using $\|\nabla_1\1_A\|_1=2\Per(A)$ gives
\begin{equation*}
 |A|^{(D-1)/D}\le C\Per(A),
 \qquad 0<|A|<\infty\,.
\end{equation*}
These are the classical group Sobolev and isoperimetric exponents; see~\cite{CoulhonSaloffCoste1993,HebischSaloffCoste1993}.

These conclusions apply, for instance, to $\mathbb Z^d$ with $D=d$ and to the discrete Heisenberg group $\mathbb H_3(\mathbb Z)$ with $D=4$ by the Bass--Guivarc'h formula~\cite{Bass1972,Guivarch1973}; see also~\cite[Section~2]{Blachere2003}.

\subsection{Strongly recurrent graphs}\label{sec:examples-strongly-recurrent}

Assume $({\rm SG}_{d_f,d_w})$ with $1\le d_f<d_w$, and put $\kappa=d_w-d_f$. Theorem~\ref{thm:wBE-strong-recurrence} and the identity $\cN_{p,1}^{\rm H}(f)=\|\nabla_pf\|_p$ give
\begin{equation}\label{eq:prescribed-bridge}
 \|f\|_{p,\beta_p(\kappa)}\simeq\|\nabla_pf\|_p,\qquad
 p\beta_p(\kappa)=\frac{d_f+(p-1)\kappa}{d_w},
 \qquad 1\le p\le2.
\end{equation}
For the heat kernel decay exponent $\delta=d_f/d_w$, one has
\begin{equation*}
 p\beta_p(\kappa)-\delta=\frac{(p-1)\kappa}{d_w}.
\end{equation*}
Thus the $L^1$ inequality is critical, while the inequalities for $1<p\le2$ lie in the supercritical range of Corollary~\ref{cor:abstract-GN-specializations}.

\begin{corollary}\label{cor:strongly-recurrent-inequalities}
Assume~\eqref{eq:LB} and $({\rm SG}_{d_f,d_w})$ with $1\le d_f<d_w$, and put $\kappa=d_w-d_f$. For $1\le p\le2$ and $f\in\ell^p(V,m)$,
\begin{equation}\label{eq:strongly-recurrent-infinity-GN}
 \|f\|_\infty
 \le C_p\|\nabla_pf\|_p^{\frac{d_f}{d_f+(p-1)\kappa}}
 \|f\|_p^{\frac{(p-1)\kappa}{d_f+(p-1)\kappa}}.
\end{equation}
At $p=1$, this reads
\begin{equation}\label{eq:strongly-recurrent-L1-Sobolev}
 \|f\|_\infty\le C\|\nabla_1f\|_1,
 \qquad f\in\ell^1(V,m).
\end{equation}
For $1<p\le2$ and $f\in\ell^1(V,m)$,
\begin{equation}\label{eq:strongly-recurrent-Nash}
 \|f\|_p\le C_p\|\nabla_pf\|_p^{\eta_p}\|f\|_1^{1-\eta_p},
 \qquad
 \eta_p=\frac{(p-1)d_f}{pd_f+(p-1)\kappa}.
\end{equation}
If $\operatorname{supp}f\subset\Omega$ and $m(\Omega)<\infty$, then
\begin{equation}\label{eq:strongly-recurrent-supported}
 \|f\|_p\le C_p m(\Omega)^{\frac{d_f+(p-1)\kappa}{pd_f}}
 \|\nabla_pf\|_p,\qquad 1\le p\le2.
\end{equation}
\end{corollary}

\begin{proof}
By Proposition~\ref{prop:volume-growth}, the masses at the vertices are bounded below. Apply Corollary~\ref{cor:abstract-GN-specializations}, parts~{\rm(ii)} and~{\rm(iv)}, with
\begin{equation*}
 \alpha=\beta_p(\kappa),\qquad \delta=\frac{d_f}{d_w},
 \qquad \cN_{p,1}^{\rm H}(f)=\|\nabla_pf\|_p.
\end{equation*}
The interpolation parameters are
\begin{equation*}
 \frac{\delta}{p\alpha}=\frac{d_f}{d_f+(p-1)\kappa},
 \qquad
 \frac{(p-1)\delta}{p(\alpha+\delta)-\delta}
 =\frac{(p-1)d_f}{pd_f+(p-1)\kappa}.
\end{equation*}
These give~\eqref{eq:strongly-recurrent-infinity-GN} and~\eqref{eq:strongly-recurrent-Nash}. For the supported estimate, insert
\begin{equation*}
 \|f\|_p\le m(\Omega)^{1/p}\|f\|_\infty
\end{equation*}
into~\eqref{eq:strongly-recurrent-infinity-GN} and rearrange.
\end{proof}

At $p=2$, the Nash and supported inequalities become
\begin{equation}\label{eq:strongly-recurrent-L2-Nash}
 \|f\|_2^{2+2d_w/d_f}\le C\cE(f,f)\|f\|_1^{2d_w/d_f}
\end{equation}
and
\begin{equation}\label{eq:strongly-recurrent-Faber-Krahn}
 \|f\|_2^2\le C m(\Omega)^{d_w/d_f}\cE(f,f),
 \qquad \operatorname{supp}f\subset\Omega.
\end{equation}
These are the usual Nash and Faber--Krahn inequalities for sub-Gaussian heat kernel decay. The $L^\infty$ interpolation estimate has a continuous Dirichlet-space counterpart in~\cite[Corollary~4.3 and Remark~4.4]{AlonsoRuizBaudoin2021}.

Under~\eqref{eq:p0} and $(\mathrm V_{d_f})$, the supported estimate of~\cite[Theorem~1.2]{Chen2019} has exponent $1/p+(p-1)/(pd_f)$. Its difference from the exponent in~\eqref{eq:strongly-recurrent-supported} is
\begin{equation*}
 \frac{p-1}{pd_f}(1-\kappa).
\end{equation*}
Hence~\eqref{eq:strongly-recurrent-supported} improves that estimate when $1<p\le2$ and $\kappa<1$; compare Remark~\ref{R:comparison_ball_vs_path}.

\subsection{Fractal graphs}\label{sec:examples-fractals}

The examples in this section concern the semigroup $P=(I+P_0)/2$, where $P_0$ corresponds to that associated with simple random walk. Both have the same reversible measure and the same harmonic functions, while $P$ satisfies~\eqref{eq:p0} and~\eqref{eq:LB}. Mean exit times from balls for $P$ are twice those for $P_0$. The characterization $({\rm VD})+({\rm H})+({\rm E}_\beta)\Leftrightarrow({\rm UE}_\beta)+({\rm LE}_\beta)$ of~\cite[Theorem~3.1]{GrigoryanTelcs2002} thus transfers the estimates proved for $P_0$ in the mentioned references to $P$ with the same exponents. Therefore, each of the graphs treated below has controlled weights and satisfies $({\rm SG}_{d_f,d_w})$ with $1\le d_f<d_w$. 

\medskip

We begin with the Vicsek graph, for which the critical exponents can be computed for every finite $p$.

\begin{example}
The planar Vicsek graph has
\begin{equation}\label{eq:Vicsek-dimensions}
 d_f=\frac{\log5}{\log3},\qquad d_w=d_f+1,\qquad \kappa=1;
\end{equation}
see~\cite[Section~4, Lemma~4.2]{Barlow2004} and~\cite{Chen2019}. For every $1\le p<\infty$,
\begin{equation}\label{eq:critical-Vicsek}
 \alpha_p^{\rm H}=\frac{d_f+p-1}{pd_w},
 \qquad
 \alpha_p^{\rm B}=\frac{d_f+p-1}{p}.
\end{equation}
Both exponents are attained, and both families satisfy weak monotonicity at their critical exponents.

For $p\ge2$, these assertions follow from Corollary~\ref{cor:critical-dw-df-one}. For $1\le p<2$, take a cell of diameter $R=3^n$ and the function $f_R$ from~\cite[Section~4, proof of Proposition~1.3]{Chen2019}. It equals one at the centre and zero at the boundary corners, is linear on the paths joining the centre to those corners, and is constant on the branches attached to these paths. Extend it by zero outside the cell. On the central subcell it is bounded below by a positive constant, while
\begin{equation*}
 \#\{\{x,y\}\in E:f_R(x)\ne f_R(y)\}\le CR,
 \qquad
 \sup_{x\sim y}|f_R(x)-f_R(y)|\le C/R.
\end{equation*}
Since the vertex masses and nonzero edge conductances are uniformly comparable to one,
\begin{equation*}
 m(\operatorname{supp}f_R)\le CR^{d_f},\qquad
 \|f_R\|_p^p\ge cR^{d_f},\qquad
 \|\nabla_pf_R\|_p^p\le CR^{1-p}.
\end{equation*}
Proposition~\ref{prop:critical-Poincare} yields $\alpha_p^{\rm B}\ge(d_f+p-1)/p$. The reverse bound and critical weak monotonicity follow from Theorem~\ref{thm:BL-WM}; Theorem~\ref{thm:ball-heat-comparison} and Proposition~\ref{prop:critical-comparison} give the corresponding heat kernel assertions.

For $1\le p\le2$, the supported inequality~\eqref{eq:strongly-recurrent-supported} coincides with~\cite[Theorem~1.2]{Chen2019}, which is optimal on Vicsek graphs by~\cite[Proposition~1.3]{Chen2019}. For $1<p<\infty$, the exponent $\alpha_p^{\rm H}$ in~\eqref{eq:critical-Vicsek} also coincides with the critical exponent obtained in~\cite{Feneuil2026} for the inequalities $\|\nabla_pf\|_p\le C\|L^\gamma f\|_p$ on Vicsek graphs.
\end{example}

\begin{example}
For the infinite Sierpi\'nski gasket graph, the heat kernel estimates of~\cite{Jones1996} give
\begin{equation}\label{eq:gasket-dimensions}
 d_f=\frac{\log3}{\log2},\qquad
 d_w=\frac{\log5}{\log2},\qquad
 \kappa=\frac{\log(5/3)}{\log2}\in(0,1).
\end{equation}
Expanding gasket cells $A_n$ satisfy $m(A_n)\simeq3^n$ and $\sup_n\Per(A_n)<\infty$. By Proposition~\ref{prop:critical-perimeter}, the critical exponents at $p=1$ are attained and satisfy
\begin{equation}\label{eq:critical-gasket}
 \alpha_1^{\rm H}=\frac{d_f}{d_w},\qquad \alpha_1^{\rm B}=d_f\,.
\end{equation}
The values and attainment at $p=2$ follow from~\eqref{eq:critical-two-and-limit}. For the other values of $p$, apply~\eqref{eq:critical-strongly-recurrent-small-p} and~\eqref{eq:critical-strongly-recurrent-large-p}.

Weak monotonicity on continuous nested fractals and their blow-ups is studied in~\cite{GaoYuZhang2026}. On the continuum gasket, the critical exponents involve nonlinear $p$-energy renormalization~\cite[Section~6.2]{Baudoin2024KS}. 
\end{example}

\begin{example}
For the Sierpi\'nski carpet, set $F_0=[0,1]^2$ and
\begin{equation*}
 F_{n+1}=\bigcup_{i\in\{0,1,2\}^2\setminus\{(1,1)\}}
 \frac{F_n+i}{3},\qquad
 F_{\rm pre}=\bigcup_{n\ge0}3^nF_n.
\end{equation*}
The vertices of the graphical carpet are the centres of the unit squares with integer vertices contained in $F_{\rm pre}$; two vertices are adjacent if their squares share a side. The degree is at most four, so the walk $P=(I+P_0)/2$ has controlled weights.

The estimates in~\cite{BarlowBass1999,Barlow2017} give $({\rm SG}_{d_f,d_w})$, with
\begin{equation}\label{eq:carpet-dimensions}
 d_f=\frac{\log8}{\log3},\qquad
 d_w\approx2.097,\qquad
 \kappa=d_w-d_f\in(0,1).
\end{equation}
Here $d_w=\log(8\rho)/\log3$, where $\rho$ is the resistance scale factor of the carpet; its exact value is unknown, and the numerical estimate $\rho\approx1.2515$ reported in~\cite[Section~2]{BarlowBass1992} gives the displayed value; see also~\cite{GaoYang2026}. For the critical exponents, apply~\eqref{eq:critical-two-and-limit}, \eqref{eq:critical-strongly-recurrent-small-p}, and~\eqref{eq:critical-strongly-recurrent-large-p}. As on the gasket, the continuum nonlinear energies of~\cite{MuruganShimizu2023} require a separate comparison with the graph energies before they can identify the remaining exponents.
\end{example}

\subsection{Products of strongly recurrent graphs}\label{sec:products}

We use the  product $G^{\otimes N}$ from Section~\ref{sec:critical-exponents}, with measure $m_N=m^{\otimes N}$, transition kernel $\pi_N$ and Markov operator $P_N=P^{\otimes N}$. Its graph distance is the maximum of the coordinate distances. The heat kernel and regularity estimates tensorize as follows.

\begin{proposition}\label{prop:tensor-product}
Let $G$ have controlled weights and satisfy $({\rm SG}_{d_f,d_w})$ and $({\rm wBE}(\kappa))$, where $0<\kappa<d_w$. For every integer $N\ge2$, the product $G^{\otimes N}$ has controlled weights and satisfies $(\mathrm V_{Nd_f})$, $({\rm SG}_{Nd_f,d_w})$ and $({\rm wBE}(\kappa))$.
\end{proposition}

\begin{proof}
Reversibility follows from $c_N=m_N\pi_N$. For adjacent distinct vertices,
\begin{equation}\label{eq:p0_LB_for_product}
 \pi_N(\mathbf x,\mathbf y)\ge
 \min\{\varepsilon_0,\varepsilon_{\rm L}\}^N,
 \qquad
 \pi_N(\mathbf x,\mathbf x)\ge\varepsilon_{\rm L}^N.
\end{equation}
Equation~\eqref{eq:product-distance-balls} and Proposition~\ref{prop:volume-growth} give
\begin{equation*}
 m_N(B_N(\mathbf x,r))=\prod_{i=1}^N\Vol(x_i,r)\simeq r^{Nd_f}.
\end{equation*}
Moreover,
\begin{equation}\label{eq:product-heat-kernel}
 P_N^n=(P^n)^{\otimes N},\qquad
 p_n^{(N)}(\mathbf x,\mathbf y)=\prod_{i=1}^Np_n(x_i,y_i).
\end{equation}
Taking the product of the upper estimates gives
\begin{equation}\label{eq:SG-upper_product}
 p_n^{(N)}(\mathbf x,\mathbf y)
 \le Cn^{-Nd_f/d_w}
 \exp\left[-c\left(
 \frac{d_N(\mathbf x,\mathbf y)}{n^{1/d_w}}
 \right)^{d_w/(d_w-1)}\right].
\end{equation}
For the lower estimate, condition~\eqref{eq:LB} gives
\begin{equation*}
 p_{n+1}(x,y)\ge
 \frac{\varepsilon_{\rm L}}{1+\varepsilon_{\rm L}}
 \bigl(p_n(x,y)+p_{n+1}(x,y)\bigr).
\end{equation*}
When $d_N(\mathbf x,\mathbf y)\le n$, it follows that
\begin{align*}
 p_n^{(N)}(\mathbf x,\mathbf y)+p_{n+1}^{(N)}(\mathbf x,\mathbf y)
 &\ge c\prod_{i=1}^N
 \bigl(p_n(x_i,y_i)+p_{n+1}(x_i,y_i)\bigr)\\
 &\ge cn^{-Nd_f/d_w}
 \exp\left[-C_N\left(
 \frac{d_N(\mathbf x,\mathbf y)}{n^{1/d_w}}
 \right)^{d_w/(d_w-1)}\right].
\end{align*}
This proves $({\rm SG}_{Nd_f,d_w})$.

Finally, change one coordinate at a time in $P_N^nF(\mathbf x)-P_N^nF(\mathbf y)$, as in~\cite[Proposition~3.8]{AlonsoEtAl2021III}. After integration in the other coordinates, each difference is bounded by the factor estimate~\eqref{eq:wBE}. Therefore
\begin{align*}
 |P_N^nF(\mathbf x)-P_N^nF(\mathbf y)|
 &\le C\sum_{i=1}^N
 \min\left\{1,\left(\frac{d(x_i,y_i)}{n^{1/d_w}}\right)^\kappa\right\}
 \|F\|_\infty\\
 &\le C_N\min\left\{1,
 \left(\frac{d_N(\mathbf x,\mathbf y)}{n^{1/d_w}}\right)^\kappa\right\}
 \|F\|_\infty.
\end{align*}
\end{proof}

Proposition~\ref{prop:critical-tensor} shows that the critical exponents are unchanged by these tensor powers. Thus a product of planar Vicsek graphs has the exponents in~\eqref{eq:critical-Vicsek}, where $d_f,d_w$ are the dimensions of one factor. The functional inequalities, however, use the product volume exponent $Nd_f$.

For a strongly recurrent factor, \eqref{eq:critical-two-and-limit} and the upper bound in~\eqref{eq:critical-heat-bounds} give
\begin{equation}\label{eq:product-dimension-range}
 \frac{d_w}{2}=\alpha_2^{\rm B}\le\frac{d_f+1}{2}\,.
\end{equation}
In particular, $0<\kappa\le1$. On the product,
\begin{equation*}
 \delta=\frac{Nd_f}{d_w},\qquad
 \delta-p\beta_p(\kappa)
 =\frac{(N-1)d_f-(p-1)\kappa}{d_w}.
\end{equation*}
Since $d_f\ge1$ and $\kappa\le1$, this difference is positive for $N\ge2$ and $1\le p\le2$ except when $(N,p,d_f,d_w)=(2,2,1,2)$.

\begin{corollary}\label{cor:strongly-recurrent-products}
Let $G$ have controlled weights and satisfy $({\rm SG}_{d_f,d_w})$ with $1\le d_f<d_w$. Put $\kappa=d_w-d_f$ and let $N\ge2$.

For $1\le p\le2$ with $(N,p,d_f,d_w)\ne(2,2,1,2)$,
\begin{equation}\label{eq:strongly-recurrent-product-Sobolev}
 \|F\|_{q_{p,N}}\le C_{p,N}\|\nabla_p^{(N)}F\|_p,
 \qquad
 q_{p,N}=\frac{pNd_f}{(N-1)d_f-(p-1)\kappa},
 \quad F\in\ell^p(V^N,m_N).
\end{equation}
In the exceptional case, for every $2<q<\infty$,
\begin{equation}\label{eq:critical-product-GN}
 \|F\|_q\le C_q\|\nabla_2^{(2)}F\|_2^{1-2/q}\|F\|_2^{2/q},
 \qquad F\in\ell^2(V^2,m_2).
\end{equation}
If, in this case, $\operatorname{supp}F\subset\Omega$ and $0<m_2(\Omega)<\infty$, then
\begin{equation}\label{eq:critical-product-logarithmic}
 \|F\|_\infty
 \le C\|\nabla_2^{(2)}F\|_2
 \left[\log\left(
 \frac{e\,m_2(\Omega)}{\inf_{\mathbf x\in V^2}m_2(\mathbf x)}
 \right)\right]^{1/2}.
\end{equation}
For every $N\ge2$ and $A\subset V^N$ with $0<m_N(A)<\infty$,
\begin{equation}\label{eq:product-set-isoperimetry}
 m_N(A)^{(N-1)/N}
 \le C_N\operatorname{Per}_{G^{\otimes N}}(A).
\end{equation}
For $1<p\le2$ and $F\in\ell^1(V^N,m_N)$,
\begin{equation}\label{eq:product-p-Nash}
 \|F\|_p\le C_{p,N}
 \|\nabla_p^{(N)}F\|_p^{\eta_{p,N}}\|F\|_1^{1-\eta_{p,N}},
 \qquad
 \eta_{p,N}=
 \frac{(p-1)Nd_f}{(p-1)Nd_f+d_f+(p-1)\kappa}.
\end{equation}
\end{corollary}

\begin{proof}
Theorem~\ref{thm:wBE-strong-recurrence}, Proposition~\ref{prop:tensor-product} and Theorem~\ref{thm:HK-WM} give $\mathrm P(p,\beta_p(\kappa))$ on the product. Apply Corollary~\ref{cor:abstract-GN-specializations} with $\delta=Nd_f/d_w$ and
\begin{equation*}
 \cN_{p,1}^{{\rm H},G^{\otimes N}}(F)=\|\nabla_p^{(N)}F\|_p.
\end{equation*}
Parts~{\rm(i)} and~{\rm(iii)} give~\eqref{eq:strongly-recurrent-product-Sobolev} and~\eqref{eq:critical-product-GN}; Corollary~\ref{cor:abstract-Trudinger} gives~\eqref{eq:critical-product-logarithmic}. At $p=1$, apply~\eqref{eq:strongly-recurrent-product-Sobolev} to $\1_A$ and use
\begin{equation*}
 \|\nabla_1^{(N)}\1_A\|_1
 =2\operatorname{Per}_{G^{\otimes N}}(A)
\end{equation*}
to obtain~\eqref{eq:product-set-isoperimetry}. Part~{\rm(ii)} gives~\eqref{eq:product-p-Nash}, since
\begin{equation*}
 \frac{(p-1)\delta}{p(\beta_p(\kappa)+\delta)-\delta}
 =\frac{(p-1)Nd_f}{(p-1)Nd_f+d_f+(p-1)\kappa}.
\end{equation*}
\end{proof}

At $p=2$, the Nash parameter is $Nd_f/(Nd_f+d_w)$; outside the exceptional case, the Sobolev exponent is $2Nd_f/(Nd_f-d_w)$. Both agree with the heat kernel decay exponent $Nd_f/d_w$. At $p=1$, the Sobolev and perimeter exponents depend only on $N$. They agree with the continuum nested-fractal product estimates in~\cite[Remark~4.19]{AlonsoEtAl2021III}; see also the tensor-product discussion in~\cite[Remark~{\rm(3)} after Theorem~4.4, equation~(4.18)]{Kaimanovich1992}, where the $p=1$ case is attributed to~\cite{Varopoulos1985}.

The $L^p$ inequalities contain more information than the $L^1$ estimate alone. The power-function argument of~\cite[Theorem~4.2]{Kaimanovich1992}, applied to $|F|^{p(N-1)/(N-p)}$, gives the exponent $pN/(N-p)$ when $1<p<N$. In contrast,
\begin{equation}\label{eq:product-Sobolev-comparison}
 \frac1{q_{p,N}}-\frac{N-p}{pN}
 =\frac{(p-1)(d_f-\kappa)}{pNd_f}\ge0,
 \qquad 1<p\le2,\quad p<N.
\end{equation}
If $\kappa<d_f$, the exponent $q_{p,N}$ is strictly smaller and the resulting global discrete Sobolev inequality is stronger, since the vertex masses are bounded below. When $N=p=2$ and $\kappa<d_f$, the estimate still has the finite exponent
\begin{equation*}
 q_{2,2}=\frac{4d_f}{d_f-\kappa}.
\end{equation*}

\subsection*{Acknowledgments}
F.B.\ is partially supported by grant 10.46540/4283-00175B from the Independent Research Fund Denmark, by the Villum Investigator grant \emph{Stochastic Analysis in Aarhus}, and by the European Research Council (ERC) under the European Union's Horizon Europe research and innovation programme (RanGe project, Grant Agreement No.~101199772).  L.C.\ is partially supported by grant 10.46540/4283-00175B from the Independent Research Fund Denmark.

\subsection*{AI disclosure statement}

Generative AI tools (ChatGPT and Claude) were used at a late stage of the project for proofreading and mathematical exploration, when the manuscript was already nearly final. This exploration led to the following mathematical improvements:
\begin{itemize}
\item Theorem~\ref{thm:WM-wBE} was originally proved by the authors for $p=1$ and $p=2$. ChatGPT 5.6 proposed a strategy for the full range $1\le p\le2$. The proposed proof was simplified and rewritten by the authors.
\item A weaker form of Theorem~\ref{thm:HK-WM-p-ge-2} was originally proved using a discrete analogue of~\cite[Theorem~5.1]{AlonsoEtAl2020I} and  duality. ChatGPT 5.6  found a proof of the present  statement that the authors simplified  using  the moment estimate~\eqref{eq:moment_estimate_Naor} from~\cite[Theorem~4.5]{NaorPeresSchrammSheffield2006}.
\end{itemize}
 The tools were also used for proofreading, improving proof presentation, bibliographic research, language editing, and ensuring consistency in the \LaTeX{} source. The authors take full responsibility for the content of the manuscript.

\end{document}